\documentclass[11pt,twoside,reqno]{amsart}
\usepackage{amsmath,amsfonts}
\usepackage{amssymb}
\usepackage{amsthm}
\usepackage[all,cmtip,line]{xy}
\usepackage{array}
\usepackage{psfrag}
\usepackage{overpic}
\newtheorem{theorem}{Theorem}%[section]
\newtheorem{lemma}{Lemma}[section]
\newtheorem{proposition}[lemma]{Proposition}%[section]
\newtheorem{corollary}[lemma]{Corollary}%[section]
\newtheorem{definition}{Definition}%[section]
\newtheorem{remark}[lemma]{Remark}%[section]
\newtheorem{problem}{{\bf{Problem}}}%[section]
\theoremstyle{definition}

\newcommand \alp{\alpha}
\newcommand \eps{\varepsilon}
\newcommand \vphi{\varphi}

\newcommand \Gam{\Gamma}

\newcommand \gam{\gamma}

\newcommand \R{\mathbb{R}}

\newcommand \til{\tilde}
\newcommand \wtil{\widetilde}

\newcommand \grad{D}

\newcommand \der{\partial}

\newcommand \mcl{\mathcal}

\newcommand \ol{\overline}
\newcommand \vtheta{\vartheta}

\newcommand \gradi{\nabla}
\newcommand{\divg}{ \mbox{div}}

\newcommand \Om{\Omega}

\newcommand \corners{\der S_0\cup\der\Gam_{ex}}
\newcommand \localcorners{\ol{\Sigma_{S_0,R}}\cap\ol{\Sigma_{\Gam_w,R}}}

\numberwithin{equation}{section}

\begin{document}

 %\tableofcontents
\title[Transonic shocks in divergent nozzles]
{Transonic shocks in multidimensional divergent nozzles}
\author{Myoungjean Bae}
\address{M. Bae, Department of Mathematics\\
         Northwestern University\\
         Evanston, IL 60208, USA}
\email{bae@math.northwestern.edu}

\author{Mikhail Feldman}
\address{M. Feldman, Department of Mathematics\\
         University of Wisconsin\\
         Madison, WI 53706-1388, USA}
\email{feldman@math.wisc.edu}

\date{}
\begin{abstract}
\iffalse
We consider transonic shocks in divergent nozzles(in $\R^n, n\ge 2$) provided the exit pressure is given appropriately. Many experimental results or a quasi-linear approximation show that if the exit pressure is fixed in a certain range, then a transonic shock appears in a divergent nozzle, and its location is uniquely determined by the exit flow parameters and the exist pressure.
But yet, there is no rigorous proof for the case of multidimensional divergent nozzles with arbitrary cross-sections.
In this paper
\fi
We establish existence, uniqueness and stability of transonic shocks for %inviscid
 steady compressible non-isentropic potential flow system in a
multidimensional divergent nozzle with an arbitrary smooth cross-section, for a prescribed exit pressure.
The proof is based on
solving a free boundary problem for a system of partial differential equations consisting of an elliptic equation and a
transport equation.
In the process, we obtain unique solvability for a class of transport equations with velocity fields of weak
regularity(non-Lipschitz), an infinite dimensional weak implicit mapping theorem which does not require
continuous Fr\'{e}chet differentiability, and regularity theory for a class of elliptic partial differential
equations with discontinuous oblique boundary conditions.
\end{abstract}
%\vfill\eject
\maketitle

\section{Introduction}
We consider inviscid compressible steady flow of ideal polytropic gas.
\iffalse:
\begin{align*}
&div(\rho \vec u)=0\\
&div(\rho \vec u\otimes\vec u+pI_n)=0\\
&div(\rho \vec u(\frac 12|\vec u|^2+\frac{\gam p}{(\gam-1)\rho}))=0
\end{align*}
where $\rho, \vec u, p$ represent the density, velocity and the pressure respectively. Here, the constant $\gam>1$ is the \emph{adiabatic exponent} of the flow.
\fi

One of important but difficult subjects in the study of transonic flow is to understand a global feature of flow through a convergent-divergent nozzle so called a \emph{de Laval nozzle}. According to the quasi-linear approximation in \cite[Chapter V. section 147]{CF}, given incoming subsonic flow at the entrance, the flow is always accelerated through the convergent part of the nozzle unless the exit pressure exceeds the pressure at the entrance. The flow pattern through the divergent part, however, varies depending on the exit pressure and the shape of the nozzle.  In the divergent part, the flow may remain subsonic all the way to the exit, or it may be accelerated to a supersonic state after passing the throat of the nozzle and have a transonic shock across which the velocity jumps down from supersonic to subsonic. In particular, the approximation implies that if the exit pressure is given by a constant $p_c$ satisfying $0<p_{min}<p_c<p_{max}<\infty$ for some constants $p_{min}$ and $p_{max}$,  then a transonic shock occurs in the divergent part of the nozzle.  One may refer \cite{CF} for the quasi-linear approximation in curved nozzles. In section \ref{subsection:radial}, for any given constant exit pressure $p_c\in(p_{min}, p_{max})$, we rigorously compute the transonic flow of the Euler system in a multidimensional nozzle expanding cone-shaped.

The issues above motivated recent works by several authors on existence and stability of steady transonic shocks in nozzles using the models of potential flow or compressible Euler system
\cite{ChenFeldman1, Ch-F2, Ch-F3,zh1, Sxchen05, Yu-06, Jch-Ch-F, zh2}. In these works transonic shocks were studied in cylindrical nozzles
and its perturbations. A physically natural setup for transonic shock problem is to prescribe the
parameters of flow on the
nozzle entrance
and the pressure on the nozzle exit. In particular, an important issue is to show that the shock
location is uniquely determined by these parameters. However, such problem is not well posed in the cylindrical
nozzles. Indeed,
the flat shock between uniform states can be translated along the cylindrical  nozzle,
and this transform does not change the
flow parameters on the nozzle entrance and the pressure on the nozzle exit.
This provides an explicit example of non-uniqueness.
Also, for uniform states in cylindrical nozzles, the flow on the nozzle entrance determines the pressure on the nozzle exit for
a transonic shock solution. This degeneracy leads to non-existence, i.e. one cannot prescribe an arbitrary
(even if small) perturbation of pressure on the
nozzle exit. On the other hand, as the de Laval nozzle example suggests, it is natural to study transonic
shocks in the diverging nozzles. This removes the translation invariance in the case of Euler system. Recent works
by Liu, Yuan\cite{Liu-Y},  S.-X. Chen \cite{Schen3}, and by Li, Xin, Yin \cite{LiXinYin09}
showed well-posedness of transonic shock problem with prescribed flow on entrance and pressure on exit of a divergent
nozzle, which are small perturbation of the corresponding parameters of the background solution,
 in the case of compressible Euler system in dimension two, under some additional restrictions
on the perturbed pressure on the exit. In particular, it was assumed that for the prescribed pressure on exit,
the normal derivatives on the nozzle walls vanish. While this condition is not physically natural,
mathematically
it leads to substantial simplifications: it allows to show $C^2$ regularity of the shock up to the nozzle walls,
which allows to work with $C^1$ vector fields in the transport equations.

In this paper we study transonic shocks in diverging nozzles in any dimension, for the case
of perturbed cone-shaped nozzle of arbitrary cross-section. Moreover, we do not make assumptions
on the pressure on the nozzle exit, other than the smallness of the perturbation. This requires to consider
shocks with regularity deteriorating near the nozzle walls, and then vector fields in transport equations
are of low regularity (non-Lipschitz).

Furthermore, we study this problem in the framework of potential flow. A surprising feature of this problem
is that it is not well-posed for the ``standard'' potential flow equation. Namely, for radial solutions
in a cone-shaped nozzle, the spherical transonic shock can be placed in any location between the nozzle entrance and exit,
without changing the parameters of flow on the nozzle entrance and exit. This follows from
a calculation by G.-Q Chen, M. Feldman \cite[pp. 488-489]{ChenFeldman1}, see more details below.
Thus, spherical shocks separating radial flows in the cone-shaped divergent nozzles can be translated
along the nozzle for potential flow equation.
This is similar
to the case of flat shocks separating uniform flows in straight nozzles.

In order to remove this degeneracy, we introduce a new potential flow model, which we call {\em non-isentropic
potential flow system}. It is based on the full Euler system (as opposed to the ``standard'' potential
flow equation, which is based on the isentropic Euler system), and allows entropy jumps across the shock.
The system consists of an equation of second order, which is elliptic in the subsonic regions and hyperbolic
in the supersonic regions, and
two transport equations.

The main purpose of this paper is to establish the existence, uniqueness and stability of a transonic shock
for non-isentropic
potential flow system
in a multidimensional divergent nozzle with an arbitrary smooth cross-section provided that we smoothly perturb the
cone-shaped divergent nozzle, incoming radial supersonic flow and the constant  pressure at the exit.

In Section \ref{Section:nozzle1}, we introduce the non-isentropic potential flow model, and define a transonic shock
solution of this new model. Then we state two main theorems of this paper and present a framework of the proof for the theorems. Section \ref{sec:free-bdry} is devoted to a free boundary problem for the velocity potential $\vphi$ where the exit normal velocity is fixed . In section \ref{sec:nozzle3}, we plug a solution $\vphi$ of the free boundary problem in section \ref{sec:free-bdry} to a transport equation to find the corresponding pressure $p$. In section \ref{sec:nozzle4}, we prove the existence of a transonic shock solution of the non-isentropic potential flow for a fixed exit pressure $p_{ex}$. The uniqueness of the transonic shock solution is proven in section \ref{sec:nozzle5}.
\section{Main theorems}
\label{Section:nozzle1}
\subsection{Compressible Euler system and  shocks}
We study transonic shocks of steady inviscid compressible flow of ideal polytropic gas. Such flow is governed by the conservation of mass, momentum and the conservation of energy. These three conservation laws provide a system of PDEs for the density $\rho$, velocity $\vec u=(u_1,\cdots,u_n)$ and the pressure $p$ as follows:
\begin{equation}
\label{2-1}
\begin{split}
&div(\rho \vec u)=0\\
&div(\rho \vec u \otimes \vec u+pI_n)=0\\
&div\bigl(\rho \vec u(\frac 12|\vec u|^2+\frac{\gam p}{(\gam-1)\rho})\bigr)=0
\end{split}
\end{equation}
where $I_n$ is $n\times n$ identity matrix.
The quantity $c$, given by
\begin{equation}
\label{soundspeed}
c:=\sqrt{\frac{\gam p}{\rho}},
\end{equation}
is called the \emph{sound speed}. The flow type can be classified by the ratio, \emph{called the Mach number}, of $|\vec u|$ to $c$. If $\frac{|\vec u|}{c}>1$ then the flow is said \emph{supersonic}, if $\frac{|\vec u|}{c}<1$ then the flow is said \emph{subsonic}. If $\frac{|\vec u|}{c}=1$ then the flow is said \emph{sonic}. Depending on entrance and exit data of flow or a shape of a channel, there may be a jump transition of the flow type across a curve or a surface. Such a transition can be understood in a weak sense as follows.

For  an open set $\Om\subset \R^n$, if $(\rho,\vec u,p)\in L^1_{loc}(\Om)$ satisfies
\begin{equation}
\label{1}
\begin{split}
\int_{\Om}\rho \vec u\cdot\nabla \xi dx
=\int_{\Om}(\rho u_k\vec u+p\hat{e}_k)\cdot\nabla\xi dx
=\int_{\Om}\rho\vec u\bigl(\frac 12|\vec u|^2+\frac{\gam p}{(\gam-1)}\bigr)\cdot\nabla \xi dx=0
\end{split}
\end{equation}
for any $\xi\in C_0^{\infty}(\Om)$ and all $k=1,\cdots,n$ where $\hat{e}_k$ is a unit normal in k-th direction for $\R^n$, $(\rho,\vec u,p)$ is called a \emph{weak solution} to (\ref{2-1}) in $\Om$.

Suppose that a smooth (n-1)-dimensional surface $S$ divides $\Om$ into two disjoint subsets $\Om^{\pm}$. If
 $(\rho,\vec u,p)\in L^1_{loc}(\Om)$ is in $C^{1}(\Om^{\pm})\cap C^0(\ol{\Om^{\pm}})$ then $(\rho,\vec u,p)$ is a weak solution of \eqref{2-1} if and only if $(\rho,\vec u,p)$
satisfies \eqref{2-1} in each of $\Om^{\pm}$ and the
\emph{Rankine-Hugoniot jump conditions}(which is abbreviated as \emph{R-H conditions} hereafter)
\begin{align}
\label{2-2}
&[\rho\vec u \cdot \nu_s]_S=0\\
\label{2-2-1}
&[\rho(\vec u\cdot \nu_s)\vec u+p\nu_s]_S=\vec 0\\
\label{2-2-2}
&[\frac 12|\vec u|^2+\frac{\gam p}{(\gam-1)\rho}]_S=0
\end{align}
for a unit normal $\nu_s$ on $S$ where $[F]_S$ is defined by
$$[F(x)]_S:=F(x)|_{\overline{\Om^-}}-F(x)|_{\overline{\Om^+}}\quad\text{for}\;\;x\in S.$$

To study a transonic shock in $\R^n$ for $n\ge 2$, we consider an irrotational flow in which the velocity $\vec u$ has an expression of $\vec u=\nabla\vphi$ for a scalar function $\vphi$. Such a flow is called a \emph{potential flow} and $\vphi$ is called a \emph{(velocity) potential}.

\subsection{Steady isentropic potential flow equation and radial transonic shocks}
\label{isentrPotentFl}
A widely used steady potential flow model(see e.g. \cite{ChenFeldman1, zh1} et al.) consists of the conservation law of
mass, the
Bernoulli law for the velocity. So it can be written into the following
second-order,
nonlinear elliptic-hyperbolic equation of mixed type
for the velocity potential $\varphi: \Omega\subset\R^n\rightarrow\R^1$:
\begin{equation}\label{PotenEulerCompres}
\divg(\rho(|\nabla\varphi|^2)\nabla\varphi)=0,
\end{equation}
where, with a particular choice of the Bernoulli's invariant, the density function $\rho(q^2)$ is
\begin{equation}\label{densFunc}
\rho(q^2)=\left(1-\frac{\gamma-1}{2} q^2\right)^{\frac{1}{\gamma-1}}.
\end{equation}

We note that equation (\ref{PotenEulerCompres}) can be derived formally from the isentropic Euler system.
Steady isentropic Euler system is obtained from (\ref{2-1})
by setting $p=\kappa\rho^\gamma$, where $\kappa>0$ is a constant, and dropping the third line in (\ref{2-1}).
Then equation  (\ref{PotenEulerCompres}) can be obtained by
a formal calculation,  substituting $\vec u=\nabla\vphi$ into the  isentropic Euler system,
and  choosing $\kappa=\frac{1}{\gam}$.
 \begin{equation}\label{isentropic-ellipticity}
 \begin{split}
&\text{\emph{The second-order nonlinear equation \eqref{PotenEulerCompres} is}}\\
&\text{\emph{ strictly elliptic if}}\;\;
|\nabla\vphi|>c,\;\;
\text{\emph{and is strictly hyperbolic if}}\;\;
|\nabla\vphi|<c
\end{split}
\end{equation}
with $c=\sqrt{\frac{\gam p}{\rho}}=\sqrt{\rho^{\gam-1}}$.
The elliptic regions of the equation \eqref{PotenEulerCompres}
correspond to the subsonic flow, and the hyperbolic regions of
 \eqref{PotenEulerCompres} to the supersonic flow.

From the divergence structure of (\ref{PotenEulerCompres}), it is easy to derive the following RH conditions
on shocks. Suppose that a smooth (n-1)-dimensional surface $S$ divides $\Om$ into two disjoint subsets $\Om^{\pm}$. If
 $\varphi\in C^{0,1}(\Om)$ is in $C^{2}(\Om^{\pm})\cap C^1(\Om^{\pm}\cup S)$ then $\varphi$ is a
a weak solution of \eqref{PotenEulerCompres} if and only if $\varphi$
satisfies \eqref{PotenEulerCompres} in each of $\Om^{\pm}$ and
\begin{equation}\label{FBCondition-0}
\bigg[\rho(|\nabla\varphi|^2)\nabla\varphi\cdot \nu\bigg]_S=0.
\end{equation}
Note that $[\varphi]_S=0$ since $\varphi\in C^{0,1}(\Om)$.

Shock $S$ is called transonic if it separates subsonic and supersonic regions.

\medskip
Now we discuss radial transonic shocks, i.e. transonic shock solutions of the form
$$
\varphi(x)=w(r)\quad  \text{where}\quad r=|x|.
$$
Let $0<r_1<r_2$
$$
\Om=\{r_1<|x|<r_2\}.
$$
Then calculation in \cite[pp. 488-489]{ChenFeldman1} show that there exists (radial)
functions $\varphi_0^\pm(x)=w_0^\pm(|x|$ (defined by \cite[eqn (7.4)]{ChenFeldman1}
with e.g. $R_0:=r_1$) such that
\begin{enumerate}\renewcommand{\theenumi}{\roman{enumi}}
\item
$\varphi_0^\pm\in C^\infty(\Om)$ and satisfy (\ref{PotenEulerCompres}) in $\Om$,
\item
$\varphi_0^-$ is supersonic and $\varphi_0^+$ is subsonic in $\Om$.
Moreover $0<\max_{r\in[r_1, r_2]}\frac{d\,w_0^+}{d\,r}(r)<\min_{r\in[r_1, r_2]}\frac{d\,w_0^-}{d\,r}(r)$.
\item for any $r_s\in (r_1, r_2)$, choosing a constant $c=c(r_s)$ so that $w_0^+(r_s)=w_0^-(r_s)+c$,
there holds:
function
$$
\varphi_0(x)=\min(\varphi_0^+(x), \varphi_0^-(x)+c(r_s))\quad \text{for}\quad
x\in \Omega
$$
is a transonic shock solution of (\ref{PotenEulerCompres}) with shock $S=\{|x|=r_s\}$,
and $\varphi_0=\varphi_0^-$ in $\Om^-_{r_s}:=\{r_1<|x|<r_s\}$; and
 $\varphi_0=\varphi_0^+$ in $\Om^+_{r_s}:=\{r_s<|x|<r_2\}$. In particular
$\varphi_0$ is supersonic in $\Om^-_{r_s}$ and subsonic in $\Om^+_{r_s}$.
\end{enumerate}
Note that $\partial_r\varphi_0^\pm>0$. Thus flow enters $\Om$ through $\{|x|=r_1\}$
and exits $\Om$ through $\{|x|=r_2\}$.
By part(iii) above, $\nabla\varphi$ does not change on the entrance and exit when the position $r_s$
of the shock changes within $(r_1, r_2)$. Since $\vec{u}=\nabla\varphi$ determines the
density $\rho$ by (\ref{densFunc})
and pressure $p=\kappa\rho^\gamma$, it follows that parameters of the flow on the entrance and exist of the
annulus-shaped domain $\Om$ do not determine position $r=r_s$ of the shock.

Restricting the annulus-shaped domain $\Om$ to a cone-shaped (nozzle) domain $\mcl{N}$
as in Section \ref{subsection:radial}, we see that the the same example shows non-uniqueness of
the shock position in the cone-shaped nozzle.

Thus we propose a new potential model, to fix this degeneracy.

\subsection{The non-isentropic potential flow system and transonic shocks}

\iffalse
The main purpose of this paper is to accomplish the unique existence and stability of a transonic shock in a multidimensional divergent nozzle provided that the incoming supersonic flow and the exit pressure are given in an appropriate range. As we will show later, the \emph{isentropic potential flow}, however, cannot be used for our purpose because it causes an ill-posedness of the problem. So it is inevitable to introduce a new potential flow model which is not isentropic.
\fi

The new model, which we call \emph{the non-isentropic potential flow system},
is obtained by substituting $\vec{u}=\nabla\varphi$ into the the full steady compressible Euler system (\ref{2-1}).
Unknown functions in the new model are the scalar functions $(\varphi, \rho, p)$.
Note that second line of  (\ref{2-1}), representing conservation of momentum,
consists of $n$ equations. Thus we have overdeterminacy
in the case of potential model, i.e. we need to retain one equation
from the conservation of momentum equations. For that, we note that the component in the direction of $\vec{u}$
of the vector determined by the second line of  (\ref{2-1})
 represents conservation of entropy along streamlines (in the smooth flow regions).
This motivates us to introduce the following non-isentropic potential flow system:
\begin{align}
\label{2-3}
&div(\rho \gradi\vphi)=0,\\
\label{2-4}
&\gradi\vphi\cdot div(\rho\gradi\vphi\otimes \gradi\vphi+pI_n)=0,\\
\label{2-5}
&div\bigl(\rho\gradi\vphi(\frac 12|\gradi\vphi|^2+\frac{\gam p}{(\gam-1)\rho})\bigr)=0.
\end{align}
\eqref{2-3} and \eqref{2-5}  represent the conservation of mass and energy respectively, and
\eqref{2-4} concerns the conservation of the entropy along each streamline but it allows the entropy to change
between streamlines.

Since \eqref{2-4} is not in a divergence form, we cannot define a shock solution for this new model in sense of \eqref{1}. Instead, we employ all the R-H conditions \eqref{2-2}-\eqref{2-2-2} to define a shock solution of the new model as follows:
\begin{definition}[Shock solutions]
\label{definition2-1}
$(\rho,\vphi,p)$ is a shock solution of the non-isentropic potential flow in $\Om$ with a shock $S$ if
\begin{itemize}
\item[(i)]  $\rho,p\in C^0(\overline{\Om^{\pm}})\cap C^1(\Om^{\pm})$, $\vphi\in C^1(\overline{\Om^{\pm}})\cap C^2(\Om^{\pm})$;
%\item $|\gradi\vphi|>c$ in $\Om^-$ and $|\gradi\vphi|<c$ in $\Om^+$;
\item[(ii)] $(\rho,\nabla\vphi,p)|_{S\cap\ol{\Om^-}}-(\rho,\nabla\vphi,p)|_{S\cap\ol{\Om^+}}\neq \vec 0$;
\item[(iii)]  $(\rho,\vphi,p)$ satisfies {\rm{(\ref{2-3})-(\ref{2-5})}} pointwise in $\Om^{\pm}$, and {\rm{(\ref{2-2})-(\ref{2-2-2})}} on $S$ with $\vec u=\gradi\vphi$.
\end{itemize}
\end{definition}

For $(\rho,\varphi, p)\in C^{\infty}$, if $\rho>0$ then \eqref{2-3}-\eqref{2-5} are equivalent to
\begin{align}
\label{2-3'}
&div\bigl((B-\frac 12|\gradi{\vphi}|^2)^{\frac{1}{\gam-1}}\gradi{\vphi}\bigr)=0 \\
\label{2-4'}
&\gradi{\vphi}\cdot\gradi\frac{{p}}{(B-\frac 12|\gradi{\vphi}|^2)^{\frac{\gam}{\gam-1}}}=0 \\
\label{2-5'}
&\nabla\vphi\cdot\nabla B=0
\end{align}
where the \emph{Bernoulli's invariant} $B$ is defined by
\begin{equation}
\label{bernoulli}
B=\frac 12|\nabla\vphi|^2+\frac{\gam p}{(\gam-1)\rho}.
\end{equation}
An explicit computation shows that, as a second-order PDE for $\vphi$, \eqref{2-3'} with \eqref{bernoulli} is
\begin{equation*}
\text{strictly elliptic if}\;\;|\nabla\varphi|<c, \;\; \text{and is strictly hyperbolic if}\;\;|\nabla\varphi|>c
\end{equation*}
with $c$ defined in \eqref{soundspeed}. Note that this is similar to \eqref{isentropic-ellipticity}.

\begin{definition}[A transonic shock solution]
\label{def-ts}
A shock solution $(\rho,\vphi,p)$, defined in Definition \ref{definition2-1}, is a transonic shock solution if (i)$\der_{\nu_s}\vphi|_{\ol{\Om^{-}}\cap S}>\der_{\nu_s}\vphi|_{\ol{\Om^+}\cap S}>0$ for the unit normal $\nu_s$ on $S$ pointing toward $\Om^+$, and (ii) $|\nabla\vphi|>c$ in $\Omega^-$, $|\nabla\vphi|<c$ in $\Omega^+$ where $c$ is defined in \eqref{soundspeed}.
\end{definition}
%For $\Omega\subset\R^n$, \eqref{2-2-1} consists of $n$ equations.
Let $\tau_s$ be a tangent vector field on $S$, then, by \eqref{2-2-1}, there holds
\begin{equation}
\label{2-6}
[\rho(\gradi\vphi\cdot\nu_s)(\gradi\vphi\cdot\tau_s)]_S=0.
\end{equation}
Assuming no vacuum state in $\Om$ i.e., $\rho\neq 0$ in $\Om$, \eqref{2-2} and \eqref{2-6} imply either $\nabla\vphi\cdot\nu_s=0$ or $[\gradi\vphi\cdot\tau_s]_S=0$. If $\nabla\vphi\cdot\nu_s=0$ then Definition \ref{def-gs} (ii) implies  $[\gradi\vphi\cdot\tau_s]_S\neq 0$. In that case, $S$ becomes a vortex sheet. If $[\nabla\vphi\cdot\tau_s]_S=0$ holds on $S$ then this implies $[\vphi]_S=k$ for some constant $k$. We are interested in the case of $\nabla\vphi\cdot\nu_s\neq 0$ on $S$ so \eqref{2-6} provides
\begin{equation}
\label{2-7}
[\vphi]_S=0
\end{equation}
choosing the constant $k=0$ without loss of generality. In addition, $(\ref{2-2-1})\cdot\nu_s$ provides
\begin{equation}
\label{2-8}
[\rho(\gradi\vphi\cdot\nu_s)^2+p]_S=0.
\end{equation}
\begin{remark}
\label{remark-entropy}
The specific entropy $\epsilon$ of ideal polytropic gas satisfies the constitutive relation
$$
\kappa\exp(\epsilon/c_v)=\frac{p}{\rho^{\gamma}}
$$
for constants $\kappa$ and $c_v$. A straightforward calculation using \eqref{2-2},\eqref{2-2-2},\eqref{2-7} and \eqref{2-8} indicates that the entropy $\epsilon$ increases across a shock $S$ for the non-isentropic potential flow model. The entropy jump across a shock plays an essential role for the well-posedness of the problem considered in this paper.
\end{remark}

\subsection{Radial transonic shock solutions}
\label{subsection:radial}
%As the simplest transonic shock solution in a straight divergent nozzle, we first consider a radial transonic shock solution.

Fix an open connected set $\Lambda\subset \mathbb{S}^{n-1}(n\ge 2)$ with a smooth boundary $\der\Lambda$, and define a straight divergent nozzle $\mcl{N}$ by
\begin{equation*}
\mcl{N}:=\{x: r_0<|x|<r_1, \frac{x}{|x|}\in \Lambda\}\quad\text{for}\;\;0<r_0<r_1<\infty.
\end{equation*}

%Since our interest is in transonic flow through $\mcl{N}$, it may be convenient to use a spherical coordinate system from time to time.
Regarding $\Lambda$ as a n-1 dimensional submanifold of $\mathbb{R}^{n-1}$, let $x'=(x_1',\!\cdots,\!x_{n-1}')$ be a coordinate system in $\Lambda$. Particularly, if $\Lambda\subsetneq\mathbb{S}^{n-1}$, then $\Lambda$ can be described by a single smooth diffeomorphism(e.g. stereographic projection), say, $T$ i.e., there is $\Lambda_T\subset\R^{n-1}$ with a smooth boundary so that $\Lambda_T$ can be described by
\begin{equation}
\label{lambda-t}
\Lambda=\{T(x'):x'\in\Lambda_T\}.
\end{equation}
Then, for any $x\in\mcl{N}$, there exists unique $x'$ satisfying $\frac{x}{|x|}=T(x')$. For convenience, we denote as
\begin{equation}
\label{spherical-coord}
x=_T(r,x')\quad x'\in_T\Lambda\quad
\text{if}\;\;r=|x|,\quad T(x')=\frac{x}{|x|}\in\Lambda.
\end{equation}
To simplify notations, we denote as $x'\in\Lambda$ instead of $x'\in_T\Lambda$ hereafter.
The coordinates $(r,x')$ in \eqref{spherical-coord} is regarded as a spherical coordinate system in $\mcl{N}$.

For later use, we set
\begin{equation*}
\begin{split}
&\Gam_{ent}:=\{x=_T(r,x'):r=r_0, x'\in\Lambda\},\quad\Gam_{ex}:=\{x=_T(r,x'):r=r_1, x'\in\Lambda\},\\
&\Gam_w:=\der\mcl{N}\setminus(\Gam_{ent}\cup\Gam_{ex}).
\end{split}
\end{equation*}

To simplify notations, for a function $\psi(x)$ defined for $x\in\mcl{N}$, if $x=_T(r,x')$ then we identify $\psi(r,x')$ with $\psi(x)$ without any further specification for the rest of paper. Similarly, any function $f(\frac{x}{|x|})$ defined for $\frac{x}{|x|}\in \Lambda$, we identify $f(x')$ with $f(\frac{x}{|x|})$ if $\frac{x}{|x|}=T(x')$. Furthermore, any function $\vphi_0$ varying only by the radial variable $|x|$, we write as $\vphi_0(r)$. Also, in $(r,x')$ coordinates, we denote as
\begin{equation}\label{15}
\mcl{N}=(r_0,r_1)\times \Lambda_T,\;\;\mcl{N}^-_{a}=\mcl{N}\cap\{r<a\},\;\;\mcl{N}^+_a=\mcl{N}\cap \{r>a\}.
\end{equation}

For a fixed constant vector $(\rho_{in},v_{in},p_{in})\in\R^3_+$ with $v_{in}>c_{in}(=\sqrt{\frac{\gam p_{in}}{\rho_{in}}})$,
let us find a transonic shock solution $(\rho_0,\vphi_0,p_0)$ satisfying
\begin{align}
\label{2-3-1}
&(\rho_0,\der_r\vphi_0,p_0)=(\rho_{in},v_{in},p_{in})\quad\text{on}\;\Gam_{ent},\\
\label{2-3-2}
&\der_{\nu_w}\vphi_0=0\quad\text{on}\;\Gam_w
\end{align}
where $\nu_w$ is the unit normal on $\Gam_w$ toward the interior of $\mcl{N}$.
Since (\ref{2-3})-(\ref{2-5}) are invariant under rotations, we expect $(\rho_0,\vphi_0, p_0)$ to be functions of $r$ only. Then (\ref{2-3-2}) automatically holds on $\Gam_w$. For a continuously differentiable radial solution $(\rho,\vphi,p)$, if $\vphi_r\neq 0$ in $\mcl{N}$,  then (\ref{2-3})-(\ref{2-5}) are equivalent to
\begin{equation}
\label{A}
\begin{split}
&\frac{d}{dr}(r^{n-1}\rho\vphi_r)=0,\\
&\rho\frac{d}{dr}(\frac 12\vphi_r^2)+\frac{dp}{dr}=0,\\
&\frac 12\vphi_r^2+\frac{\gam p}{(\gam-1)\rho}=\frac 12 v_{in}^2+\frac{\gam p_{in}}{(\gam-1)\rho_{in}}=:B_0.
\end{split}
\end{equation}
\eqref{A} can be considered as a system of ODEs for $(\vphi_r,\rho)$ because of
\begin{equation}
\label{2-3-3}
\rho=\frac{\gam p}{(\gam-1)(B_0-\frac 12\vphi_r^2)}.
\end{equation}
To find $(\vphi_r,p)$, we rewrite (\ref{A}) as
\begin{align}
\label{2-3-4}
&\frac{d\vphi_r}{dr}=\frac{2(n-1)(\gam-1)\vphi_r(B_0-\frac 12\vphi_r^2)}{(\gam+1)r(\vphi_r^2-K_0)},\\
\label{2-3-5}
&\frac{dp}{dr}=-\frac{2(n-1)\gam\vphi_r^2 p}{(\gam+1)r(\vphi_r^2-K_0)}
\end{align}
with
\begin{equation}\label{k0}K_0:=\frac{2(\gam-1)B_0}{\gam+1}.\end{equation}

%To find a radial transonic shock solution with the shock at $r=r_s$ for $r_s\in(r_0,r_1)$,
Fix $r_s\in(r_0,r_1)$, and let us find a radial transonic shock solution with the shock on $\{r=r_s\}$. For that purpose, we need to solve (\ref{2-3-4}) and (\ref{2-3-5}) in two separate regions of supersonic state and subsonic state.

For the supersonic region $\mcl{N}^-_{r_s}$, the initial condition is
\begin{equation}
\label{2-d2-1}
(\vphi_r,p)(r_0)=(v_{in},p_{in}).
\end{equation}
Then (\ref{2-3-4}), \eqref{2-3-5} with (\ref{2-d2-1}) have a unique solution $(\der_r\vphi_0^-,p_0^-)$ in $\mcl{N}^-_{r_s}$. For the subsonic region $\mcl{N}^+_{r_s}$, we solve the system of the algebraic equations (\ref{2-2}), (\ref{2-8}) and (\ref{2-3-3}) with $\der_r\vphi_0^-$ and $p_0^-$ on the side of $\mcl{N}^-_{r_s}$ for $(\der_r\vphi,p)(r_s)$. Then there are two solutions for $(\der_r\vphi,p)(r_s)$ but the only admissible solution in sense of Definition \ref{def-ts} is
\begin{align}
\label{2-2-5}
&\der_r\vphi(r_s)=\frac{K_0}{\der_r\vphi_0^-(r_s)},\\
\label{2-2-6}
&p(r_s)=(\rho_0^-(\der_r\vphi_0^-)^2+p_0^--\rho_0^-K_0)(r_s)=:p_{s,0}(r_s).
\end{align}
So the initial condition for the subsonic region is given by (\ref{2-2-5}) and (\ref{2-2-6}).

Let $(\der_r\vphi_0^+,p_0^+)$ be the solution to (\ref{2-3-4}), (\ref{2-3-5}) with (\ref{2-2-5}), (\ref{2-2-6}).

We claim that $(\der_r\vphi_0^-,p_0^-)$ is indeed supersonic in $\mcl{N}^-_{r_s}$ while $(\der_r\vphi_0^+,p_0^+)$ is subsonic in $\mcl{N}^+_{r_s}$.
Plug
$\frac{\gam p}{\rho}=c^2$ into the equation of Bernoulli's law $\frac 12(\der_r\vphi)^2+\frac{\gam p}{(\gam-1)\rho}=B_0$ to obtain
\begin{equation}
\label{2-3-6}
\vphi_r^2-c^2=\frac{\gam+1}{2}(\vphi_r^2-K_0).
\end{equation}
Denote $\sqrt{\frac{\gam p_0{^{\pm}}}{\rho_0^{\pm}}}$ as $c_{\pm}$ then, $v_{in}^2>c_{in}^2$ and (\ref{2-3-6}) imply $\frac{d\vphi_r}{dr}> 0$ for $\vphi_r=\der_r\vphi_0^-$ thus $(\der_r\vphi_0^-)^2>c_-^2$ for $r\ge r_0$ so $(\der_r\vphi_0^-,p_0^-)$ is supersonic in $\mcl{N}^-_{r_s}$. On the other hand, by (\ref{2-2-5}), $\der_r\vphi_0^+(r_s)$ satisfies the inequality
\begin{equation}
\label{2-2-8}
(\der_r\vphi_0^+)^2-c_+^2=\frac{\gam+1}{2}\frac{K_0}{(\der_r\vphi_0^-)^2}\bigl(K_0-(\der_r\vphi_0^-)^2\bigr)<0
\end{equation}
at $r=r_s$, and this implies $\frac{d\vphi_r}{dr}<0$ for $\vphi_r=\der_r\vphi_0^+$. So $(\der_r\vphi_0^+)^2<c_+^2$ holds for $r\ge r_s$.
Therefore, we get a family of radial transonic shock solutions as follows:
\begin{definition}[\textbf{Background solutions}]\label{background}
For $r_s\in(r_0,r_1)$, define
\begin{equation*}
(\rho_0,\vphi_0,p_0)(r;r_s):=
\begin{cases}
(\rho_0^-,\vphi_0^-,p_0^-)(r;r_s)&\text{in}\; [r_0,r_s]\\
(\rho_0^+,\vphi_0^+,p_0^+)(r;r_s)&\text{in}\;(r_s,r_1]
\end{cases}
\end{equation*}
with $\rho_0^{\pm}(r;r_s):=\frac{\gam p_0^{\pm}(r)}{(\gam-1)(B_0-\frac 12(\der_r\vphi_0^{\pm}(r))^2)}$, and
\begin{equation*}
\vphi_0^-(r;r_s):=\int_{r_0}^{r}\der_r\vphi_0^-(t)dt,\quad \vphi_0^+(r;r_s)=\int_{r_s}^r\der_r\vphi_0^+(t)dt+\vphi^-_0(r_s)
\end{equation*}
where $(\der_r\vphi_0^{\pm},p_0^{\pm})(r)$ is the solution to {\rm{(\ref{2-3-4})-(\ref{2-3-5})}} with {\rm{(\ref{2-2-5})-(\ref{2-2-6})}} and \eqref{2-d2-1} respectively.
We call $(\rho_0,\vphi_0,p_0)(r;r_s)$ the background solution with the transonic shock on $\{r=r_s\}$.
\end{definition}
\begin{remark}
\label{remark2-2-1}
For a small constant $\delta>0$,
$(\der_r\vphi_0^{\pm},p_0^{\pm})$ can be extended to $\mcl{N}^+_{r_s-2\delta}$ and $\mcl{N}^-_{r_s+2\delta}$ as solutions to \eqref{2-3-4}-\eqref{2-3-5} with {\rm{(\ref{2-2-5})-(\ref{2-2-6})}} and {\rm{(\ref{2-d2-1})}} respectively where $\mcl{N}^{\pm}_{r_s\mp \delta}$ is defined by \eqref{15}.
\end{remark}
\begin{remark}
\label{remark2-d2-1}
For a fixed $r_s$, $\rho_0^-$ and $p_0^-$ monotonically decrease while $\der_r\vphi_0^-$ monotonically increases by $r$. On the other hand, by {\rm{(\ref{2-2-8})}},  $\rho_0^+$ and $p_0^+$ monotonically increase while $\der_r\vphi_0^+$ monotonically decreases.
\end{remark}

Another important property of the background solutions is the monotonicity of the exit values of $(\rho_0,\der_r\vphi_0,p_0)$ depending on the location of shocks $r_s$. In \cite{Yu}, it is proven that, for the full Euler system, $\rho_0^+(r_1;r_s)$ and $p_0^+(r_1;r_s)$ monotonically decrease with respect to $r_s$ while $\der_r\vphi_0^+(r_1;r_s)$ monotonically increases. If $\der_r\vphi_0(r;r_s)\neq 0$ in $\mcl{N}$, then the non-isentropic potential flow model is equivalent to the full Euler system for smooth flow. So the monotonicity properties of the exit data applies to the background solution $(\rho_0,\vphi_0,p_0)$ of the non-isentropic potential flow model so we have:
\begin{proposition}[{\cite[Theorem 2.1]{Yu}}]
\label{proposition2-2-1}
Let us set $p_{min}=p_0^+(r_1;r_1)$ and $p_{max}:=p_0^+(r_1;r_0)$. Then $p_{min}<p_{max}$ holds. Moreover, for any constant $p_c\in(p_{min},p_{max})$, there exists unique $r^*_s\in(r_0,r_1)$ so that the background solution $(\rho_0,\vphi_0,p_0)(r;r^*_s)$ satisfies $p_0^+(r_1;r_s^*)=p_c$.
\end{proposition}
Readers can refer \cite{Yu} for the proof although we would like to point out that the following lemma, which will be used later in this paper,  plays an important role for the proof of Proposition \ref{proposition2-2-1}.
\begin{lemma}
\label{lemma-1}
For any given $r_s\in(r_0,r_1)$, let $(\rho_0^+,\vphi_0^+,p_0^+)(r)$ be the background transonic shock solution with the shock on $\{r=r_s\}$ defined in Definition \ref{background}, and let us define $\mu_0$ by
$$
\mu_0:=\frac{\frac{d}{dr}(\frac{K_0}{\der_r\vphi_0^-}-\der_r\vphi_0^+)(r_s)}{\der_r(\vphi_0^--\vphi_0^+)(r_s)}.
$$

\begin{itemize}
\item[(a)] Then, there holds
 $\mu_0>0$.
\item[(b)]Also, we have, for $p_{s,0}$ defined in \eqref{2-2-6}
\begin{equation}
\label{2-3-8}
\bigl[\frac{dp_{s,0}}{dr}-\frac{dp_0^+}{dr}\bigr](r_s)=-\frac{(n-1)\gam p_0^-[(\der_r\vphi_0^-)^2+
\frac{\gam-1}{\gam+1}K_0]}{(\gam-1)r(B_0-\frac 12(\der_r\vphi_0^-)^2)}|_{r=r_s}<0.
\end{equation}
\end{itemize}
\begin{proof}
By  \eqref{2-3-4}, we have
\begin{equation*}
\begin{split}
&\frac{d}{dr}\bigl(\frac{K_0}{\der_r\vphi_0^-}-\der_r\vphi_0^+\bigr)(r_s)\\
&=
-\frac{C(n,\gamma)}{r_s}\bigl[\frac{K_0(B_0-\frac 12(\der_r\vphi_0^-)^2)}{\der_r\vphi_0^-((\der_r\vphi_0^-)^2-K_0)}
+\frac{\der_r\vphi_0^+(B_0-\frac 12 (\der_r\vphi_0^+)^2)}{(\der_r\vphi_0^+)^2-K_0}\bigr](r_s).
\end{split}
\end{equation*}
with the constant $C(n,\gamma)=\frac{2(n-1)(\gam-1)}{\gam+1}>0$. Applying \eqref{2-2-5} to replace $\der_r\vphi_0^+(r_s)$ by $\frac{K_0}{\der_r\vphi_0^-(r_s)}$, we obtain
\begin{equation}
\label{3-1-21}
\frac{d}{dr}\bigl(\frac{K_0}{\der_r\vphi_0^-}-\der_r\vphi_0^+\bigr)(r_s)
=\frac{2(n-1)\gam K_0}{(\gam+1)r_s\der_r\vphi_0^-(r_s)}>0.
\end{equation}
By \eqref{2-2-5} and \eqref{2-2-8}, it is easy to see $\der_r(\vphi_0^--\vphi_0^+)(r_s)>0$, and thus (a) holds true for all $r_s\in(r_0,r_1)$.

Similarly, by (\ref{2-3-3})-(\ref{2-3-5}), (\ref{2-2-5}) and (\ref{2-2-6}), we obtain \eqref{2-3-8}.
\end{proof}
\end{lemma}

\iffalse

\begin{remark}
\label{2-2-3}
In the isentropic potential flow, the explicit calculation found in \cite[pp488--489]{ChenFeldman1} shows
\begin{equation*}
p_0^+(r_1;r_s)=p_0^+(r_1;r_0)
\end{equation*}
for all $r_s\in[r_0,r_1]$. That is, the exit pressure of the radial transonic shock solution is determined only by $(\rho_{in},v_{in},p_{in}),r_0$ and $r_1$ regardless of $r_s$.
\end{remark}
According to Remark {\ref{2-2-3}}, the transonic shock problem  with a fixed exit pressure for the isentropic potential flow model is ill-posed.
Intuitively, the ill-posedness above indicates that a jump of the entropy across a shock is essential for well-posedness of the transonic shock problem with a fixed exit pressure in a divergent nozzle.
\fi
\subsection{Main theorems}
According to Proposition \ref{proposition2-2-1}, for any given constant $p_c\in(p_{min},p_{max})$, there exists a transonic shock solution whose exit pressure is $p_c$. Our goal is to achieve the existence of a transonic shock solution when we smoothly perturb incoming supersonic flow, the exit pressure and the nozzle $\mcl{N}$.

Let $\wtil{\mcl{N}}$ be a nozzle smoothly perturbed from $\mcl{N}$ by a diffeomorphism $\Psi:\mcl{N}\to \R^n$, and $(\til{\rho},\til{\vphi},\til p)$ be a transonic shock solution of a non-isentropic potential flow in $\wtil{\mcl{N}}$ with a shock $\wtil S$. Then $\wtil S$ separates $\wtil{\mcl{N}}$ into a supersonic region $\wtil{\mcl{N}}^-$ and a subsonic region $\wtil{\mcl{N}}^+$. Let us denote $(\til{\rho},\til{\vphi},\til p)|_{\wtil{\mcl{N}}^{\pm}}$ as $(\til{\rho}_{\pm},\til{\vphi}_{\pm},\til p_{\pm})$. If the entrance data for $(\til{\rho},\til{\vphi},\til p)$ is given so that the Bernoulli's invariant $B$ is a constant at the entrance of the nozzle $\wtil{\mcl{N}}$ then, by \eqref{2-2}, (\ref{2-2-2}) and \eqref{2-5'}, $(\til{\rho}_{\pm},\til{\vphi}_{\pm},\til p_{\pm})$ is a solution to the system of (\ref{2-3'}), (\ref{2-4'}) and
\begin{equation}
\label{2-3-11}
\frac 12|\gradi{\vphi}|^2+\frac{\gam {p}}{(\gam-1){\rho}}=B_0\quad\text{in}\quad \wtil{\mcl{N}}^{\pm}
\end{equation}
where the constant $B_0>0$ is determined by the entrance data.

On the shock $\wtil S$, (\ref{2-7}) and (\ref{2-3-11}) imply
\begin{equation}
\label{2-3-13}
\bigl[\frac 12(\gradi\til{\vphi}\cdot\til{\nu}_s)^2+\frac{\gam \til{p}}{(\gam-1)\til{\rho}}\big]_{\wtil{S}}=0
\end{equation}
for a unit normal $\til{\nu}_s$ on $\wtil S$.

Let us set $\wtil{K}_s:=\frac{2(\gam-1)}{\gam+1}\bigl(\frac 12(\gradi\til{\vphi}_-\cdot\til{\nu}_s)^2+\frac{\gam \til{p}_-}{(\gam-1)\til{\rho}_-}\bigr)$, and solve (\ref{2-2}), (\ref{2-8}) and (\ref{2-3-13}) for $\gradi\til{\vphi}_+\cdot\til{\nu}_s$ and $\til{p}_+$ to obtain
\begin{align}
\label{2-3-15}
&\gradi\til{\vphi}_+\cdot\til{\nu}_s=\frac{\wtil{K}_s}{\gradi\til{\vphi}_-\cdot\til{\nu}_s},\\
\label{2-3-16}
&\til{p}_+=\til{\rho}_-(\gradi\til{\vphi}_-\cdot\til{\nu}_s)^2+\til{p}_--\til{\rho}_-\wtil{K}_s\quad\text{on}\quad\wtil{S}.
\end{align}

For a transonic shock solution $(\til{\rho},\til{\vphi},\til p)$ in $\wtil{\mcl{N}}$, let
$\wtil{\mcl{N}}^-:=\{|\nabla\til\vphi|>c\}$ be the supersonic region, and $\wtil{\mcl{N}}^+:=\{|\nabla\til\vphi|<c\}$ be the subsonic region where $c$ is defined in \eqref{soundspeed} by replacing $\rho, p$ by $\til \rho, \til p$ on the right-hand side.
\begin{problem}
[Transonic shock problem in a perturbed nozzle $\wtil{\mcl{N}}=\Psi(\mcl{N})$]
\label{problem1}
Let $(\til{\rho}_-,\til{\vphi}_-,\til p_-)$ be a supersonic solution upstream with $B=B_0$  at the entrance $\wtil\Gam_{ent}=\Psi(\Gam_{ent})$ for a constant $B_0$,  and suppose that $\nabla\til{\vphi}_-\cdot\nu_w=0$ on $\wtil\Gam_w=\Psi(\Gam_w)$ for a unit normal $\nu_w$ on $\wtil\Gam_w$.

Given an exit pressure function $\til{p}_{ex}$ on $\wtil\Gam_{ex}=\Psi(\Gam_{ex})$, locate a transonic shock $\wtil S$ i.e., find a function $f$ satisfying
$$
\wtil{\mcl{N}}^-=\wtil{\mcl{N}}\cap\{|x|<f(\frac{x}{|x|})\},\quad
\wtil{\mcl{N}}^+=\wtil{\mcl{N}}\cap\{|x|>f(\frac{x}{|x|})\},$$
and then find a corresponding subsonic solution $(\til{\rho}_+,\til{\vphi}_+,\til p_+)$ downstream so that
\begin{itemize}
\item[(i)] $(\til{\rho}_+,\til{\vphi}_+,\til p_+)$ satisfies {\rm{(\ref{2-3'}), (\ref{2-4'})}} and \eqref{2-3-11} in $\wtil{\mcl{N}}^+$ ;
\item[(ii)] $\til{\vphi}_+$ satisfies the slip boundary condition $\nabla\til{\vphi}_+\cdot\nu_w=0$ on $\wtil \Gam_w^+=\der\wtil{\mcl{N}}^+\cap \wtil \Gam_w$;
\item[(iii)] $(\til{\vphi}_+,\til p_+)$ satisfies {\rm{(\ref{2-7}), (\ref{2-3-15})}} and {\rm{(\ref{2-3-16})}} on $\wtil S$;
\item[(iv)] $\til p_+=\til p_{ex}$ holds on $\wtil \Gam_{ex}$.
\end{itemize}
\end{problem}

For convenience, we reformulate Problem \ref{problem1} for $(\rho,\vphi,p)=(\til\rho,\til\vphi,\til p)\circ \Psi$ in the unperturbed nozzle $\mcl{N}$. If $(\til\rho,\til\vphi,\til p)$ solves \eqref{2-3'}, \eqref{2-4'}, \eqref{2-3-11} then $(\rho,\vphi,p)$ satisfies
\begin{align}
\label{3-0-1}
&div A(x,D\Psi,\grad\vphi)=0\\
\label{2-p}
&(D\Psi^{-1})^T(D\Psi^{-1})\grad \vphi\cdot \grad \frac{p}{(B_0-\frac 12|D\Psi^{-1}\grad\vphi|^2)^{\frac{\gam}{\gam-1}}}=0\\
\label{2-b}
&\frac 12|D\Psi^{-1}\grad\vphi|^2+\frac{\gam p}{(\gam-1)\rho}=B_0
\end{align}
for $x\in\Psi^{-1}(\wtil{\mcl{N}}^{\pm})=:\mcl{N}^{\pm}$ with $D\vphi=(\der_{x_j}\vphi)_{j=1}^n$, $D\Psi=(\der_{x_i}\Psi_j)_{i,j=1}^n$, $D\Psi^{-1}=(D\Psi)^{-1}$. Here, $A(x,D\Psi,D\vphi)$ in \eqref{3-0-1} is defined by
\begin{equation}
\label{3-0-2}
A(x,m,\eta)=\det m\; (B_0-\frac 12|m^{-1}\eta|^2)^{\frac{1}{\gam-1}}(m^{-1})^Tm^{-1}\eta
\end{equation}
for $(x,m,\eta)\in \mcl{N}\times B^{(n\times n)}_{d}(I_n)\times B^{(n)}_R(0)$ where we set $B^{(k)}_{\mcl{R}}(a)$ as a closed ball of the radius $\mcl{R}$ in $\R^k$ with the center at $a\in\R^k$. Here, $d$ is a sufficiently small so that \eqref{3-0-2} is well defined for all $m\in B_d^{n\times n}(I_n)$.

By \eqref{2-7}, the unit normal on $\wtil{S}$ toward $\wtil{\mcl{N}}^+$ has the expression of
\begin{equation}
\label{3-0-3}
\til{\nu}_s=\frac{D_y(\til{\vphi}_--\til{\vphi}_+)}{|D_y(\til{\vphi}_--\til{\vphi}_+)|}\;\;\text{with}\;\;D_y=(\der_{y_1},\cdots,\der_{y_n}).
\end{equation}
So the R-H conditions for $(\rho,\vphi,p)$ on $S:=\Psi^{-1}(\wtil{S})$ are
\begin{align}
\label{3-0-4}
&\vphi_+=\vphi_-\\
\label{3-0-5}
&Q^*(D\Psi,D\vphi_-,D\vphi_+)D\vphi_+\cdot\nu_s
=\frac{K_s}{Q^*(D\Psi,D\vphi_-,D\vphi_+)\grad\vphi_-\cdot\nu_s}\\
\label{2-p-bc}
&p=\rho_-[(D\Psi^{-1})^T(D\Psi^{-1})\grad\vphi_-\cdot\nu_s]^2+p_--\rho_-K_s
\end{align}
with $D=(\der_{x_1},\cdots,\der_{x_n})$ where we set
\begin{align}
\label{3-0-20}
&Q^*(m,\xi,\eta):=\frac{|\xi-\eta|}{|m^{-1}(\xi-\eta)|}(m^{-1})^T(m^{-1}),\\
\label{3-0-21}
&K_s=\frac{2(\gam-1)}{\gam+1}\bigl(\frac 12(\frac{|\grad(\vphi_--\vphi)|(D\Psi^{-1})^TD\Psi^{-1}\grad\vphi_-\cdot\nu_s}{|D\Psi^{-1}\grad(\vphi_--\vphi)|})^2+\frac{\gam p_-}{(\gam-1)\rho_-}\bigr),\\
\label{normal}
&\nu_s=\frac{D(\vphi_--\vphi_+)}{|D(\vphi_--\vphi_+)|}.
\end{align}
We note that, by \eqref{3-0-4}, $\nu_s$ is the unit normal on $S$ toward $\mcl{N}^+$. Problem \ref{problem1} is equivalent to:
\begin{problem}
[Transonic shock problem in the unperturbed nozzle $\mcl{N}$]
\label{problem2}
Let $\Psi$ be a smooth diffeomorphism in $\mcl{N}$. Given a supersonic solution $(\rho_-,\vphi_-,p_-)=(\til{\rho}_-,\til{\vphi}_-,\til p_-)\circ\Psi$ in $\mcl{N}^-$ where $(\til{\rho}_-,\til{\vphi}_-,\til p_-)$ is as in Problem \ref{problem1}, and an exit pressure function $p_{ex}=\til p_{ex}\circ \Psi$, locate a transonic shock $S$ i.e., find a function $f$ satisfying $\mcl{N}^-=\{|D\Psi^{-1}D\vphi|>c\}=\mcl{N}\cap\{r<f(x')\}$and $\mcl{N}^+=\{|D\Psi^{-1}D\vphi|<c\}=\mcl{N}\cap \{r>f(x')\}$ for the $(r,x')$ coordinates defined in \eqref{spherical-coord}, and find a corresponding subsonic solution $(\rho_+,\vphi_+,p_+)$ in $\mcl{N}^+$ so that
\begin{itemize}
\item[(i)] $(\rho_+,\vphi_+,p_+)$ satisfies {\rm{(\ref{3-0-1})-(\ref{2-b})}} in $\mcl{N}^+$;
\item[(ii)]On $\Gam_w^+:=\Gam_w\cap \der\mcl{N}^+,$ $\vphi_+$ satisfies
\begin{equation}
\label{3-0-6}
A(x,D\Psi,\grad\vphi)\cdot \nu_w=0
 \end{equation}
for a unit normal $\nu_w$ on $\Gam_w$;

\item[(iii)] $(\vphi_+,p_+)$ satisfies {\rm{(\ref{3-0-4})-(\ref{2-p-bc})}} on $S$;
\item[(iv)] $p=p_{ex}$ holds on $\Gam_{ex}$.
\end{itemize}
\end{problem}

Our claim is that if a background supersonic solution $(\rho_0^-,\vphi_0^-,p_0^-)$ and a constant exit pressure $p_c\in(p_{min},p_{max})$ are perturbed sufficiently small, and also if $\Psi$ is a small perturbation of the identity map, then there exists a corresponding transonic shock solution for a non-isentropic potential flow.

To study Problem \ref{problem2}, we will use weighted H\"{o}lder norms. For a bounded connected open set $\Om\subset\R^n$, let $\Gam$ be a closed portion of $\der\Om$. For $x,y\in\Om$, set
\begin{equation*}
\delta_x:=dist(x,\Gam)\quad \delta_{x,y}:=\min(\delta_x,\delta_y).
\end{equation*}
For $k\in\R$, $\alp\in(0,1)$ and $m\in \mathbb{Z}_+$, we define
\begin{align*}
&\|u\|_{m,0,\Om}^{(k,\Gam)}:=\sum_{0\le|\beta|\le m}\sup_{x\in \Om}\delta_x^{\max(|\beta|+k,0)}|D^{\beta}u(x)|\\
&[u]_{m,\alp,\Om}^{(k,\Gam)}:=\sum_{|\beta|=m}\sup_{x,y\in\Om, x\neq y}\delta_{x,y}^{\max(m+\alp+k,0)}\frac{|D^{\beta}u(x)-D^{\beta}u(y)|}{|x-y|^{\alp}}\\
&\|u\|_{m,\alp,\Om}^{(k,\Gam)}:=\|u\|_{m,0,\Om}^{(k,\Gam)}+[u]_{m,\alp,\Om}^{(k,\Gam)}
\end{align*}
where we write $D^{\beta}=\der_{x_1}^{\beta_1}\cdots\der_{x_n}^{\beta_n}$ for a multi-index $\beta=(\beta_1,\cdots,\beta_n)$ with $\beta_j\in\mathbb{Z}_+$ and $|\beta|=\sum_{j=1}^n\beta_j$. We use the notation $C^{m,\alp}_{(k,\Gam)}(\Om)$ for the set of functions whose $\|\cdot\|_{m,\alp,\Om}^{(k,\Gam)}$ norm is finite.

\begin{remark}
\label{remark-holdernorms}
Let $\Om$ be a connected subset of $\mcl{N}$, and $\Gamma$ be a portion of  $\der\Om$. For $T$ in \eqref{spherical-coord}, let us set
$
\Om_*:=\{(r,x'): x=_T(r,x'), x\in\Om\}$ and $\Gamma_*:=\{(r,x'):x=_T(r,x'), x\in\Gam\}.
$
Then there is a constant $C_*(\mcl{N},\Om,n,m,\alp,k)$ so that, for any $u\in C^{m,\alp}_{k,\Gamma}(\Om)$, there holds
\begin{equation}
\label{est-spherical-euclid}
\frac{1}{C_*}\|u\|_{m,\alp,\Om_*}^{(k,\Gamma_*)}\le\|u\|_{m,\alp,\Om}^{(k,\Gamma)}\le C_*\|u\|_{m,\alp,\Om_*}^{(k,\Gamma_*)}.
\end{equation}
\end{remark}

For a fixed constant data $(\rho_{in},v_{in},p_{in})\in\R^3_+$ satisfying $v_{in}^2>\frac{\gam p_{in}}{\rho_{in}}$ and $p_c\in(p_{min},p_{max})$, let $(\rho_0,\vphi_0,p_0)$ be the background solution with $(\rho_0,\der_r\vphi_0,p_0)(r_0)=(\rho_{in},v_{in},p_{in})$ and $p_0(r_1)=p_c$, and let $S_0=\{r=r_s\}\cap \mcl{N}$ be the shock of $(\rho_0,\vphi_0,p_0)$. Now, we state the main theorems.
\begin{theorem}
[\emph{Existence}]
\label{main-thm'}
For any given $\alp\in(0,1)$, there exist constants $\sigma_1,\delta, C>0$ depending on $(\rho_{in},v_{in},p_{in}),p_c,r_0,r_1,\gam,n,\Lambda$ and $\alp$ so that whenever $0<\sigma\le\sigma_1$, if
\begin{itemize}
\item[(i)] a diffeomorphism $\Psi:\ol{\mcl{N}}\to \R^n$ satisfies
$$
\varsigma_1:=|\Psi-Id|_{2,\alp,\mcl{N}}\le \sigma,
$$
\item[(ii)] $(\rho_-,\vphi_-,p_-)$ is a solution to {\rm{(\ref{3-0-1})-(\ref{2-b})}} in $\mcl{N}^-_{r_s+2\delta}$, and satisfies the boundary condition {\rm{(\ref{3-0-6})}} on $\Gam_w\cap\mcl{N}^-_{r_s+2\delta}$ and the estimate
\begin{equation*}
\varsigma_2:=|\rho_--\rho_0^-|_{2,\alp,\mcl{N}^-_{r_s+2\delta}}+|p_--p_0^-|_{2,\alp,\mcl{N}^-_{r_s+2\delta}}+|\vphi_--\vphi_0^-|_{3,\alp,\mcl{N}^-_{r_s+2\delta}}\le \sigma,
\end{equation*}
\item[(iii)] a function $p_{ex}$ defined on $\Gam_{ex}$ satisfies
$$
\varsigma_3:=\|p_{ex}-p_c\|_{1,\alp,\Lambda}^{(-\alp,\der\Lambda)}\le \sigma,
$$
\end{itemize}
then there exists a transonic shock solution $(\rho,\vphi,p)$ for a non-isentropic potential flow with a shock $S$ satisfying that
\begin{itemize}
\item[(a)] there is a function $f:\Lambda\to \R^+$ so that $S$, $\mcl{N}^-\!\!=\{x\in\mcl{N}:|D\Psi^{-1}D\vphi|>c\}$ and $\mcl{N}^+\!\!=\{x\in\mcl{N}:|D\Psi^{-1}D\vphi|<c\}$ are given by
    \begin{equation*}
    \begin{split}
    &S=\{(r,x')\in\mcl{N}:r=f(x'),x'\in\Lambda\},\\
    &\mcl{N}^-=\mcl{N}\cap\{(r,x'):r<f(x')\},\quad
    \mcl{N}^+=\mcl{N}\cap\{(r,x'):r>f(x')\},
    \end{split}
    \end{equation*}
    and $f$ satisfies
    \begin{equation*}
|f-r_s|_{1,\alp,\Lambda}\le C(\varsigma_1+\varsigma_2+\varsigma_3)\le C\sigma;
\end{equation*}

\item[(b)] $(\rho,\vphi,p)=(\rho_-,\vphi_-,p_-)$ holds in $\mcl{N}^-$;
\item[(c)] $(\rho,\vphi,p)|_{\mcl{N}^+}$ satisfies {\rm{(\ref{3-0-4})-{\rm{(\ref{2-p-bc})}}, (\ref{3-0-6})}} and the estimate
$$
\|\rho-\rho_0^+\|_{1,\alp,\mcl{N}^+}^{(-\alp,\Gam_w)}+\|\vphi-\vphi_0^+\|_{2,\alp,\mcl{N}^+}^{(-1-\alp,\Gam_w)}
+\|p-p_0^+\|_{1,\alp,\mcl{N}^+}^{(-\alp,\Gam_w)}\le C(\varsigma_1+\varsigma_2+\varsigma_3);
$$

\item[(d)] $p=p_{ex}$ holds on $\Gam_{ex}$.

\end{itemize}
\end{theorem}
\begin{remark}
\label{remark-cons}
In this paper, we phrase any constant depending on $\rho_{in},v_{in},p_{in},$ $p_c,r_0,r_1,\gam,n, \Lambda$  and $\alp$ as a constant depending on the data unless otherwise specified.
\end{remark}
\begin{theorem}
[Uniqueness]
\label{theorem-uniq}
Under the same assumptions with Theorem \ref{main-thm'}, for any  $\alp\in(\frac 12,1)$, there exists a constant $\sigma_2>0$ depending on the data in sense of Remark \ref{remark-cons} so that whenever $0<\sigma\le \sigma_2$, the transonic shock solution in Theorem \ref{main-thm'} is unique in the class of weak solutions satisfying Theorem \ref{main-thm'}(a)-(d).
\end{theorem}
\begin{remark}
By Theorem \ref{main-thm'}(c), Theorem \ref{theorem-uniq} implies that the transonic shock solutions are stable under a small and smooth perturbation of the incoming supersonic flow, exit pressure and the nozzle.
\end{remark}
\subsection{Framework of the proof of Theorem \ref{main-thm'} and Theorem \ref{theorem-uniq}}
\label{subsec-framework}
By the decomposition \eqref{3-0-1}-\eqref{2-b},  if we solve \eqref{3-0-1} for $\vphi$, then we can plug it into the transport equation to solve for $p$, and solve the Bernoulli's law for $\rho$. To solve \eqref{2-p} for $p$, however, the initial condition on $S$ can be computed by \eqref{3-0-5} and \eqref{2-p-bc}. But if we also fix $p=p_{ex}$ on $\Gam_{ex}$, then the problem gets overdetermined. So, we set the framework to prove Theorem \ref{main-thm'} and \ref{theorem-uniq} as follows.\\

(Step 1) To prove Theorem \ref{main-thm'}, we first solve
    \begin{problem}[Free boundary problem for ${\vphi}$]
     \label{problem2-1}
     Given the incoming supersonic flow $(\rho_-,\vphi_-,p_-)$ and a function ${v}_{ex}$ defined on $\Gam_{ex}$, find $f\in C^{1,\alp}_{(-\alp,\der\Lambda)}(\Lambda)$ and $\vphi$ in ${\mcl{N}}^+$ so that (i) $\vphi$ satisfies $|D\Psi^{-1}D\vphi|<c$ and \eqref{3-0-1} in ${\mcl{N}}^+=\mcl{N}\cap\{(r,x'):r>f(x')\}$, (ii) also satisfies  {\rm{(\ref{3-0-4}), (\ref{3-0-5}), (\ref{3-0-6})}}
    and
    $
    M(x,\grad\vphi)=v_{ex}\;\;\text{on}\;\;\Gam_{ex}
$
    with the definition of $M(x,\grad\vphi)$ being given later.
\end{problem}
Let us emphasize that Lemma \ref{lemma-1} is the key ingredient to prove Lemma \ref{lemma3-1-1}, and this lemma provides well-posedness of Problem \ref{problem2-1}.\\

(Step 2) In Section \ref{sec:nozzle3}, we solve \eqref{2-p} for ${p}$ after plugging $\vphi$ in to \eqref{2-p} by the method of characteristics in ${\mcl{N}}^+$ with the initial condition (\ref{2-p-bc}) on $S$. However, as we will show later, $D\vphi$ is only in $C^{\alp}$ up to the boundary of ${\mcl{N}}^+$ so the standard ODE theory is not sufficient to claim solvability of \eqref{2-p}. For that reason, we need a subtle analysis of $D\vphi$ to overcome this difficulty.\\

(Step 3) For a fixed $(\Psi,\vphi_-,p_-)$, define $\mcl{P}$ by
$\mcl{P}(\Psi,\vphi_-,p_-,v_{ex})=p_{ex}$
 where $p_{ex}$ is the exit value of the solution $p$ to \eqref{2-p} where $\vphi$ in Problem \ref{problem2-1} is uniquely determined by $(\Psi, \vphi_-, p_-, v_{ex})$. In Section \ref{sec:nozzle4}, we will show that $\mcl{P}(\Psi,\vphi_-,p_-,\cdot):v_{ex}\mapsto p_{ex} $ is locally invertible in a weak sense near the background solution. The weak invertibility is proven by modifying the proof of the \emph{right inverse mapping theorem}
for finite dimensional Banach spaces in \cite{Sz}, and this proves Theorem \ref{main-thm'}.\\

(Step 4) To prove Theorem \ref{theorem-uniq}, we need to estimate the difference of two transonic shock solutions in a weaker space than the space of the transonic shock solutions in step 3. This requires a careful analysis of solutions to elliptic boundary problems with a blow-up of oblique boundary conditions at the corners of the nozzle $\mcl{N}$. See Lemma \ref{prop7-1}.

%%%%%%%%%%%%%%%%%%%%%%%%%%%%%%%%%%%%%%%%%%%%%%%%%%%%%%%%%%%%%%%%%%%%%
%%%%%%%%%%%%%%%%%%%%%%%%%%%%%(Section 3)%%%%%%%%%%%%%%%%%%%%%%%%%%%%%
%%%%%%%%%%%%%%%%%%%%%%%%%%%%%%%%%%%%%%%%%%%%%%%%%%%%%%%%%%%%%%%%%%%%%

\section{Problem \ref{problem2-1}: Free Boundary Problem for ${\vphi}$}\label{sec:free-bdry}
To solve Problem \ref{problem2-1}, we fix the boundary condition on the exit $\Gam_{ex}$ as
\begin{equation}
\label{3-0-7}
A(x,I_n,\grad\vphi)\grad\vphi\cdot\nu_{ex}=v_{ex}
\end{equation}
where $v_{ex}$ is a small perturbation of the constant $v_c$ given by
$$
v_c:=A(x,I_n,\grad\vphi_0^+)D\vphi_0^+\cdot\nu_{ex}|_{\Gam_{ex}}=(B_0-\frac 12|\der_r\vphi_0^+(r_1)|^2)^{\frac{1}{\gam-1}}\der_r\vphi_0^+(r_1).
$$
Here, $\nu_{ex}$ is the outward unit normal on $\Gam_{ex}$ so $\nu_{ex}$ is $\hat r=\frac{x}{|x|}$. Let $(\rho_0,\vphi_0,p_0)$ and $r_s$ be as in Theorem \ref{main-thm'}.
\begin{proposition}
\label{theorem3-1} Under the assumptions of Theorem \ref{main-thm'}(i) and (ii),
for any given $\alp\in(0,1)$, there exist constants $\sigma_3, \delta, C>0$ depending on the data in sense of Remark \ref{remark-cons} so that whenever $0<\sigma\le\sigma_3$, if a function $v_{ex}$ defined on $\Gam_{ex}$ satisfies
$$
\varsigma_4:=\|{v}_{ex}-v_c\|_{1,\alp,\Lambda}^{(-\alp,\der{\Lambda})}\le \sigma,
$$

then there exists a unique solution $\vphi$ to \eqref{3-0-1}, \eqref{3-0-4}, \eqref{3-0-5}, \eqref{3-0-6}, \eqref{3-0-7}  with the shock $S$, and moreover,
\begin{itemize}
\item[(a)]there exists a function $f:\Lambda\to \R^+$ so that $S$, $\mcl{N}^-=\mcl{N}\cap\{|D\Psi^{-1}D\vphi|>c\}$ and $\mcl{N}^+=\mcl{N}\cap\{|D\Psi^{-1}D\vphi|<c\}$ are given by
    \begin{equation*}
    \begin{split}
    &S=\{(r,x')\in\mcl{N}:r=f(x'),x'\in\Lambda\},\\
    &\mcl{N}^-=\mcl{N}\cap\{(r,x'):r<f(x')\},\quad
    \mcl{N}^+=\mcl{N}\cap\{(r,x'):r>f(x')\},
    \end{split}
    \end{equation*}
    and $f$ satisfies
    \begin{equation*}
|f-r_s|_{1,\alp,\Lambda}\le C(\varsigma_1+\varsigma_2+\varsigma_4)\le C\sigma,
\end{equation*}

\item[(b)] $\vphi=\vphi_-$ holds in $\mcl{N}^-$,
\item[(c)] $\vphi$ satisfies the estimate
\begin{equation}
\label{3-0-19}
\|\vphi-\vphi_0^+\|_{2,\alp,\mcl{N}^+}^{(-1-\alp,\Gam_w)}
\le C(\varsigma_1+\varsigma_2+\varsigma_4)\le C\sigma.
\end{equation}

\end{itemize}
\end{proposition}

We prove Proposition \ref{theorem3-1} by the implicit mapping theorem. For that, we first derive the equation and the boundary conditions for $\psi:=\vphi-\vphi_0^+$.
\subsection{Nonlinear boundary problem for $\psi=\vphi-\vphi_0^+$}
If $\vphi$ is a solution in Proposition \ref{theorem3-1} then $\psi:=\vphi-\vphi_0^+$ satisfies
\begin{align}
\label{3-1-1}
&\der_k(a_{jk}(x,D\psi)\der_j\psi)=\der_k F_k(x,D\Psi,D\psi)&\;\;\text{in}\;\;\mcl{N}^+,\\
\label{3-1-4}
&(a_{jk}(x,D\psi)\der_j\psi)\cdot\nu_w=[A(x,I_n,D\vphi)-A(x,D\Psi,D\vphi)]&\;\;\text{on}\;\;\Gam_w^+,\\
\label{3-1-5}
&(a_{jk}(x,D\psi)\der_j\psi)\cdot
\nu_{ex}=v_{ex}-v_c
&\;\;\text{on}\;\;\Gam_{ex}
\end{align}
with $\Gam_w^+:=\Gam_w\cap\der\mcl{N}^+$ and
\begin{align}
\label{3-1-2}
&a_{jk}(x,\eta)=\int_0^1\der_{\eta_k}A_j(x,I_n,\grad\vphi_0^++a\eta)da\\
\label{3-1-3}
&F_k(x,m,\eta)=-A_k(x,m,\grad\vphi_0^++\eta)+A_k(x,I_n,\grad\vphi_0^++\eta)
\end{align}
where $A=(A_1,\cdots,A_n)$ is defined in (\ref{3-0-2}). The boundary condition for $\psi$ on $S$ needs to be computed more carefully.
\begin{remark}
\label{remark-rhos} In this paper, we fix the Bernoulli's invariant $B$ as a constant $B_0$. So,
for given $(\Psi, \vphi_-,p_-)$,  by \eqref{2-b}, $\rho_-$ is given by
 \begin{equation*}
 \rho_-=\frac{\gam p_-}{(\gam-1)(B_0-\frac 12|
 d\Psi^{-1}D\vphi_-|^2)}.
 \end{equation*}
So we regard $(\Psi,\vphi_-,p_-)$ as factors of perturbations of the nozzle, and the incoming supersonic flow.
\end{remark}
\begin{lemma}
\label{lemma3-1-1} Let $f$ and $\vphi$ be as in Proposition
\ref{theorem3-1} and set $\psi_-:=\vphi_--\vphi_0^-$. Then there are two functions $\mu=\mu_f$, $g_1=g_1(x,D\Psi,p_-,\psi_-,\grad\psi_-,\grad\psi)$ in sense of Remark \ref{remark-rhos}, and a vector valued function ${b}_1={b}_1(x,\grad\psi_-,\grad\psi)$ so that $\psi$ satisfies
 \begin{equation}
\label{3-1-17}
{b}_1\cdot\grad\psi-\mu_f\psi=g_1\;\text{on}\;S.
\end{equation}
Moreover, $b_1,\mu_f$ and $g_1$ satisfy, for  $\mu_0>0$ defined in Lemma \ref{lemma-1},
\begin{align}
\label{h2}
&\|\mu_f-\mu_0\|_{1,\alp,\Lambda}\le C(\varsigma_1+\varsigma_2+\varsigma_4),\\
\label{h}
&\|b_1-\nu_s\|_{1,\alp,\Lambda}^{(-\alp,\der \Lambda)}\le C(\varsigma_1+\varsigma_2+\varsigma_4),\\
\label{d4-1}
&\|g_1\|_{1,\alp,\Lambda}^{(-\alp,\der \Lambda)}\le C(\varsigma_1+\varsigma_2).
\end{align}
So if $\sigma_3$(in Proposition \ref{theorem3-1}) is sufficiently small depending on the data, then there is a constant $\lambda>0$ satisfying
\begin{equation}
\label{h3}
\mu_f\ge \frac{\mu_0}{2}>0,\quad \text{and}\quad {b}_1\cdot\nu_s\ge \lambda\quad\text{on}\quad S
\end{equation}
where the constant $\lambda$ also only depends on the data as well.
\begin{proof}
\textbf{Step 1.}
In general, at $(f(x'),x')$ for $x'\in \Lambda$, $\psi=\vphi-\vphi_0^+$ satisfies
\begin{equation}
\label{sl}
\der_r\psi=\der_r\vphi-\frac{K_0}{\der_r\vphi_0^-}+(\frac{K_0}{\der_r\vphi_0^-}-\der_r\vphi_0^+)
\end{equation}
where $K_0$ is given by \eqref{k0}.
By \eqref{2-2-5} and \eqref{3-0-4}, \eqref{sl} is equivalent to
\begin{equation}
\label{new-2}
D\psi\cdot \hat r-\mu_f\psi=h
\end{equation}
with
\begin{align}
\label{3-1-16}
\mu_f(x'):=\frac{\int_0^1\frac{d}{dr}\bigl(\frac{K_0}{\der_r\vphi_0^-}-\der_r\vphi_0^+\bigr)(r_s+t(f(x')-r_s))dt}
{\int_0^1\der_r(\vphi_0^--\vphi_0^+)(r_s+t(f(x')-r_s))dt},\\
\label{new-4}
h=D\vphi\cdot\hat r-\frac{K_0}{D\vphi_0^-\cdot\hat r}-\mu_f(\vphi_--\vphi_0^-)=:h_*-\mu_f\psi_-
\end{align}
where we set $\psi_-=:\vphi_--\vphi_0^-$.

Suppose that $\Psi=Id, (\vphi_-,p_-)=(\vphi_0^-,p_0^-)$ and $\vphi$ in Proposition \ref{theorem3-1} is a radial function with a shock on $S=\{r=f_0\}$ for a constant $f_0\in(r_0,r_1)$. Then, \eqref{2-2-5} implies $h=0$ in \eqref{new-2}.
Then the boundary condition for $\psi$ on $S$ is
\begin{equation}
\label{new-52}
\der_r\psi-\mu_{f_0}\psi=0\;\;\text{on}\;\;S=\{r=f_0\}.
\end{equation}

\textbf{Step 2.} \eqref{new-52} shows how to formulate a boundary condition for $\vphi$ in general.

Since $\frac{K_0}{\der_r\vphi_0^-}-\der_r\vphi_0^+$ is smooth, by Proposition \ref{theorem3-1}, we have
 $$\|\mu_f-\mu_0\|_{1,\alp,\Lambda}\le C\|f-r_s\|_{1,\alp,\Lambda}\le C\|\psi-\psi_-\|_{1,\alp,S}\le C(\varsigma_1+\varsigma_2+\varsigma_4).$$
\eqref{h2} is verified. By Lemma \ref{lemma-1}, $\mu_0$ is a positive constant.

 We note that $h^*$ in \eqref{new-4} depends on $D\Psi, p_-, D\psi_-$ in sense of Remark \ref{remark-rhos}, and also depends on  $D\psi$. So we should carefully decompose $h^*$ in the form of $h^*=\beta\cdot D\psi+g$ so that a weighted H\"{o}lder norm of $g$ does not depend on $\|\psi\|_{2,\alp,\mcl{N}^+}^{(-1-\alp,\Gam_w)}$ if $\sigma_3$ is sufficiently small. For $h^*$ defined in \eqref{new-4}, by \eqref{3-0-5}, we set $h^*=a_1-a_2$ with
 \begin{equation}
 \label{new-5}
 a_1=D\vphi\cdot\hat r-QD\vphi\cdot\nu_s,\quad a_2=\frac{K_0}{D\vphi_0^-\cdot\hat r}-\frac{K_s}{QD\vphi^-\cdot \nu_s}\;\;
 \text{on}\;\;S.
 \end{equation}

 We note that, by \eqref{3-0-4}, $\nu_s-\hat r$ is written as
\begin{equation}
\label{3-1-12}
\nu_s-\hat{r}
=\nu(\grad\psi_-,\grad\psi)-\nu(0,0)=j_1(D\psi_-,D\psi)D\psi_-+j_2(D\psi_-,D\psi)D\psi
\end{equation}
where we set
\begin{equation}
\label{nu}
\nu(\xi,\eta):=\frac{(D\vphi_0^-+\xi)-(D\vphi_0^++\eta)}{|(D\vphi_0^-+\xi)-(D\vphi_0^++\eta)|},
\end{equation}
\begin{equation}
\label{js}
j_1(\xi,\eta):=\int_0^1D_{\xi}\nu(t\xi,\eta)dt\;,\quad
j_2(\xi,\eta):=\int_0^1 D_{\eta}\nu(0,t\eta)dt\;
\end{equation}
for $\xi,\eta\in\R^n$.
 Then, we can rewrite  $a_1$ as
 \begin{equation}
 \label{new-6}
 \begin{split}
 a_1
 &=D\vphi\cdot(\hat r-\nu_s)+(I-Q)D\vphi\cdot\nu_s\\
 &=-[j^T_2(D\psi_-,D\psi)]D(\vphi_0^++\psi)\cdot D\psi+h_1:=\beta_1\cdot D\psi+h_1
 \end{split}
 \end{equation}
 with
 $$
 h_1(D\Psi,D\psi_-,D\psi)=-D\vphi\cdot[j_1(D\psi_-,D\psi)D\psi_-]-(Q-I)D\vphi\cdot\nu_s.
 $$
 Here, all the quantities are evaluated on $(f(x'),x')\in S$.

 To rewrite $a_2$ for our purpose, let us consider $K_0-K_s$ first.
 \begin{equation}
 \label{ks}
 \begin{split}
 K_0-K_s
 &=\frac{\gam-1}{\gam+1}(\underset{=:q_*}{\underbrace{(D\vphi_0^-\cdot\hat r)^2-(QD\vphi_-\cdot\nu_s)^2}})+\frac{2\gam}{\gam+1}(\frac{p_-}{\rho_-}-\frac{p_0^-}{\rho_0^-})
 \end{split}
 \end{equation}
 $\|(\frac{p_-}{\rho_-}-\frac{p_0^-}{\rho_0^-})|_S\|_{1,\alp,\Lambda}^{(-\alp,\der\Lambda)}$ is independent of $\vphi$ but $q_*$ should be treated carefully.
 \begin{equation*}
 \begin{split}
 q_*
 &=|D\vphi_0^-|^2(1-(\hat r\cdot\nu_s)^2)+[(D\vphi_0^-\cdot\nu_s)^2-(QD\vphi_-\cdot\nu_s)^2]=:I_1+I_2
 \end{split}
 \end{equation*}

It is easy to see that  a H\"{o}lder norm of $I_2$ is uniformly bounded independent of $\vphi$. Let us consider $I_1$ now. By \eqref{normal}, a simple calculation provides
\begin{equation}
\label{new-7}
1-(\hat{r}\cdot \nu_s)^2=\frac{|D(\psi_--\psi)|^2-|\der_r(\psi_--\psi)|^2}{|D(\vphi_--\vphi)|^2}=:
\frac{(D\psi-\der_r\psi\hat r)\cdot D\psi}{|D(\vphi_--\vphi)|^2}+I_3
\end{equation}
where a H\"{o}lder norm of $I_3$ is independent of $\psi$ if $\sigma_3$ is sufficiently small. By \eqref{new-5}, \eqref{ks} and \eqref{new-7}, $a_2$ is
\begin{equation}
\label{a2}
\begin{split}
a_2
&=\bigl(\frac{(\gam-1)\der_r\vphi_0^-(D\psi-\der_r\psi\hat r)}{(\gam+1)|D(\vphi_--\vphi)|^2}+\frac{[j_2^T(D\psi_-,D\psi)]K_s\hat r}{QD\vphi_-\cdot\nu_s}\bigr)\cdot D\psi+h_2\\
&=:\beta_2\cdot D\psi+h_2
\end{split}
\end{equation}
with
\begin{equation}
\label{new-8}
\begin{split}
&\phantom{=}h_2=\frac{\gam-1}{\gam+1}(\der_r\vphi_0^- I_3+I_2)+\frac{2\gam}{(\gam+1)\der_r\vphi_0^-}(\frac{p_-}{\rho_-}-\frac{p_0^-}{\rho_0^-})\\
&\phantom{aa}+\frac{K_s}{QD\vphi_-\cdot\nu_s}(\frac{(Q-I)D\vphi_-\cdot\nu_s+D\psi_-\cdot\nu_s}{\der_r\vphi_0^-}
+[j_1^T(D\psi_-,D\psi)]\der_r\psi_-).
\end{split}
\end{equation}

By \eqref{new-6} and \eqref{new-8}, we can write $h^*=a_1-a_2$ as
\begin{equation}
\label{new-9}
h^*=(\beta_1-\beta_2)\cdot D\psi+(h_1-h_2)
\end{equation}
for $\beta_{1,2}=\beta_{1,2}(D\psi_-,D\psi), h_{1,2}=h_{1,2}(D\Psi,D\psi_-,D\psi)$ defined in \eqref{new-6} and \eqref{a2}.
From this, we define $b_1, g_1$ in \eqref{3-1-17} by
\begin{equation}
\label{new-11}
\begin{split}
b_1=\hat r-\beta_1+\beta_2,\quad
g_1=h_1-h_2-\mu_f\psi_-.
\end{split}
\end{equation}
By the definition of $g_1$ in \eqref{new-11} and Proposition \ref{theorem3-1}, we have
$$
\|g_1(D\Psi,p_-,D\psi_-,\psi_-,D\psi)\|_{1,\alp,\Lambda}^{(-\alp,\der\Lambda)}\le C(1+\sigma)(\varsigma_1+\varsigma_2)
$$
so \eqref{d4-1} holds true if we choose $\sigma_3$ satisfying $0<\sigma_3\le 1$.

It remains to verify \eqref{h}. From the definition of $\nu(\xi,\eta)$, one can check
$$D_{\eta}\nu(0,0)=\frac{I_n}{\der_r(\vphi_0^--\vphi_0^+)}
-\frac{D^T(\vphi_0^--\vphi_0^+)D(\vphi_0^--\vphi_0^+)}{|D(\vphi_0^--\vphi_0^+)|^3}$$
so we have $j_2(0,0)\hat r=0$ where $j_2$ is defined in \eqref{js} then \eqref{new-6} and \eqref{a2} imply $\beta_1(0,0)=\beta_1(0,0)=\vec 0$. As vector valued functions of $\xi,\eta\in \R^n$, $\beta_1$ and $\beta_2$ are smooth for $\xi,\eta\in B^n_d(0)$ for a constant $d>0$ sufficiently small. So if $\sigma_3$ is sufficiently small, then we obtain
\begin{equation}
\label{betas}
\|\beta_{1,2}(D\psi_-,D\psi)\|_{1,\alp,\Lambda}^{(-\alp,\der\Lambda)}\le C(\varsigma_1+\varsigma_2+\varsigma_4),
\end{equation}
and this implies
$$
\|b_1(D\psi_-,D\psi)-\nu_s\|_{1,\alp,\Lambda}^{(-\alp,\der\Lambda)}\le
 C(\varsigma_1+\varsigma_2+\varsigma_4)+\|\nu_s-\hat r\|_{1,\alp,\Lambda}^{(-\alp,\der\Lambda)}.
$$
By \eqref{normal} and Proposition \ref{theorem3-1}, $\nu_s-\hat r$ satisfies
$$
\|\nu_s-\hat r\|_{1,\alp,\Lambda}^{(-\alp,\der\Lambda)}\le C(\varsigma_1+\varsigma_2+\varsigma_4),
$$
therefore \eqref{h} is verified provided that we choose $\sigma_3$ in Proposition \ref{theorem3-1} sufficiently small depending on the data.
\end{proof}

\end{lemma}

\subsection{Proof of Proposition \ref{theorem3-1}}
Fix $\alp\in(0,1)$. For a constant $R>0$, keeping Remark \ref{remark-rhos} in mind, we define neighborhoods of $Id$, $(\vphi_0^-,p_0^-)$, $v_c$, $r_s$.
\begin{equation}
\label{sets}
\begin{split}
&{\mcl{B}}^{(11)}_R(Id):=\{\Psi\in C^{2,\alp}(\ol{\mcl{N}},\R^n): |\Psi-Id|_{2,\alp,\mcl{N}}\le R\},\\
&\mcl{B}^{(12)}_R(\vphi_0^-,p_0^-):=\{(\vphi_-,p_-)\in C^{3,\alp}(\ol{\mcl{N}^-_{r_s+\delta}})\times C^{2,\alp}(\ol{\mcl{N}^-_{r_s+\delta}}):\\
&\phantom{aaaaaaaaaaaaaaa}|\vphi^--\vphi_0^-|_{3,\alp,\mcl{N}^-_{r_s+\delta}}+|p^--p_0^-|_{2,\alp,\mcl{N}^-_{r_s+\delta}}\le R\},\\
&\mcl{B}^{(13)}_R(v_c):=\{v_{ex}\in C^{1,\alp}_{(-\alp,\der\Lambda)}(\Lambda):\|v_{ex}-v_c\|_{1,\alp,\Lambda}^{(-\alp,\der\Lambda)}\le R\},\\
&\mcl{B}^{(1)}_R(Id,\vphi_0^-,p_0^-,v_c):=\mcl{B}_R^{(11)}(Id)\times \mcl{B}_R^{(12)}(\vphi_0^-,p_0^-)\times \mcl{B}_R^{(13)}(v_c),\\
&\mcl{B}_R^{(2)}(r_s):=\{f\in C^{2,\alp}_{(-1-\alp,\der\Lambda)}(\Lambda):\|f-r_s\|_{2,\alp,\Lambda}^{(-1-\alp,\der\Lambda)}\le R\}.
\end{split}
\end{equation}

For a sufficiently small $R>0$, $\mcl{B}_{R}^{(1)}(Id,\vphi_0^-,p_0^-,v_c)$ is a set of small perturbations of the nozzle, the incoming supersonic flow and the exit boundary condition for $\vphi$. For two constants $C^*,\sigma^*>0$, define $\mathfrak{J}:\mcl{B}_{\sigma^*}^{(1)}(Id,\vphi_0^-,p_0^-,v_c)\times \mcl{B}^{(2)}_{C^*\sigma^*}(r_s)\to C^{2,\alp}_{(-1-\alp,\der\Lambda)}(\Lambda)$ by
\begin{equation}
\label{definition of J}
\mathfrak{J}(\Psi,\vphi_-,p_-,v_{ex},f):=(\vphi_--\vphi_0^+-\psi)(f(x'),x')
\end{equation}
for $x'\in\Lambda$ where $\psi$ satisfies the following conditions: setting
$$
\mcl{N}^+_f:=\mcl{N}\cap\{r>f(x')\},\quad S_f:=\{r=f(x'), x'\in\Lambda\},
$$
\begin{itemize}
\item[(i)] $\psi$ satisfies
\begin{equation}
\label{M}
\|\psi\|_{2,\alp,\mcl{N}^+_f}^{(-1-\alp,\Gam_w)}\le M\sigma^*
\end{equation}
for a constant $M>0$ to be determined later;
\item[(ii)] $\psi$ solves the nonlinear boundary problem:
\begin{equation}
\label{ift-eqn}
\begin{split}
&N_+(\psi):=\sum_{j,k}\der_k(a_{jk}(x,D\psi)\der_j\psi)=\sum_k\der_kF_k(x,D\Psi,D\psi)\quad\text{in}\quad \mcl{N}^+_f\\
&M^1_f(\psi):=b_1(D\psi,\psi)\cdot D\psi-\mu_f\psi\\
&\phantom{aaaaaa}=g_1(x,D\Psi,p_-, \psi_-,D\psi_{-},D\psi)\quad\text{on}\quad S_f\\
&M^2(\psi):=[a_{jk}(x,D\psi)]D\psi\cdot\nu_w=g_2(x,D\Psi,D\psi)\quad\text{on}\quad \Gam_{w,f}:=\der\mcl{N}^+_f\cap\Gam_w\\
&M^3(\psi):=[a_{jk}(x,D\psi)]D\psi\cdot\nu_{ex}=g_3(v_{ex})\quad\text{on}\quad \Gam_{ex}
\end{split}
\end{equation}
with
\begin{equation}
\label{def-gs}
\begin{split}
&g_2(x,D\Psi,D\psi)=-[A(x,D\Psi,D\vphi_0^++D\psi)-A(x,I_n,D\vphi_0^++D\psi)],\\
&g_3(v_{ex})=v_{ex}-v_c
\end{split}
\end{equation}
\end{itemize}
where $A$ in \eqref{def-gs} is defined as in \eqref{3-0-7}.
For $\sigma^*>0$ sufficiently small, if $M$ in \eqref{M} is appropriately chosen depending on the data so that \eqref{ift-eqn} is a uniformly elliptic boundary problem with oblique boundary conditions, then $\mathfrak{J}$ is well defined. More precisely, we have the following lemma.
\begin{lemma}[Well-posedness \eqref{ift-eqn}]
\label{lemma-nonl}
Fix $\alp\in(0,1)$ and $C^*>0$. Then there is a constant $M, \sigma^{\sharp}$ depending on the data in sense of Remark \ref{remark-cons} and also depending on $C^*$ so that the followings hold:
\begin{itemize}
\item[(i)] for any given $(\Psi,\vphi_-,p_-,v _{ex})\in \mcl{B}_{\sigma^{\sharp}}^{(1)}(Id,\vphi_0^-,p_0^-,v_c), f\in \mcl{B}^{(2)}_{C^{*}\sigma^{\sharp}}(r_s)$, \eqref{ift-eqn} is a uniformly elliptic boundary problem, and it has a unique solution $\psi$ satisfying \eqref{M},
\item[(ii)] if $\sigma^{\sharp}$ is sufficiently small depending on the data and $C^*$, then there is a constant $C$ depending only on the data so that the solution of \eqref{ift-eqn} satisfies
    $$\|\psi\|_{2,\alp,\mcl{N}^+_f}^{(-1-\alp,\Gam_w)}\le C(\varsigma_1+\varsigma_2+\varsigma_4)$$
    where $\varsigma_1, \varsigma_2, \varsigma_4$ are defined in Theorem \ref{proposition2-2-1} and Proposition \ref{theorem3-1}.
\end{itemize}
\end{lemma}
We will prove Lemma \ref{lemma-nonl} by the contraction mapping principle. For that purpose, we first need to study a linear boundary problem related to \eqref{ift-eqn}.
Let $\sigma^*>0$ be a small constant to be determined later. For a given constant $C^*>0$, fix $(\Psi,\vphi_-,p_-,v_{ex})\in \mcl{B}^{(1)}_{\sigma^*}(Id,\vphi_0^-,p_0^-,v_c)$, $f\in \mcl{B}^{(2)}_{C^*\sigma^*}(r_s)$, and let us set $\psi_-=\vphi_--\vphi_0^-$. Reduce $\sigma^*$ sufficiently small depending on the data so that $D\Psi$ is invertible in $\mcl{N}$, and set
$$
\mcl{K}_f(M):=\{\phi\in C^{1,\alp}(\ol{\mcl{N}^+_f}):\|\phi\|_{2,\alp,\mcl{N}^+_f}^{(-1-\alp,\Gam_w)}\le M\sigma^*\}.
$$
For a fixed $\phi\in \mcl{K}_f(M)$, consider the following linear problem in $\mcl{N}^+_f$:
\begin{equation}
\label{linear-p}
\begin{split}
&\sum_{j,k}\der_k(a_{jk}(x,D\phi)\der_j u)=\sum_k\der_kF_k(x,D\Psi,D\phi)\quad\text{in}\quad \mcl{N}^+_f,\\
&b_1(D\phi,\phi)\cdot D u-\mu_f u=g_1(x,D\Psi,p_-,\psi_-,D\psi_-,D\phi)\quad\text{on}\quad S_f,\\
&(a_{jk}(x,D\phi)\der_j u)\cdot\nu_w=g_2(x,D\Psi,D\phi)\quad\text{on}\quad \der\mcl{N}^+_f\cap\Gam_w,\\
&(a_{jk}(x,D\phi)\der_j u)\cdot\nu_{ex}=g_3(v_{ex})\quad\text{on}\quad \Gam_{ex}.
\end{split}
\end{equation}
with $\psi_-:=\vphi_--\vphi_0^-$.
The definitions of $a_{jk}, b_1, \mu_f, g_1, g_2, g_3$ are given in \eqref{3-1-2}, \eqref{3-1-16}, \eqref{new-11}, and \eqref{def-gs}.
\begin{lemma}[Well posedness of \eqref{linear-p}]
\label{lemma-linear-existence}
For any given $\alp\in(0,1)$,  and $M, C^*>0$ is $\sigma^*>0$ is chosen sufficiently small depending only on the data, $C^*$ and $M$, then there exists a constant $C$ depending on the data and $\alp$ so that, for any given  $f\in \mcl{B}_{C^*\sigma^*}^{(2)}(r_s)$ and $\phi\in\mcl{K}_f(M)$, \eqref{linear-p} has a unique solution $u^{(\phi)}\in C^{2,\alp}_{(-1-\alp,\Gam_w)}(\mcl{N}^+_f)$. Moreover $u^{(\phi)}$ satisfies
\begin{equation}
\label{3-2-8}
\|u^{(\phi)}\|_{2,\alp,\mcl{N}^+_{f}}^{(-1-\alp,\Gam_w)}\le C(\varsigma_1+\varsigma_2+\varsigma_4).
\end{equation}
\end{lemma}
To prove Lemma \ref{lemma-linear-existence}, we basically repeat the proof of \cite[Theorem 3.1]{Ch-F3}.
A key ingredient of the proof in \cite{Ch-F3} is the method of continuity and the divergence structure of a linear boundary problem. In our case, however, \eqref{linear-p} is not in a divergence form because of the derivative boundary condition on $S_f$. For that reason, we need Lemma \ref{lemma-ver2} so that \eqref{linear-p} can be treated as a divergence problem. Before we proceed further, we first note two obvious lemmas for the later use.
\begin{lemma}
\label{lemma-2}
For any $C^*, M>0$, and for $f\in \mcl{B}^{(2)}_{C^*\sigma^*}(r_s), \phi\in\mcl{K}_f(M)$, there hold
\begin{align}
\label{coeff-est1}
&\|a_{jk}(x,D\phi)-a_{jk}(x,0)\|_{1,\alp,\mcl{N}^+_f}^{(-\alp,\Gam_w)}\le CM\sigma^*,\\
\label{coeff-est2}
&\|b_1(D\phi,\phi)-b_1(0,0)\|_{1,\alp,\mcl{N}^+_f}^{(-\alp,\Gam_w)}\le CM\sigma^*,\\
\label{coeff-est3}
&\|\mu_f-\mu_0\|_{1,\alp,\Lambda}\le C\|f-r_s\|_{1,\alp,\Lambda},\\
&\|F(x,D\Psi)\|_{\alp,\mcl{N}^+_f}^{(1-\alp,\Gam_w)}+\|g_1(x,D\Psi,p_-,\psi_-,D\psi_{-},D\phi)\|_{1,\alp,\Lambda}^{(-\alp,\der\Lambda)}\notag\\
\label{coeff-est4}
&\phantom{aaaaa}+
\|g_2(x,D\Psi,D\phi)\|_{\alp,\Gam_{w,f}}+\|g_3(v_{ex})\|_{1,\alp,\Lambda}^{(-\alp,\der\Lambda)}
\le C\sigma^*.
\end{align}
where the constant $C$ in \eqref{coeff-est1}-\eqref{coeff-est4} only depends on the data in sense of Remark \ref{remark-cons}.
Thus, for any $C^*, M>0$, if $\sigma^*$ is chosen sufficiently small depending on the data, $C^*$ and $M$, then the equation of \eqref{linear-p} is uniformly elliptic in $\mcl{N}^+_f$.
\begin{proof}
\eqref{coeff-est1}-\eqref{coeff-est4} are obvious by the definitions.

For any given $\zeta\in\R^n$, $a_{ij}(x,0)$ satisfies
\begin{equation}
\label{new-40}
a_{ij}(x,0)\zeta_i\zeta_j
\ge \frac{(\gam+1)}{2(\gam-1)}(B_0-\frac 12|\grad\vphi_0^+|^2)^{\frac{2-\gam}{\gam-1}}(K_0-|\grad\vphi_0^+|^2)|\zeta|^2.
\end{equation}
So, by \eqref{2-3-6}, \eqref{2-2-8} and Remark \ref{remark2-2-1}, for any $C^*, M>0$, one can choose $\sigma^*$ depending on the data, $C^*$ and $M$ so that \eqref{linear-p} is uniformly elliptic in $\mcl{N}^+_f$.
\end{proof}
\end{lemma}

\begin{lemma}
\label{lemma-3}
For any $C^*, M>0$, let $\nu_f$ be the inward unit normal on $S_f$, and $b_1$ be defined in \eqref{new-11}.  Then, for any given$f\in \mcl{B}^{(2)}_{C^*\sigma^*}(r_s)$, $\phi\in \mcl{K}_f(M)$, we have
\begin{align}
\label{obliqueness}
&b_1(D\phi,\phi)\cdot\nu_f\ge \frac 12\quad\text{on}\quad S_f,\\
\label{zeroorder}
&\mu_f\ge \frac 12 \mu_0>0\quad\text{on}\quad S_f\quad\text{for any}\;\;f\in \mcl{B}_{C^*\sigma^*}^{(2)}(r_s)
\end{align}
where $\mu_0$ is defined in Lemma \ref{lemma-1} if $\sigma^*$ is chosen sufficiently small depending on the data, $C^*$ and $M$.
\begin{proof}
For $\nu$ defined in \eqref{nu}, we have
\begin{align*}
&\phantom{=}b_1(D\phi,\phi)\cdot\nu_f\\
&=(b_1(D\phi,\phi)-b_1(0,0))\cdot\nu_f+b_1(0,0)\cdot(\nu_f-\nu(0,0))+b_1(0,0)\cdot\nu(0,0)\\
&\ge -C(M+C^*)\sigma^*+1
\end{align*}
for a constant $C$ depending only on the data. By \eqref{coeff-est2}, for any $C^*, M>0$, if we choose $\sigma^*$ to satisfy
 $-C(M+C^*)\sigma^*+1\ge\frac 12$, then \eqref{obliqueness} is obtained. Reducing $\sigma^*$ if necessary, \eqref{coeff-est4} provides \eqref{zeroorder}.
\end{proof}
\end{lemma}

By Lemma \ref{lemma-2} and \ref{lemma-3}, for any $ C^*, M>0$,  choosing $\sigma^*$ sufficiently small depending on the data, $C^*$ and $M$, \eqref{3-1-1} is a uniformly elliptic linear boundary problem with oblique boundary conditions on the boundary.
\begin{proof}[\textbf{Proof of Lemma \ref{lemma-linear-existence}}]\quad\\
\textbf{Step 1.} Suppose that $u$ solves \eqref{linear-p} and then we first verify \eqref{3-2-8} by following the method of the proof for \cite[Theorem 3.1]{Ch-F3}. First, we apply the change of variables $x\mapsto y=\Phi(x)$  with $\Phi$ defined by
\begin{equation}
\label{ver3-1}
|y|=r_s+\frac{r_1-r_s}{r_1-f{(x')}}(|x|-f(x')),\qquad \frac{y}{|y|}=\frac{x}{|x|},
\end{equation}
This change of variables transforms $\mcl{N}^+_f$ to $\mcl{N}^+_{r_s}$, and $S_f=\{r=f(x')\}\cap\mcl{N}^+_f$ to $S_0=\{r=r_s\}\cap \mcl{N}^+_{r_s}$ in the spherical coordinates in \eqref{spherical-coord}. Then $w(y)=u\circ \Phi^{-1}(y)$ solves a linear boundary problem of the same structure as \eqref{linear-p}. For simplicity, we write the the equation and the boundary conditions for $w$ as
\begin{equation}
\label{d4-2'}
\begin{split}
%\label{3-2-3}
&\der_k({\til a}_{jk}(y,f,D\phi)\der_j w)=\der_k \til F_k(y,f,D\Psi,D\phi)\;\;\text{in}\;\;\mcl{N}^+_{r_s},\\
%\label{3-2-4}
&\til b_1(y,f,\phi,D\phi)\cdot\!\grad w-\mu_fw={\til g}_1(y,f,D\Psi,p_-,\psi_-,D\psi_-,D\phi)\;\;\text{on}\;\; S_0\\
%\label{3-2-5}
&({\til a}_{jk}(y,f,D\phi)\der_j w)\cdot\nu_w={\til g}_2(D\Psi,f,D\phi)\;\;\text{on}\;\;\Gam_{w,r_s},\\
%\label{3-2-6}
&({\til a}_{jk}(y,f,D\phi)\der_j w)\cdot\nu_{ex}=\til g_3(v_{ex})\;\;\text{on}\;\;\Gam_{ex}.
\end{split}
 \end{equation}
with $\Gam_{w,r_s}=\Gam_w\cap \der\mcl{N}^+_{r_s}$ and $\psi_-=\vphi_--\vphi_0^-$, where $\til a_{jk}, \til b_1, \til F_k, \til g_m$ are obtained from $a_{jk}, b_1, F_k, g_m$ in \eqref{linear-p} through the change of variables $x\mapsto y=\Phi(x)$ for $j,k=1,\cdots, n$, and $m=1,2,3$. Next, we rewrite \eqref{d4-2'} as a conormal boundary problem with coefficients $a_{jk}(y,0)$ and $b_1(y,0)=\hat r$ for the equation and the boundary conditions of $w$.

We remind that, by \eqref{3-1-2}, the matrix $[a_{jk}(x,0)]$ is
\begin{equation}
\label{matrixat0}
[a_{jk}(y,0)]=(B_0-\frac 12|D\vphi_0^+|^2)^{\frac{1}{\gam-1}}I_n-
\frac{(B_0-\frac 12|D\vphi_0^+|^2)^{\frac{2-\gam}{\gam-1}}}{\gam-1}(D\vphi_0^+)^TD\vphi_0^+.
\end{equation}
\begin{lemma}[Boundary condition for $w$ on $S_0$]
 \label{lemma-ver2}
 Suppose that $w$ solves \eqref{d4-2'}, and let $M,C^*,\sigma^*$ be as in Lemma \ref{lemma-linear-existence}.
 For any given $\alp\in(0,1)$, $f\in B^{(2)}_{C^*\sigma^*}(r_s),$ $\phi\in\mcl{K}_f(M)$, there is a constant $C$ depending only on the data in sense of Remark \ref{remark-cons} so that $w$ satisfies

\begin{equation}
\label{BC-shock}
(a_{jk}(y,0)\der_j w)\cdot\nu_s=h_1\quad\text{on}\;\;S_0
\end{equation}
for a function $h_1\in C^{1,\alp}_{(-\alp,\der\Lambda)}(\Lambda)$, and $h_1$ satisfies
 \begin{equation}
 \label{new-43}
 \|h_1\|_{1,\alp,\Lambda}^{(-\alp,\der\Lambda)}
 \le C((M+C^*)\sigma^*
 \|w\|_{2,\alp,\mcl{N}^+_{r_s}}^{(-1-\alp,\Gam_w)}
+\|w\|_{1,\alp,\mcl{N}^+_{r_s}}^{(-\alp,\Gam_w)}+\|\til g_1\|_{1,\alp,\Lambda}^{(-\alp,\der\Lambda)}).
\end{equation}

 \begin{proof}
The unit normail on $S_0$ toward $\mcl{N}^+_{r_s}$ is $\hat r$. Then, by \eqref{matrixat0} and the boundary condition of $w$ on $S_0$ in \eqref{d4-2'}, there holds
\begin{equation*}
\begin{split}
\phantom{=}(a_{jk}(y,0)\der_j u)\cdot\nu_s
&=(B_0-\frac 12|D\vphi_0^+|^2)^{\frac{1}{\gam-1}}(1-\frac{|D\vphi_0^+|^2}{(\gam-1)(B_0-\frac 12|D\vphi_0^+|^2)})\der_r w\\
&=\frac{(\gam-1)B_0-\frac{\gam+1}{2}|D\vphi_0^+|^2}{(\gam-1)(B_0-\frac 12|D\vphi_0^+|^2)^{1-\frac{1}{\gam-1}}}[(\hat r-\til b_1)\cdot Dw+\mu_fw+\til g_1]=: h_1
\end{split}
\end{equation*}
By \eqref{coeff-est2} and the definition of $\til b_1$, \eqref{new-43} is easily obtained.
 \end{proof}
 \end{lemma}
Back to the proof of Lemma \ref{lemma-linear-existence}, using Lemma \ref{lemma-ver2}, $w$ solves the conomal boundaru problem
\begin{equation}
\label{neweq}
\begin{split}
&\der_k(a_{jk}(y,0)\der_j w)=\der_k[\til F_k(y,f,D\Psi,D\phi)+\delta \til a_{jk}\der_j w]=:\der_k G_k\;\;\text{in}\;\;\mcl{N}^+_{r_s},\\
&(a_{jk}(y,0)\der_j u)\cdot\nu_s=h_1\;\;\text{on}\;\;S_0,\\
&\der_k(a_{jk}(y,0)\der_j w)\cdot \nu_w=\til g_2+(\delta \til a_{jk}\der_j w)\cdot \nu_w=:h_2\;\;\text{on}\;\;\Gam_{w,r_s},\\
&\der_k(a_jk(y,0)\der_j w)\cdot \nu_{ex}=\til g_3+(\delta \til a_{jk}\der_j w)\cdot \nu_{ex}=:h_3\;\;\text{on}\;\;\Gam_{ex}
\end{split}
\end{equation}
with $\delta \til a_{jk}:=a_{jk}(y,0)-\til a_{jk}(y,f,D\phi)$. We note that, by \eqref{coeff-est1}, $\delta \til a_{jk}$ satisfies
\begin{equation}
\label{new-42}
\|\delta \til a_{jk}\|_{1,\alp,\mcl{N}^+_{r_s}}^{(-\alp,\Gam_w)}\le C(\|\phi\|_{2,\alp,\mcl{N}^+_{f}}^{(-1-\alp,\Gam_w)}+\|f-r_s\|_{2,\alp,\Lambda}^{(-1-\alp,\der\Lambda)})
\end{equation}
for a constant $C$ depending on the data in sense of Remark \ref{remark-cons}.

\textbf{Step 2.} Using a weak formulation of the equation for $w$ in \eqref{neweq} with \eqref{matrixat0}, one can check
\begin{equation}
\label{00}
\begin{split}
&\der_{k}(a_{jk}(y,0)\der_{j}w)=\der_{k}G_k\\
&\Leftrightarrow
L_0:=\der_r(k_1(r)\omega_1(x')\der_r w)+\der_{x'_l}(k_2(r)b_{lm}(x')\der_{x'_m}w)=\der_r(\til G\cdot \hat r)+\sum_{l=1}^{n-1}\der_{x'_l}(\til G\cdot \hat {x'_l})
\end{split}
\end{equation}
with
\begin{equation}
\label{d8-1}
\begin{split}
&k_1(r):=\frac{(\gam+1)r^{n-1}(K_0-(\der_r\vphi_0^+(r))^2)}{2(\gam-1)(B_0-\frac 12(\der_r\vphi_0^+(r))^2)^{\frac{\gam-2}{\gam-1}}},\;\;
k_2(r):=r^{n-3}(B_0-\frac 12(\der_r\vphi_0^+)^2)^{\frac{1}{\gam-1}}\\
&r^{n-1}\omega_1(x')=\det \frac{\der x}{\der(r,x')}
\end{split}
\end{equation}
for $(r,x')$ is given as in \eqref{spherical-coord}, and every $b_{lm}(x')$  is smooth in $\Lambda$. Here we set $\til G=\omega_1[\frac{\der(r,x')}{\der x}]^T G$.

We note that $\omega_1$ is bounded below by a positive constant in $\ol{\Lambda}$, and the linear operator $L_0$ is uniformly elliptic in $\mcl{N}^+_{r_s}$ ,and bounded by \eqref{new-40} and \eqref{2-2-8}.

\textbf{Step 3.} To obtain the desired estimate of $u$ in \eqref{linear-p}, we use \eqref{neweq}-\eqref{d8-1}, Lemma \ref{lemma-ver2} and the interpolation inequality of H\"{o}lder norms to follow the method of the proof for \cite[Theorem 3.1]{Ch-F3}, and prove the following lemma.
\begin{lemma}
\label{lemma-u-est} There exists a constant $\kappa>0$ depending only on the data in sense of Remark \ref{remark-cons} so that whenever $M\sigma^*,C^*\sigma^*\le \kappa$,
for any $f\in B^{(2)}_{\sigma}(r_s)$, $\phi\in\mcl{K}_f(M)$ with $\alp\in(0,1)$, there is a constant $C$ depending on the data so that, if $w$ solves \eqref{d4-2'} then there holds
\begin{equation}
\label{estimate1}
|w|_{1,\alp,\mcl{N}^+_{r_s}}\le C(|w|_{0,\mcl{N}^+_{r_s}}+|\til F|_{0,\alp,\mcl{N}^+_{r_s}}+|\til g_1|_{0,\alp,S_{0}}+|\til g_2|_{0,\alp,\Gam_{w,r)s}}+|\til g_3|_{0,\alp,\Gam_{ex}}).
\end{equation}
\end{lemma}
The only difference to prove Lemma \eqref{lemma-u-est} from \cite[Theorem3.1]{Ch-F3} is that, to apply the method of reflection as in \cite[Theorem3.1]{Ch-F3}, we need to extend $a_{jk}(y,0)$ in $r$-direction by reflection on $S_0$ and $\Gam_{ex}$. For that purpose, we use \eqref{d8-1}. Except for that, the proof of Lemma \ref{lemma-u-est} is basically same as the proof for \cite[Theorem3.1]{Ch-F3} so we skip details.

Once \eqref{estimate1} is obtained, then the standard scaling provides $C^{2,\alp}_{(-1-\alp,\Gam_w)}$ estimate of $w$. For a fixed point $y_0\in{\ol{\mcl{N}^+_{r_s}}\setminus \Gam_{w,r_s}}$ let $2d:=dist(y_0,\Gam_{w,r_s})$, and we set a scaled function $w^{(y_0)}(s)$ near $y=y_0$ by
    $$
    w^{(y_0)}(s):=\frac{w(y_0+ds)-w(y_0)}{d^{1+\alp}}
    $$
    for $s\in B_1^{(y_0)}:=\{s\in B_1(0):y_0+ds\in{\ol{\mcl{N}^+_{r_s}}\setminus \Gam_{w,r_s}}\}$. By \cite[Theorem 6.29]{GilbargTrudinger} and Lemma \ref{lemma-ver2}, for each $y_0\in\ol{\mcl{N}^+_{r_s}}\setminus \Gam_{w,r_s}$, $w^{(y_0)}$ satisfies
    \begin{equation}
    \label{est6}
    \begin{split}
    |w^{(y_0)}|_{2,\alp,B^{(y_0)}_{1/2}}
    &\le C(|w|_{1,\alp,\mcl{N}^+_{r_s}}
    +\|\til{F}\|^{(-\alp,\Gam_w)}_{1,\alp,\mcl{N}^+_{r_s}}+\|\til g_1\|_{1,\alp,\Lambda}^{(-\alp,\der \Lambda)}
    +\|\til g_3\|_{1,\alp,\Lambda}^{(-\alp,\der\Lambda)})\\
    &\le C(|w|_{0,\mcl{N}^+_{r_s}}+\varsigma_1+\varsigma_2+\varsigma_4)
    \end{split}
    \end{equation}
    where the second inequality above holds true by \eqref{d4-1}, \eqref{def-gs} and \eqref{coeff-est4}.

    By \eqref{coeff-est3} and Lemma \ref{lemma-1}, reducing $\sigma^*$ if necessary, we can have $\mu_f\ge \frac 12\mu_0>0$ on $\ol{\Lambda}$. Then, by \eqref{est6}, Remark \ref{remark-holdernorms} and \cite[Corollary 2.5]{Lie-1},  we conclude that if $w$ is a solution to \eqref{linear-p} then it satisfies \eqref{3-2-8}.

\textbf{Step 4.} Finally, we prove the existence of a solution to \eqref{d4-2'}. Consider the following auxiliary problem:
\begin{equation}
\label{auxiliary}
\begin{cases}
\der_k(a_{jk}(y,0)\der_j u)-u=\der_k F_k&\;\text{in}\;\mcl{N}^+_{r_s}\\
(a_{jk}(y,0)\der_j u)\cdot\nu_s=g_1&\;\text{on}\;S_{r_s}\\
(a_{jk}(y,0)\der_j u)\cdot\nu_w=g_2&\;\text{on}\;\Gam_{w,r_s}\\
(a_{jk}(y,0)\der_j u)\cdot\nu_{ex}=g_3&\;\text{on}\;\Gam_{ex}
\end{cases}
\end{equation}
for $F\in C^{\alp}(\ol{\mcl{N}^+_{r_s}})$, $g_1, g_3\in C^{\alp}(\ol{\Lambda})$ and $g_2\in C^{\alp}(\ol{\Gam_{w,r_s}})$. \eqref{auxiliary} has a unique weak solution because the functional
$$I[u]=
\frac 12\int_{\mcl{N}^+_{r_s}}a_{jk}(y,0)\der_j u\der_k u+u^2-2F\cdot Du
+\int_{\der\mcl{N}^+_{r_s}}(F\cdot \nu_{out}+g)u
$$
has a minimizer over $W^{1,2}(\mcl{N}^+_{r_s})$ with the outward unit normal $\nu_{out}$ of $\der\mcl{N}^+_{r_s}$ where we write
$$g:=g_1\chi_{S_0}+g_2\chi_{\Gam_{w,r_s}}-g_3\chi_{\Gam_{ex}}$$
with $\chi_{\mcl{D}}$ defined by
$$\chi_{\mcl{D}}(y)=1\;\;\text{for}\;\;y\in \mcl{D},\quad\chi_{\mcl{D}}(y)=0\;\;\text{otherwise}\quad\text{for}\;\;\;\mcl{D}\subset \R^n.$$

For $F\in C^{1,\alp}_{(-\alp,\Gam_w)}(\mcl{N}^+_{r_s}),$  $g_1, g_3\in C^{1,\alp}_{(-\alp,\der \Lambda)}(\Lambda)$, the weak solution $u$ also satisfies \eqref{3-2-8}. By the method of continuity, we conclude that the linear boundary problem \eqref{linear-p} has a solution $u\in C^{2,\alp}_{(-1-\alp,\Gam_w)}(\mcl{N}^+_{f})$. By the comparison principle and \cite[Corollary 2.5]{Lie-1}, the solution to \eqref{linear-p} is unique.
\end{proof}

\begin{remark}
\label{laplacian-sphere} According to the definition of the laplacian $\vartriangle_{S^{n-1}}$ on a sphere $\mathbb{S}^{n-1}$, we note
$$
\frac{1}{\omega^1}\der_{y'_m}(b_{lm}\der_{y'_l}u)=\vartriangle_{S^{n-1}}u
$$
where we replace $-\vartriangle_{S^{n-1}}$ by $\vartriangle_{S^{n-1}}$ in the widely used convention(e.g. \cite{B}) of $\vartriangle_{S^{n-1}}$.
\end{remark}
%%%%%%%%%%%%%%%%%%%%%%%%%%%%%%%%%%%%%%%%%
Lemma \ref{lemma-nonl} is a direct result from Lemma \ref{lemma-linear-existence}.
\begin{proof}[\textbf{Proof of Lemma \ref{lemma-nonl}}]
For $M,C^*>0$, let $\sigma^*$ be as in Lemma \ref{lemma-linear-existence}. Fix $(\Psi,\rho_-,\vphi_-,p_-,v_{ex})\in \mcl{B}^{(1)}_{\sigma_*}(Id,\rho_0^-,\vphi_0^-,p_0^-,v_c)$ and $f\in \mcl{B}^{(2)}_{C^*\sigma^*}(r_s)$, and define a nonlinear mapping $\mcl{G}_f$ by
$$\mcl{G}_f:\phi\in\mcl{K}_f(M)\mapsto u^{(\phi)}\in C^{2,\alp}_{(-1-\alp,\Gam_w)}(\mcl{N}^+_f)$$
where $u^{(\phi)}$ is a unique solution to \eqref{linear-p} satisfying \eqref{3-2-8}. By Lemma \ref{lemma-linear-existence}, if we choose $M$ as
\begin{equation}
\label{choice-M}
M=6C
\end{equation}
for the constant $C$ in \eqref{3-2-8}, then $\mcl{G}_f$ maps $\mcl{K}_f(M)$ into itself which is a Banach space.

For $k=1,2$, let us set $w_k=\mcl{G}_f(\phi_k)$. By subtracting the boundary problem \eqref{d4-2'} for $w_2$ from the boundary problem \eqref{d4-2'} for $w_1$, one can directly show that there is a constant $\til C$ depending on the data so that there holds
\begin{equation}
\label{contraction}
\|\mcl{G}_f(\phi_1)-\mcl{G}_f(\phi_2)\|_{2,\alp,\mcl{N}^+_f}^{(-1-\alp,\Gam_w)}\le \til C\sigma^*\|\phi_1-\phi_2\|_{2,\alp,\mcl{N}^+_f}^{(-1-\alp,\Gam_w)}.
\end{equation}
So if we reduce $\sigma_*>0$ satisfying Lemma \ref{lemma-2}, \eqref{obliqueness} ,\eqref{zeroorder} and $\til C\sigma^*<1$, then the contraction mapping principle applies to $\mcl{G}_f$. Therefore $\mcl{G}_f$ has a unique fixed point, say, $\psi$. Obviously, $\psi$ is a solution to \eqref{ift-eqn} and it also satisfies Lemma \ref{lemma-nonl}(ii). From this, we choose $\sigma^{\sharp}=\sigma^*$. We note that the choice of $\sigma^{\sharp}$  depends on the data( in sense of Remark \ref{remark-cons}) and $C^{*}$ .
\end{proof}

By Lemma \ref{lemma-nonl}, $\mathfrak{J}$, defined in \eqref{definition of J}, is well-defined when we replace $\sigma^*$ by $\sigma^{\sharp}$ in \eqref{definition of J}.

For any given $(\Psi,\vphi_-,p_-,v_{ex})\in \mcl{B}_{\sigma^{\sharp}}^{(1)}(Id,\vphi_0^-,p_0^-,v_c)$ and $f\in \mcl{B}^{(2)}_{C^{*}\sigma^{\sharp}}(r_s)$, let $\psi^f$ be the fixed point of $\mcl{G}_f$. Suppose that $f^*\in\mcl{B}^{(2)}_{C^{*}\sigma^{\sharp}}(r_s)$ satisfies
$$
\mathfrak{J}(\Psi,\vphi_-,p_-,v_{ex},f^*)=(\vphi_--\vphi_0^+-\psi^{f^*})(f^*(x'),x')=0 \;\;\text{in}\;\;\Lambda.
$$
Then
\begin{equation}
\label{new-45}
\vphi^{f^*}:=\begin{cases}
\vphi_-&\text{in}\;\;{\ol{\mcl{N}^-_{f^{*}}}}\\
\vphi_0^++\psi^{f^*}&\text{in}\;\;{\mcl{N}^+_{f^*}}
\end{cases}
\end{equation}
is a transonic shock solution of Problem \ref{problem2} with a shock $S=\{r=f^*(x'), x'\}\cap \mcl{N}.$
Moreover, $\vphi^{f^*}$ satisfies Proposition \ref{theorem3-1}(a)-(c). Therefore, if there is a mapping $\mathfrak{S}:q=(\Psi,\vphi_-,p_-,v_{ex})\mapsto f$ so that $\mathfrak{J}(q,\mathfrak{S}(q))=0$ holds, then Proposition \ref{theorem3-1} is proven by Lemma \ref{lemma-nonl}.

\begin{lemma}\label{lemma-frechet}
Fix $\alp\in(0,1)$ and $C^*>0$. For $\mathfrak{J}$ defined in \eqref{definition of J}, there holds
\begin{equation}
\label{new-44}
\mathfrak{J}(Id,\vphi_0^-,p_0^-,v_c,r_s)=0.
\end{equation}
Also, there exists a constant $\delta_*>0$ depending on the data and $C^*$ so that
\begin{itemize}
\item[(i)] $\mathfrak{J}$ is continuously Fr\'{e}chet differentiable in $\mcl{B}^{(1)}_{\delta_*}(Id,\vphi_0^-,p_0^-,v_c)\times \mcl{B}^{(2)}_{C^*\delta_*}(r_s)$,
\item[(ii)] $D_f\mathfrak{J}(Id,\vphi_0^-,p_0^-,v_c,r_s):C^{2,\alp}_{(-1-\alp,\der \Lambda)}(\Lambda)\to C^{2,\alp}_{(-1-\alp,\der \Lambda)}(\Lambda)$ is invertible.
\end{itemize}
\begin{proof}
\eqref{new-44} is obvious.

Fix $q_0=(\Psi,\vphi_-,p_-,v_{ex})\!\in \mcl{B}_{\sigma^{\sharp}/2}^{(1)}(Id,\vphi_0^-,p_0^-,v_c), f\!\in \mcl{B}^{(2)}_{C^{*}\sigma^{\sharp}/2}(r_s)$, $\til q\!=(\til\Psi,\til\vphi_-,\til p_-,\til v_{ex})\in C^{2,\alp}(\ol{\mcl{N}},\R^n)\times C^{3,\alp}(\ol{\mcl{N}^-_{r_s+\delta}})\times C^{2,\alp}(\ol{\mcl{N}^-_{r_s+\delta}})\times C^{1,\alp}_{(-\alp,\der\Lambda)}(\Lambda), \til f\in C^{2,\alp}_{(-1-\alp,\der\Lambda)}(\Lambda)$ with
\begin{equation}
\label{qf}
\begin{split}
&\|\til q\|_1:=\|\til\Psi\|_{2,\alp,\mcl{N}}+\|\til\vphi_-\|_{3,\alp,\mcl{N}^-_{r_s+\delta}}+\|\til p_-\|_{2,\alp,\mcl{N}^-_{r_s+\delta}}
+\|\til v_{ex}\|_{1,\alp,\Lambda}^{(-\alp,\der\Lambda)}=1,\\
&\|\til f\|_2:=\|\til f\|_{2,\alp,\Lambda}^{(-1-\alp,\der\Lambda)}=1.
\end{split}
\end{equation}
To prove the Fr\'{e}chet differentiability of $\mathfrak{J}$, we compute
\begin{equation}
\label{formal}
\frac{\der}{\der \eps}|_{\eps=0}\mathfrak{J}(q_0+\eps \til q, f+\eps \til f)
=\frac{\der}{\der\eps}|_{\eps=0}(\vphi_-+\eps\til\vphi_--\vphi_0^+-\psi_{\eps})(f(x')+\eps\til f(x'),x')
\end{equation}
where $\psi_{\eps}$ satisfies \eqref{linear-p} in $\mcl{N}^+_{f+\eps\til f}$ with replacing $f, \phi, \Psi, p_-, v_{ex}$ by $f+\eps\til f, \psi_{\eps}, \Psi+\eps \til\Psi, p_-+\eps \til p_-, v_{ex}+\eps \til v_{ex}$ respectively. Since the case of $(\til\rho_-,\til{\vphi}_-,\til p_-,\til v_{ex})\neq(0,0,0,0)$ can be handled similarly, we assume $(\til\rho_-,\til{\vphi}_-,\til p_-,\til v_{ex})=(0,0,0,0)$ for simplification.

To compute $\frac{\der}{\der\eps}|_{\eps=0}\psi_{\eps}(f(x')+\eps\til f(x'),x')$, we use $\til\Phi_{f+\eps\til f}$, given by
\begin{equation}
\label{changeofvar}
\til\Phi_{f+\eps\til f}:(r,x')\mapsto (\frac{r_1-f}{r_1-f-\eps \til f}(r-f(x')-\eps\til f(x'))+f(x'),x'),
\end{equation}
so that $\til\psi_{\eps}:=\psi_{\eps}\circ (\Phi_{f+\eps\til f})^{-1}$ is defined in $\mcl{N}^+_f$ for every $\eps$ with $|\eps|<\frac{\sigma^{\sharp}}{2}$. The corresponding boundary problem for $\til\psi_{\eps}$ is as in \eqref{d4-2'} with replacing $f, \phi, \Psi, p_-, v_{ex}$ by $f+\eps\til f, \psi_{\eps}, \Psi+\eps \til\Psi, p_-+\eps \til p_-, v_{ex}+\eps \til v_{ex}$ respectively.

$a_{jk}(x,\eta), F(x,m,\eta)$, $b_1(x,\xi,\eta)$, $g_1(x,m,p,\xi,\eta)$, $g_2(x,m,\eta)$, $g_3(x,z)$, defined in \eqref{3-1-2}, \eqref{3-1-3}, \eqref{new-11} and \eqref{def-gs}, are smooth with respect to their variables. So, subtraction the boundary problem for $\psi_0$ from the boundary problem for $\til\psi_{\eps}$, one can show, by the elliptic estimates, that $\frac{\til\psi_{\eps}-\psi_0}{\eps}$ converges in $C^{2,\alp}_{(-1-\alp,\Gam_w)}(\mcl{N}^+_f)$ to $u$ where $u$ solves the following linear boundary problem

\begin{equation}
\label{zero-bd}
\begin{split}
&\der_k((a_{jk}(x,D\psi^0)+o^{(1)}_{jk})\der_j u)=\der_k (a_{1}D\til f+a_{2} D\til\Psi) \;\;\text{in}\quad\mcl{N}^+_f,\\
&(b_1(D(\vphi^--\vphi_0^-),D\psi^0)+o^{(2)})\cdot Du-\mu_f u=a_{3}\til f+a_{4} D\til f+a_{5} D\til\Psi \;\;\text{on}\quad S_f,\\
&((a_{jk}(x,D\psi^0)+o_{jk}^{(3)})\der_ju)\cdot \nu_w=a_{6} D\til f+a_{7}D\til\Psi \;\;\text{on}\quad \Gam_{w,f},\\
&((a_{jk}(x,D\psi^0)+o_{jk}^{(4)})\der_ju)\cdot \nu_{ex}=\til v_{ex}+a_{8} D\til f+a_{9} D\til\Psi \;\;\text{on}\quad \Gam_{ex}.
\end{split}
\end{equation}

where  $o_{jk}^{(l)}(l=1,\cdots,4), a_i(i=1,\cdots,9)$(whose dimensions can be understood from \eqref{zero-bd} so we only give symbolical expressions on the righthand sides of \eqref{zero-bd}) are (vector valued) functions with the dependence of
\begin{equation*}
\begin{split}
&o^{(l)}_{jk}=o^{(l)}_{jk}(x,D\Psi,D\psi^0),\\
&a_i=a_i(D\Psi,(\vphi_-,p_-)-(\vphi_0^-,p_0^-),D(\vphi_--\vphi_0^-),D\psi_0,f,Df),
\end{split}
\end{equation*}
and, by Lemma \ref{lemma-nonl}, each $o_{jk}^{(l)}$ satisfies
\begin{equation}
\label{b}
\|o^{(l)}_{jk}(x,D\Psi,D\psi^0)\|_{1,\alp,\mcl{N}^+_f}^{(-\alp,\Gam_w)}\le C\|\psi^0\|_{2,\alp,\mcl{N}^+_f}^{(-1-\alp,\Gam_w)}\le C\sigma^{\sharp}.
\end{equation}
Moreover, it is not hard to check that $o^{(l)}_{jk}$, $a_i$ are smooth with respect to their arguments in a small neighborhood of the background state.
By the smoothness of $o^{(l)}_{jk}, a_i$, and the smallness of $o^{(l)}_{jk}$ obtained from \eqref{b} reducing $\sigma^{\sharp}$ if necessary, \eqref{zero-bd} has a unique solution $u\in C^2(\mcl{N}^+_f)\cap C^0(\ol{\mcl{N}}^+_f)$.

By subtracting the boundary problem for $u$ from the boundary problem for $\frac{\til\psi_{\eps}-\psi_0}{\eps}$, one can also check that there exists a constant $C$ so that there holds
$
\|\frac {1}{\eps}(\til\psi^{\eps}-\psi^0)-u\|_{2,\alp,\mcl{N}^+_f}^{(-1-\alp,\Gam_w)}\le C\eps,
$
or equivalently
\begin{equation}
\label{frechet-esti}
\|\til\psi^{\eps}-\psi^0-\eps u\|_{2,\alp,\mcl{N}^+_f}^{(-1-\alp,\Gam_w)}\le C\eps^2.
\end{equation}
We emphasize that the constant $C$ in \eqref{frechet-esti} is independent of the choice of $\til q$ and $\til f$, but only depends on the data in sense of Remark \ref{remark-cons}.

In addition, $o^{(l)}_{jk}, a_i$ continuously depend on their variables in the corresponding norm for each component, also the modulus of the continuity is uniform over $\mcl{B}_{\sigma^{\sharp}}^{(1)}(Id,\vphi_0^-,p_0^-,v_c)\times \mcl{B}_{C^{*}\sigma^{\sharp}}^{(2)}(r_s)$.

Back to \eqref{formal}, we have shown
\begin{equation}
\label{formal-deriv}
\begin{split}
&\phantom{=}\frac{\der}{\der\eps}|_{\eps=0}\mathfrak{J}(q_0+\eps\til q,f+\eps\til f)(x')\\
&=\til f(x')\der_r(\vphi_--\vphi_0^+-\psi_0)(f(x'),x')+(\til\vphi_--u)(f(x'),x').
\end{split}
\end{equation}
By \eqref{zero-bd}, $u$ linearly depends on $\til q$ and $\til f$. So the mapping $\mcl{L}$ defined by
$$\mcl{L}:(\til q, \til f)\mapsto\frac{\der}{\der\eps}|_{\eps=0}\mathfrak{J}(q_0+\eps\til q, f+\eps\til f)$$ is a bounded linear mapping from $C^{2,\alp}(\ol{\mcl{N}},\R^n)\times C^{3,\alp}(\ol{\mcl{N}^-_{r_s+\delta}})\times C^{2,\alp}(\ol{\mcl{N}^-_{r_s+\delta}})\times C^{1,\alp}_{(-\alp,\der\Lambda)}(\Lambda)\times C^{2,\alp}_{(-1-\alp,\der\Lambda)}(\Lambda)$ to $C^{2,\alp}_{(-1-\alp,\der\Lambda)}(\Lambda)$ where we understand $C^{2,\alp}_{(-1-\alp,\der\Lambda)}(\Lambda)$ in sense of Remark \ref{remark-holdernorms}.

It is easy to see from \eqref{formal-deriv}, \eqref{frechet-esti} that, for $(q,f)\in \mcl{B}_{\sigma^{\sharp}}^{(1)}(Id,\rho_0^-,\vphi_0^-,p_0^-,v_c)\times \mcl{B}_{C^{*}\sigma^{\sharp}}^{(2)}(r_s)$, $D\mathfrak{J}(q,f)$ is given by
\begin{equation}
\label{formula-frechet}
D\mathfrak{J}(q,f): (\til q, \til f)\mapsto \frac{\der}{\der\eps}|_{\eps=0}\mathfrak{J}(q_0+\eps\til q,f+\eps\til f),
\end{equation}
and $D\mathfrak{J}$ is continuous in $(q,f)$. This verifies Lemma \ref{lemma-frechet}(i).

Next, we prove Lemma \ref{lemma-frechet}(ii).
By \eqref{3-1-3}, \eqref{new-11}, \eqref{def-gs}, we have, for any $\psi$,
\begin{equation}
\label{homog}
F_k(x,Id,D\psi)=g_1(x,Id,0,D\psi)=g_2(x,Id,D\psi)=g_3(v_c)=0
\end{equation}

for $k=1,\cdots,n$. Lemma \ref{lemma-nonl}, \eqref{formal} and \eqref{homog} provide
$$
\frac{\der}{\der \eps}|_{\eps=0}\mathfrak{J}(Id,\vphi_0^-,p_0^-,v_c,r_s+\eps \til f)
=\frac{\der}{\der \eps}|_{\eps=0}(\vphi_0^--\vphi_0^+)(r_s+\eps \til f)=\til f\der_ r(\vphi_0^--\vphi_0^+)(r_s).$$
Thus, we obtain
\begin{equation}
\label{partial-f-deriv}
D_f\mathfrak{J}(Id,\vphi_0^-,p_0^-,v_c,r_s):\til f\mapsto \til f\der_r(\vphi_0^--\vphi_0^+)(r_s)
\end{equation}

 By \eqref{2-2-5} and \eqref{2-2-8}, we have $\der_r(\vphi_0^--\vphi_0^+)(r_s)>0$. Thus we conclude that $D_f\mathfrak{J}$ is invertible at $(Id,\vphi_0^-,p_0^-,v_c,r_s)$. Here, we choose $\delta=\frac{\sigma^\sharp}{2}$ where  $\sigma^{\sharp}>0$ depends on the data and $C^*$.
\end{proof}
\end{lemma}
\begin{proof}
[\textbf{Proof of Proposition \ref{theorem3-1}}] Let us fix $C^*=10$ then the choice of $\sigma^*$ is fixed accordingly, and moreover the choice of $\sigma_*$ only depends on the data in sense Remarkn \ref{remark-cons}.
By Lemma \ref{lemma-frechet} and the implicit mapping theorem, there exists a constant $0<\sigma_3<\delta_*$ depending on the data in sense of Remark \ref{remark-cons} so that there is a unique mapping $\mcl{S}:\mcl{B}^{(1)}_{\sigma_3}(Id,\vphi_0^-,p_0^-,v_c)\to \mcl{B}^{(2)}_{C^{*}\sigma^{*}}(r_s)$ satisfying
\begin{equation}
\label{zeros}
\mathfrak{J}(q,\mcl{S}(q))=0\quad\text{for all}\quad q\in \mcl{B}^{(1)}_{\sigma_3}(Id,\vphi_0^-,p_0^-,v_c),
\end{equation}
and $\mcl{S}$ is continuously Fr\'{e}chet differentiable in $\mcl{B}^{(1)}_{\sigma_3}(Id,\vphi_0^-,p_0^-,v_c)$.

By \eqref{2-7}, Lemma \ref{lemma-nonl}, and the observation made after the proof of Lemma \ref{lemma-nonl}, \eqref{zeros} implies
that $f^{(q)}=\mcl{S}(q)$ and $\vphi^{f^{(q)}}$ in \eqref{new-45} satisfy Proposition \ref{theorem3-1} (a)-(c).
The uniqueness of a solution of Problem \ref{problem2-1} also follows from the implicit mapping theorem. We note that the choice of $\sigma_3>0$ only depends on the data because $\sigma_3$ depends on the data and $\delta_*$ which depends on the data and $C^*$ but we fixed $C^*$ as $C^*=10$.
\end{proof}
%%%%%%%%%%%%%%%%%%%%%%%%%%%%%%%%%%%%%%%%%%%%%%%%%%%%%%%%%%%%%%%%%%%%%
%%%%%%%%%%%%%%%%%%%%%%%%%%%%%(Section 4)%%%%%%%%%%%%%%%%%%%%%%%%%%%%%
%%%%%%%%%%%%%%%%%%%%%%%%%%%%%%%%%%%%%%%%%%%%%%%%%%%%%%%%%%%%%%%%%%%%%
\section{Transport Equation for Pressure $p$}
\label{sec:nozzle3}
This section is devoted to  step 2 in section \ref{subsec-framework}.
\subsection{The method of characteristics}
\label{subsec-char}
By the choice of $\sigma_3$ in Proposition \ref{theorem3-1}, if $\vphi$ is a solution to Problem \ref{problem2-1}, then $\der_r\vphi$ is bounded below by a positive constant in $\mcl{N}^+$. So $r$ can be regarded as a time-like variable and then we apply the method of characteristics to solve \eqref{2-p} for $p$. According to \eqref{3-0-19}, however, $D\vphi$ is not Lipschitz continuous up to $\der\mcl{N}^+$. So the unique solvability of \eqref{2-p} needs to be established first, and we will prove the following proposition.

\begin{proposition}
\label{proposition4-0-1} Let $\sigma_3$ be as in Proposition \ref{theorem3-1}. For any given $\alp\in(0,1)$ there exist constant $C>0$ and $\sigma'\in(0,\sigma_3]$ depending on the data in sense of Remark \ref{remark-cons} so that if $\varsigma_1+\varsigma_2+\varsigma_4\le \sigma'$, and $\vphi$ is a solution of Problem \ref{problem2-1} with a shock $S=\{r=f(x')\}$, then \eqref{2-p} with the initial condition \eqref{2-p-bc} on $S$ has a unique solution $p$ satisfying
\begin{equation}
\label{4-0-11}
\|p-p_0^+\|_{1,\alp,\mcl{N}^+}^{(-\alp,\Gam_w)}\le C(\varsigma_1+\varsigma_2+\varsigma_4).
\end{equation}
\end{proposition}
We denote as $\mcl{N}^+$ for the subsonic domain $\{|D\Psi^{-1}D\vphi|<c\}=\{r>f(x')\}$.
We first note that $\vphi$ is in $C^2$ up to the wall $\Gam_w$ of $\mcl{N}^+$ away from the corners.
\begin{lemma}
\label{lemma-reg}
Given $(\Psi,\vphi_-,p_-,v_{ex})\in \mcl{B}^{(1)}_{\sigma_3}(Id,\vphi_0^-,p_0^-,v_c)$ satisfying Proposition \ref{theorem3-1} (i), (ii) reducing $\sigma_3$ if necessary, $\vphi$, the corresponding solution of Problem \ref{problem2-1}, with a shock $S$ satisfies
\begin{equation}
\label{4-0-7}
\|\vphi-\vphi_0^+\|_{2,\alp,\mcl{N}^+}^{(-1-\alp,\ol{S}\cup \ol{\Gam_{ex}})}\le C(\varsigma_1+\varsigma_2+\varsigma_4)
\end{equation}
for $\varsigma_1, \varsigma_2, \varsigma_4$ in Theorem \ref{main-thm'} and Proposition \ref{theorem3-1}.
\end{lemma}
To prove Lemma \ref{proposition4-0-1}, it suffices to show that the nonlinear boundary problem \eqref{ift-eqn} in the fixed domain $\mcl{N}^+$ has a unique solution $\psi$ so that $\vphi=\vphi_0^++\psi$ satisfies \eqref{4-0-7}, and this can be proven by a similar method to the proof of Lemma \ref{lemma-nonl} except that the scaling argument in step 4 needs to be performed for any point $y_0\in \mcl{N}^+\setminus(\ol{S}\cup\ol{\Gam_{ex}})$. So we skip the proof of Lemma \ref{proposition4-0-1}.

Considering $r$ as a time-like variable, it is more convenient to use a changes of variables to transform $\mcl{N}^+$ to $\mcl{N}^+_{r_s}$ keeping \eqref{4-0-7} in the transformed domain $\mcl{N}^+_{r_s}$. This requires a new change of variables:
\begin{equation}
\label{3-2-20}
\begin{split}
G_{\vphi}(r,x')&=\bigl((k(\vphi_--\vphi)(r,x')+r_s)(1-\chi(r))+r\chi(r),x' \bigr)\\
&=:(\upsilon ^{(\vphi)}(r,x'),x')=:(\til r,x')
\end{split}
\end{equation}
with a smooth function $\chi(r)$ satisfying
$\chi(r)=\begin{cases}
0 &\text{if}\;r\le r_s+\frac{r_1-r_s}{10}\\
1 &\text{if}\;r\ge r_1-\frac{r_1-r_s}{2}
\end{cases}$ and $\chi'(r)\ge 0$ for $r>0$. By \eqref{4-0-7}, if $\varsigma_1+\varsigma_2+\varsigma_4=\sigma$ is small depending on the data in sense of Remark \ref{remark-cons}, and we choose a constant $k$ satisfying
\begin{equation}
\label{3-2-10}
\begin{split}
&C\sigma \le \frac{1}{10}\min\{(\vphi_0^--\vphi_0^+)(r_1),\der_r(\vphi_0^--\vphi_0^+)(r_s)\},\\
& k=\frac{r_1-r_s}{8(\vphi_0^--\vphi_0^+)(r_1)},
\end{split}
\end{equation}
then we have
\begin{equation}
\label{Gg}
\der_rv^{(\vphi)}\ge \frac{r_1-r_s}{4}>0\quad\text{in}\quad \mcl{N}^+_{r_s-2\delta},
\end{equation}
thus $G_{\vphi}$ is invertible.

We choose $x'$ in \eqref{spherical-coord} so that $(\til r,x')$ forms an orthogonal coordinate system. To specify this choice, we write $(r,x')$ as $(r,\vartheta)$ with $\vartheta=(\vartheta_1,\cdots,\vartheta_{n-1})$ for the rest of paper. And then, we apply the change of variables \eqref{3-2-20} to rewrite \eqref{2-p} in $\mcl{N}^+$ as
\begin{equation}
\label{2-p-n}
V\cdot DE
=(V\cdot\hat {\til r})\der_{\til r}E+\sum_{j=1}^{n-1}(V\cdot\hat{\vartheta}_j)\der_{\vartheta_j}E=0
\;\;\text{for}\;\;(\til r,\vartheta)\in (r_s,r_1]\times \Lambda
\end{equation}
with
\begin{equation}
\label{a}
V=J_{\vphi}(D\Psi^{-1})^T(D\Psi^{-1})D\vphi,\quad
E\circ G_{\vphi}=\frac{p}{(B_0-\frac 12|D\Psi^{-1}D\vphi|^2)^{\frac{\gam}{\gam-1}}}
\end{equation}
where $J_{\vphi}$ is the Jacobian matrix associated with the change of variables in \eqref{3-2-20}. Here, $\hat {\til r}, \hat{\vartheta}_j$ are unit vectors pointing the positive directions of $\til r$ and $\vartheta_j$ respectively. By Proposition \ref{theorem3-1}, Lemma \ref{lemma-reg} and \eqref{3-2-20}, we have
\begin{lemma}
\label{lemma-V}
Set $V_0:=J_{\vphi_0^+}D\vphi_0^+$. Then, $V$ defined in \eqref{a} satisfies
\begin{equation*}
\|V-V_0\|_{1,\alp,\mcl{N}^+_{r_s}}^{(-\alp,\Gam_w)}+\|V-V_0\|_{1,\alp,\mcl{N}^+_{r_s}}^{(-\alp,\ol{S_0}\cup\ol{\Gam_{ex}})}\le C(\varsigma_1+\varsigma_2+\varsigma_4)
\end{equation*}
for $\varsigma_1, \varsigma_2, \varsigma_4$ in Theorem \ref{main-thm'} and Proposition \ref{theorem3-1}. Thus, reducing $\sigma_3$ in Proposition \ref{theorem3-1} if necessary, there is a constant $\omega_0>0$  depending on the data in sense of Remark \ref{remark-cons} so that there holds
$
V\cdot\hat{\til r}\ge \omega_0>0
$ in $\ol{\mcl{N}^+_{r_s}}$.
\end{lemma}
The proof is obtained by a straightforward computation so we skip here.
\begin{lemma}
\label{lemma-ent-0}
Let us set $E_0^+:=\frac{p_0^+}{(B_0-\frac 12|\der_r\vphi_0^+|^2)^{\frac{\gam}{\gam-1}}}$, then $E_0^+$ is a constant.
\begin{proof}
For $V_0$ defined in Lemma \ref{lemma-V}, by \eqref{2-p-nn}, $E_0^+$ satisfies
\begin{equation*}
V_0\cdot DE_0^+=(V_0\cdot\hat r)\der_rE_0^+=0.
\end{equation*}
Since $V_0\cdot \hat r>0$ for all $r\in[r_s,r_1]$, and $E_0^+(r_s)$ is a constant by \eqref{2-2-5}, \eqref{2-2-6}, therefore we conclude $E_0^+(r)=E_0^+(r_s)=constant$ for all $r\in[r_s,r_1]$.
\end{proof}
\end{lemma}
By Lemma \ref{lemma-V}, we can divide the equation \eqref{2-p-n} by $V\cdot\hat {\til r}$ so that the radial speed of all characteristics associated with $V$ is 1. From now on, we consider the equation
\begin{equation}
\label{2-p-nn}
\der_{\til r}E+\sum_{j=1}^{n-1}\frac{V\cdot\hat{\vartheta}_j}{V\cdot\hat {\til r}}\der_{\vartheta_j}E=0\;\;\text{in}\;\;\mcl{N}^+_{r_s}.
\end{equation}
We denote
\begin{equation}
\label{c}
W^*=(1,\frac{V\cdot\hat\vartheta_1}{V\cdot\hat {\til r}},\cdots,\frac{V\cdot\hat\vartheta_{n-1}}{V\cdot\hat {\til r}}).
\end{equation}

The initial condition for $E$ is given from (\ref{2-p-bc}). Setting $Q_2\!:=\!(D\Psi^{-1})^T\!D\Psi^{-1}$, $E(r,\vartheta)$ must satisfy $E(r_s,\vartheta)=E_{int}(\vartheta)$ for all $\vartheta\in\Lambda$ where we set
\begin{equation}
\label{4-0-5}
E_{int}(\vartheta):=
\frac{\rho_-(\frac{|\grad(\vphi_--\vphi)|Q_2}{|D\Psi^{-1}\grad(\vphi_--\vphi)|}\grad\vphi_-\cdot\nu_s)^2+p_--\rho_-K_s}{(B_0-\frac 12|D\Psi^{-1}(\Psi)\grad\vphi|^2)^{\frac{\gam}{\gam-1}}}\mid_{(r,\vartheta)=(f(\vartheta),\vartheta)}.
\end{equation}
Above, $f(\vartheta)$ is the location of the transonic shock for the potential $\vphi$ in Proposition \ref{theorem3-1}.
To simplify notation, we write $\til r$ as $r$ hereafter. By \eqref{2-p-nn} and \eqref{4-0-5}, Proposition \ref{proposition4-0-1} is a direct result of the following lemma.
\begin{lemma}
\label{lemma4-0-1}
Fix $\alp\in(0,1)$, and suppose that a vector valued function $W=(1,W_2,\cdots,W_{n})$ satisfies $W\cdot\nu_w=0$ on $\Gam_{w,r_s}=\der\mcl{N}^+_{r_s}\cap\Gam_w$ for a unit normal $\nu_w$ on $\Gam_w^+$, and
\begin{equation}
\label{est2}
\|W\|_{1,\alp,\mcl{N}^+_{r_s}}^{(-\alp,\Gam_w)}+\|W\|_{1,\alp,\mcl{N}^+_{r_s}}^{(-\alp,\ol{S_{r_s}}\cup\ol{\Gam_{ex}})}\le k
\end{equation}
 for some constant $k>0$. Then the initial value problem
 \begin{equation}
 \label{IVP}
 \der_r F+\sum_{j=1}^{n-1}W_{j+1}\der_{\vartheta_j}F=0\;\;\text{in}\;\;\mcl{N}^+_{r_s},\quad
 F(r_s,\vartheta)=F_0(\vartheta)\;\;\text{for}\;\;\vartheta\in\Lambda
 \end{equation}
 has a unique solution $F$ satisfying
 \begin{equation}
\label{4-0-8}
\|F\|_{1,\alp,\mcl{N}^+_{r_s}}^{(-\alp,\Gam_w)}\le C\|F_0\|_{1,\alp,\Lambda}^{(-\alp,\der \Lambda)}
\end{equation}
for a constant $C>0$ depending on $n,\Lambda,r_s,r_1,\alp$ and $k$.
\end{lemma}

\begin{proof}
[\textbf{Proof of Proposition \ref{proposition4-0-1}}]\quad\\
By \eqref{3-0-6} and Lemma \ref{lemma-V}, $W^*=(1,\frac{V\cdot\hat{\vartheta}_1}{V\cdot\hat r},\cdots,\frac{V\cdot\hat{\vartheta}_{n-1}}{V\cdot\hat r} )$ in \eqref{2-p-nn} satisfies all the assumptions of Lemma \ref{lemma4-0-1}. So, \eqref{2-p-nn} with \eqref{4-0-5} has a unique solution $E$. Then the pressure function $p$ is given by
$$
p=(B_0-\frac 12|D\Psi^{-1}D\vphi|^2)^{\frac{\gam}{\gam-1}}\cdot(E\circ G_{\vphi}).
$$
It remains to verify (\ref{4-0-11}). By Lemma \ref{lemma-ent-0}, if $E$ satisfies \eqref{2-p-nn} then so does $E-E_0^+$. Then, by \eqref{4-0-8}, we obtain
\begin{equation}
\label{new-13}
\|E-E_0^+\|_{1,\alp,\mcl{N}^+_{r_s}}^{(-\alp,\Gam_w)}\le C\|E_{int}-E_0^+\|_{1,\alp,\Lambda}^{(-\alp,\der \Lambda)},
\end{equation}
where $E_{int}$ is given by \eqref{4-0-5}.

Moreover, by \eqref{4-0-5}, one can explicitly show
\begin{equation}
\label{new-14}
\|E_{int}-E_0^+\|_{1,\alp,\Lambda}^{(-\alp,\der \Lambda)}\le C(\varsigma_1+\varsigma_2+\varsigma_4).
\end{equation}
This with \eqref{3-0-19} gives \eqref{4-0-11}. We choose $\sigma'=\sigma_3$ by reducing $\sigma_3$ if necessary to satisfy \eqref{4-0-7} and \eqref{3-2-10}.
\end{proof}
Now it remains to prove Lemma \ref{lemma4-0-1}. Consider the following ODE system
\begin{equation}
\label{4-0-9'}
\begin{split}
&\dot X(t)=-W(X(t))\;\;\text{for}\;\;r<t\le 2r-r_s\\
&X(r)=(r,\vtheta)
\end{split}
\end{equation}
for $(r,\vartheta)\in \mcl{N}^+_{r_s}\cup\Gam_{ex}$.
The slip boundary condition $W\cdot\nu_w=0$ on $\Gam_w^+$ implies that any characteristic curve $X$ does not penetrate out of $\mcl{N}^+_{r_s}$.
 If (\ref{4-0-9'}) has a unique solution $X(t;r,\vartheta)$ for any given $(r,\vartheta)\in \mcl{N}^+_{r_s}\cup\Gam_{ex}$, then we can define a mapping $\mcl{J}:\mcl{N}^+_{r_s}\cup\Gam_{ex}\to \Lambda$ by
\begin{equation}
\label{4-2-1}
\mcl{J}:(r,\vartheta)\mapsto X(2r-r_s;r,\vtheta).
\end{equation}
If $\mcl{J}$ is differentiable then $F=F_{0}\circ\mcl{J}$ is a solution to (\ref{IVP}).
\subsection{The system of ODEs \eqref{4-0-9'}}
\eqref{4-0-9'} has a solution for any given $(r,\vartheta)\in(r_s,r_1]\times \Lambda$ by the Cauchy-Peano theorem \cite[7.3]{A}. But the standard ODE theory is not sufficient to claim the uniqueness of a solution for \eqref{4-0-9'} because $W$ is only in $C^{\alp}(\ol{\mcl{N}^+_{r_s}})$. To prove the uniqueness, we will use the fact that $DW$ is integrable along each solution $X$ to $\eqref{4-0-9'}$. Denote as $X(t;r,\vartheta)$ for a solution to \eqref{4-0-9'}.
Hereafter, we call each $X(t;r,\vartheta)$ a \emph{characteristic} associated with $W$ initiated from $(r,\vartheta)$.
\begin{remark}
\label{remark-int}
Writing $X=(X_1,\cdots,X_n)$, it is easy to check
\begin{equation}
\label{first-comp}
X_1(t;r,\vartheta)=2r-t
\end{equation}
for any $(r,\vartheta)\in (r_s,r_1]\times \Lambda$. So, we restrict \eqref{4-0-9'} for $t\in[2r-r_1,2r-r_s]$ because of $X_1(2r-r_1;r,\vartheta)=r_1$, $X_1(2r-r_s;r,\vartheta)=r_s$.
\end{remark}
\begin{lemma}
\label{lemma4-0-2}
For any given $(r,\vtheta)\in(r_s,r_1]\times \Lambda$, \eqref{4-0-9'} has a unique solution. Moreover there is a constant $C>0$ depending on $\Lambda,n$ and $k$ satisfying
\begin{equation}
\label{4-1-1}
\frac 1C  dist((r,\vtheta),\Gam_w)\le dist(X(t;r,\vtheta),\Gam_w)\le C dist((r,\vtheta),\Gam_w)
\end{equation}
for $t\in[2r-r_1,2r-r_s].$
\begin{proof}
Since $W$ is in $C^{\alp}(\ol{\mcl{N}^+_{r_s}})$, (\ref{4-0-9'}) has a solution $X(t;r,\vtheta)$ for any $(r,\vtheta)\in (r_s,r_1]\times \Lambda$. Fix $r\in(r_s,r_1]$. Let $X^{(1)}$ and $X^{(2)}$ be solutions to (\ref{4-0-9'}) with
\begin{equation}
\label{bb}
X^{(1)}(r)=(r,\vtheta_1),\quad X^{(2)}(r)=(r,\vtheta_2)
\end{equation}
for $\vtheta_1,\vtheta_2\in\Lambda$. By \eqref{4-0-9'}, we have
\begin{equation*}
\begin{split}
%\label{4-1-2}
\frac{d}{dt}|X^{(1)}-X^{(2)}|^2&=2(W(X^{(2)})-W(X^{(1)}))\cdot(X^{(1)}-X^{(2)})\\
&\le C|DW(2r-t,\cdot)|_{L^{\infty}(\Lambda)}|X^{(1)}-X^{(2)}|^2.
\end{split}
\end{equation*}
 In $(r,\vartheta)-$coordinates, it is easy to see
\begin{equation*}
dist((\til{r},\til{\vtheta}),\ol{S_0}\cup\ol{\Gam_{ex}})=\min\{r_1-\til{r},\til{r}-r_s\}
\quad\text{for}\;\;(\til{r},\til{\vtheta})\in \mcl{N}^+_{r_s}.
\end{equation*}
So there is a constant $C$ depending on $n, \mcl{N}^+_{r_s}$ so that $DW$ satisfies
\begin{equation}
\label{4-1-3}
|DW(2r-t,\vartheta)|\le Ck(1+(2r-r_s-t)^{-1+\alp}+(r_1-2r+t)^{-1+\alp})=:m(t)
\end{equation}
for any $\vartheta\in\Lambda$ and $t\in[r,2r-r_s]$ where $k$ is as in \eqref{est2}. Then we apply the method of the integrating factor to obtain
\begin{equation}
\label{4-1-4}
\frac{d}{dt}(e^{-\int_{r}^tm(t')dt'}|X^{(1)}-X^{(2)}|^2(t))\le 0
\quad\text{for}\;\;t\in [r,2r-r_s].
\end{equation}

If $\vtheta_1=\vtheta_2$ then \eqref{bb} and \eqref{4-1-4} imply $X^{(1)}(t)=X^{(2)}(t)$ for all $t\in [r,2r-r_s]$. We have accomplished the unique solvability of (\ref{4-0-9'}).

Since $W\cdot \nu_w=0$ on $\Gam^+_w$, any characteristic, initiated from a point on the wall of $\mcl{N}^+_{r_s}$, remains on the wall i.e.,
if $\vtheta_2\in \der\Lambda$, then $X^{(2)}$ lies on $\Gam_w$. Then (\ref{4-1-4}) implies
\begin{equation*}
\begin{split}
dist(X^{(1)}(t),\Gam_w)&=\inf_{\vtheta_2\in \der\Lambda}|X^{(1)}(t)-X^{(2)}(t)|\\
&\le C\inf_{\vtheta_2\in\der\Lambda}|(r,\vtheta_1)-(r,\vtheta_2)|=C dist((r,\vtheta_1),\Gam_w).
\end{split}
\end{equation*}
This verifies the second inequality of (\ref{4-1-1}). One can argue similarly to verify the first inequality of (\ref{4-1-1}).
\end{proof}
\end{lemma}
\subsection{Regularity of $\mcl{J}:(r,\vartheta)\mapsto X(2r-r_s;r,\vartheta)$ in \eqref{4-2-1}}
\begin{lemma}
\label{lemma4-2-1}
Fix $\alp\in(0,1)$, and let $k$ be as in Lemma \ref{lemma4-0-1}. The mapping $\mcl{J}$ defined in \eqref{4-2-1} is continuously differentiable, and it satisfies
\begin{equation*}
%\label{4-2-2}
\|\mcl{J}\|_{1,\alp,\mcl{N}^+_{r_s}}^{(-1,\Gam_w)}\le C
%|\mcl{J}|_{1,\mcl{N}^+_{r_s}}+[D\mcl{J}]_{\alp;\mcl{N}^+_{r_s}}^{(0),\Gam_w}\le C
\end{equation*}
for a constant $C>0$ depending on $n,\Lambda,r_s,r_1,\alp$ and $k$.
\end{lemma}
\begin{proof}
First of all, let us assume that $\mcl{J}$ satisfies
\begin{equation}
\label{est1}
|\mcl{J}|_{1,0,\mcl{N}^+_{r_s}}\le C,
\end{equation}
and we estimate $[D\mcl{J}]_{\alp,\mcl{N}^+_{r_s}}^{(0,\Gam_w)}$. If $\mcl{J}$ is differentiable then (\ref{4-0-9'}) implies
\begin{equation}
\label{diff1}
\begin{split}
&\text{(a)}\;\;\der_r\mcl{J}(r,\vartheta)=Z(2r-r_s;r,\vartheta)-2W(X(2r-r_s;r,\vartheta)),\\
&\text{(b)}\;\;\der_{\vartheta}\mcl{J}=Z(2r-r_s;r,\vartheta)
\end{split}
\end{equation}
where $Z(t;r,\vartheta)$ is the solution to the ODE system
\begin{equation}
\label{IVP2}
\begin{split}
&\dot{Z}(t;r,\vartheta)=-DW(X(t;r,\vartheta))Z(t;r,\vartheta),\\
&Z(r;r,\vartheta)=
\begin{cases}
(1,0,\cdots,0)+W(r,\vartheta)&\text{for  (a)}\\
(0,\delta_{1j},\cdots,\delta_{n-1,j})&\text{for  (b).}
\end{cases}
\end{split}
\end{equation}
Fix $(r,\vartheta),(r',\vartheta')\in(r_s,r_1]\times \Lambda$ with $(r,\vartheta)\neq (r',\vartheta')$, and consider $J:=\frac{|D\mcl J(r,\vartheta)-D\mcl{J}(r',\vartheta'))|}{|(r,\vartheta)-(r',\vartheta')|^{\alp}}$. For simplicity, we assume $r_s<r<r'<r_1$ and $\vartheta=\vartheta'$. Under this assumption, we first estimate $\frac{|\mcl{J}_r(r,\vartheta)-\mcl{J}_r(r',\vartheta)|}{|r-r'|^{\alp}}$ because the case of $\vtheta\neq\vtheta'$ can be handled similarly. Also, the estimate of $\frac{|\mcl{J}_{\vtheta}(r,\vtheta)-\mcl{J}_{\vtheta}(r',\vtheta')|}{|(r,\vtheta)-(r',\vtheta')|^{\alp}}$ is similar and even simpler then the estimate of $\frac{|\mcl{J}_r(r,\vtheta)-\mcl{J}_r(r',\vtheta)|}{|r-r'|^{\alp}}$.

By (\ref{est2}) and (\ref{est1}), it is easy to check
\begin{equation}
\label{4-2-15}
\frac{2|W(X(2r-r_s;r,\vtheta))-W(X(2r'-r_s;r',\vtheta))|}{|r-r'|^{\alp}}\le C
\end{equation}
for a constant $C$ depending on $n,\Lambda,r_s,r_1,\alp$ and $k$. For the rest of the proof, a constant $C$ depends on
$n,\Lambda,r_s,r_1,\alp$ and $k$ unless otherwise specified.

To estimate $\frac{|Z(2r-r_s;r,\vtheta)-Z(2r'-r_s;r',\vtheta)|}{|r-r'|^{\alp}}$, we split
the numerator as
\begin{equation*}
\begin{split}
&\phantom{\le}|Z(2r-r_s;r,\vtheta)-Z(2r'-r_s;r',\vtheta)|\\
&\le |Z(2r-r_s;r,\vtheta)-Z(2r-r_s;r',\vtheta)|+|Z(2r-r_s;r',\vtheta)-Z(2r'-r_s;r',\vtheta)|,
\end{split}
\end{equation*}
and consider them separately.

By (\ref{4-1-3}) and (\ref{IVP2}), one can easily check
\begin{equation}
\label{est4}
|Z(2r-r_s;r',\vtheta)-Z(2r'-r_s;r',\vtheta)|\le C(r'-r)^{\alp}.
\end{equation}
Let us denote as $Z_1(t)$ and $Z_2(t)$ for $Z(t;r,\vtheta)$ and $Z(t;r',\vtheta)$ respectively. Then, by (\ref{IVP2}), $z:=|Z_1-Z_2|^2$ satisfies
\begin{equation}
\label{4-2-20}
\dot{z}
\le C(|DW(X(t;r,\vtheta))|z+|DW(X(t;r,\vtheta))-DW(X(t;r',\vtheta))|z^{1/2})
\end{equation}
on the interval $I_{r,r'}:=[2r'-r_1,2r'-r_s]\cap[2r-r_1,2r-r_s]$. By (\ref{est1}), we have
\begin{equation*}
\begin{split}
&\phantom{\le}|DW(X(t;r,\vtheta))-DW(X(t;r',\vtheta))|\\
&\le C|r-r'|^{\alp}\frac{|DW(X(t;r,\vtheta))-DW(X(t;r',\vtheta))|}{|X(t;r,\vtheta)-X(t;r',\vtheta)|^{\alp}}.
\end{split}
\end{equation*}
Set $d:=\min\{dist((r,\vtheta),\Gam_w),dist((r',\vtheta),\Gam_w)\}$ and
\begin{align*}
\mcl{Z}(t):=\frac{|DW(X(t;r,\vtheta))-DW(X(t;r',\vtheta))|}{|X(t;r,\vtheta)-X(t;r',\vtheta)|^{\alp}}.
 \end{align*}

 For a constant $\beta\in(0,1)$, we estimate $\mcl{Z}^{\beta}$ and $\mcl{Z}^{1-\beta}$ using two different weighted H\"{o}lder norms of $W$. By (\ref{4-1-1}), for any $t\in I_{r,r'}$, we have
\begin{equation}
\label{4-2-18}
\mcl{Z}^{\beta}\le Cd^{-\beta}(\|W\|_{1,\alp,\mcl{N}^+_{r_s}}^{(-\alp,\Gam_w)})^{\beta}\le Ck^{\beta}d^{-\beta}.
\end{equation}
Next, set $d_*(t):=\min\{dist(X(t;r,\vtheta),\ol{S_{0}}\cup\ol{\Gam}_{ex}),dist(X(t;r',\vtheta),\ol{S_{0}}\cup\ol{\Gam_{ex}})\}$. Then $d_*(t)$ can be expressed as
\begin{equation*}
d_*(t)=\min\{2r-t-r_s,r_1-2r'+t\}\;\;\text{for}\;\;t\in I_{r,r'}.
\end{equation*}
So $(d_*(t))^{-1+\beta}$ is bounded by $m_*(t)$ which is integrable over $I_{r,r'}$, then we have
\begin{equation}
\label{4-2-19}
\mcl{Z}^{1-\beta}(t)
\le Cm_*(t)(\|W\|_{1,\alp,\mcl{N}^+_{r_s}}^{(-\alp,\ol{S_0}\cup\ol{\Gam}_{ex})})^{1-\beta}\le Ck^{1-\beta}m_*(t).
\end{equation}
Back to (\ref{4-1-3}), by (\ref{4-1-4}), (\ref{IVP2}), (\ref{4-2-18}), (\ref{4-2-19}) and the method of integrating factor, we obtain
\begin{equation*}
\begin{split}
&\frac{d}{dt}(e^{-\int_r^tm(t')dt'} z)\le (r'-r)^{\alp}Ckd^{-\beta}m_*z^{1/2}\\
&\Rightarrow z(t)\le C((r'-r)^2+(r'-r)^{\alp}d^{-\beta}\sup_{t\in I_{r,r'}}z^{1/2})\;\;\text{for all}\;t\in I_{r,r'}\\
&\Rightarrow \sup_{t\in I_{r,r'}}z\le C((r'-r)^2+d^{-2\beta}(r'-r)^{2\alp})
\end{split}
\end{equation*}
which provides
\begin{equation}
\label{est3}
\frac{|Z(2r-r_s;r,\vartheta)-Z(2r-r_s;r',\vartheta)|}{|r'-r|^{\alp}}\le Cd^{-\beta}
\end{equation}
for any $\beta\in(0,1)$ and a constant $C$ depending on $n, \mcl{N}^+_{r_s}, k$ and $\beta$. The case of $\vartheta\neq \vartheta'$ can be handled in the same fashion. The estimate of $\frac{|\der_{\vartheta}\mcl{J}(r,\vartheta)-\der_{\vartheta}\mcl{J}(r',\vartheta')|}{|(r,\vartheta)-(r',\vartheta')|^{\alp}}$ is even simpler then the estimate of $\frac{|\der_{r}\mcl{J}(r,\vartheta)-\der_{r}\mcl{J}(r',\vartheta')|}{|(r,\vartheta)-(r',\vartheta')|^{\alp}}$ so we skip the details and we conclude that, for any constant $\beta\in(0,1)$, there is a constant $C$ so that there holds
\begin{equation*}
[D\mcl{J}]_{\alp,\mcl{N}^+_{r_S}}^{(\beta-\alp,\Gam_w)}\le C.
\end{equation*}
By means of similar arguments, one can directly prove $\mcl{J}\in C^1(\ol{\mcl{N}^+_{r_s}})$.
\end{proof}
\begin{remark}
 $\mcl{J}$ is in $C^{1,\alp}_{(-1-\alp+\beta,\Gam_w)}(\mcl{N}^+_{r_s})$ for any $\beta\in(0,1)$.
\end{remark}
Now, we are ready to prove Lemma \ref{lemma4-0-1}.
\begin{proof}
[\textbf{Proof of Lemma \ref{lemma4-0-1}}] By Lemma \ref{lemma4-2-1}, $F:=F_{0}\circ\mcl{J}$ is a solution to (\ref{IVP}). By (\ref{4-1-1}), we easily get the estimate
$
\|F\|_{1,0,\mcl{N}^+_{r_s}}^{(-\alp,\Gam_w)}\le C\|F_0\|_{1,\alp,\Lambda}^{(-\alp,\der \Lambda)}
$
so it is remains to estimate $[F]_{\alp,\mcl{N}^+_{r_s}}^{(0,\Gam_w)}$.
Fix $(r,\vtheta)$ and $(r',\vtheta')$ in $(r_s,r_1]\times \Lambda$, and assume $(r,\vtheta)\neq(r',\vtheta')$. Let us set
$$
d:=\frac 12\min\{dist((r,\vtheta),\Gam_w),dist((r',\vtheta'),\Gam_w)\}>0.
$$ Setting $[\zeta]^b_a:=\zeta(b)-\zeta(a)$, we have
\begin{align*}
&|DF(r,\vtheta)-DF(r',\vtheta')|\\
&\le |[D\mcl{J}]^{(r,\vtheta)}_{(r',\vtheta')}||DF_{0}(\mcl J(r,\vtheta))|
+|D\mcl{J}(r',\vtheta')||[DF_{0}\circ \mcl{J}]^{(r,\vartheta)}_{(r',\vtheta')}|.
\end{align*}
By (\ref{4-1-1}) and Lemma \ref{lemma4-2-1}, we obtain
\begin{align*}
&|[D\mcl{J}]^{(r,\vtheta)}_{(r',\vtheta')}||DF_{0}(\mcl J(r,\vtheta))|\le C\;d^{-\alp}[D\mcl{J}]_{\alp,\mcl{N}^+_{r_s}}^{(0,\Gam_w)}\;d^{-1+\alp}K_{F_0},\\
&|D\mcl{J}(r',\vtheta')||[DF_{0}\circ \mcl{J}]^{(r,\vartheta)}_{(r',\vtheta')}|\le C\;d^{-1}K_{\mcl{J}}\;K_{F_0}(K_{\mcl{J}}|(r,\vtheta)-(r',\vtheta')|)^{\alp}
\end{align*}
with $K_{\mcl{J}}=|\mcl{J}|_{1,\mcl{N}^+_{r_s}}$ and $K_{F_0}=\|F_{0}\|_{1,\alp,\Lambda}^{(-\alp,\der \Lambda)}$.
Thus $F=F_{0}\circ \mcl{J}$ satisfies (\ref{4-0-8}).

The uniqueness of a solution can be checked easily. If $F^{(1)}$ and $F^{(2)}$ are solutions to (\ref{IVP}) then $F^{(1)}-F^{(2)}$ is identically  $0$ along every characteristic $X(t;r,\vartheta)$ associated with $W$, and so $F^{(1)}=F^{(2)}$ holds in $\mcl{N}^+_{r_s}$.
\end{proof}

%%%%%%%%%%%%%%%%%%%%%%%%%%%%%%%%%%%%%%%%%%%%%%%%%%%%%%%%%%%%%%%%%%%%%
%%%%%%%%%%%%%%%%%%%%%%%%%%%%%(Section 5)%%%%%%%%%%%%%%%%%%%%%%%%%%%%%
%%%%%%%%%%%%%%%%%%%%%%%%%%%%%%%%%%%%%%%%%%%%%%%%%%%%%%%%%%%%%%%%%%%%%

\section{Proof of Theorem \ref{main-thm'}:Existence}
\label{sec:nozzle4}
To prove Theorem \ref{main-thm'}, we prove a weak implicit mapping theorem for infinite dimensional Banach spaces.
\subsection{Weak implicit mapping theorem}

Fix a constant $R>0$ and set
\begin{align}
\label{gf}
&\mcl{B}_{R}^{(2)}(h_0):=\{h\in C^{1,\alp}_{(-\alp,\der\Lambda)}(\Lambda):\|h-h_0\|_{1,\alp,\Lambda}^{(-\alp,\der\Lambda)}\le R\}.
\end{align}
For $\mcl{B}_{\sigma}^{(1)}(Id,\vphi_0^-,p_0^-,v_c)$ in \eqref{sets}, define $\mcl{P}:\mcl{B}_{\sigma}^{(1)}(Id,\vphi_0^-,p_0^-,v_c)\to \mcl{B}_{C\sigma}^{(2)}(p_c)$ by
\begin{equation}
\label{def-p}
\mcl{P}:(\Psi,\vphi_-,p_-,v_{ex})\mapsto p|_{\Gam_{ex}}
\end{equation}
where $p$ is the solution to (\ref{2-p}) with $\vphi$ in Proposition \ref{theorem3-1} for the given $(\Psi,\vphi_-,p_-,v_{ex})$. By Proposition \ref{theorem3-1} and \ref{proposition4-0-1}, there exists positive constants $\sigma_3$ and $C$ so that whenever $0<\sigma\le \sigma_3$, $\mcl{P}$, defined in \eqref{def-p}, is well defined. To prove theorem \ref{main-thm'}, we need to show $\mcl{P}^{-1}(p_{ex})\neq \emptyset$ for any $p_{ex}\in \mcl{B}_{\til{\sigma}}^{(2)}(p_c)$ for a sufficiently small constant $\til{\sigma}>0$. So, we need the following lemma.
\begin{proposition}[Weak Implicit Mapping Theorem]
\label{lemma5-1}
Let $\mathfrak{C}_1$ and $\mathfrak{C}_2$ be Banach spaces compactly imbedded in Banach spaces $\mathfrak{B}_1$ and $\mathfrak{B}_2$ respectively. Also, suppose that  there are two Banach spaces $\mathfrak{C}_3, \mathfrak{B}_3$ with $\mathfrak{C}_3\subset\mathfrak{B}_3$.  For a point $(x_0,y_0)\in\mathfrak{C}_1\times \mathfrak{C}_2$, suppose that a mapping $\mcl{F}$ satisfies the followings:
\begin{itemize}
\item[(i)] $\mcl{F}$ maps a small neighborhood of $(x_0,y_0)$ in $\mathfrak{B}_1\times\mathfrak{B}_2$ to $\mathfrak{B}_3$, and maps a neighborhood of $(x_0,y_0)$ in $\mathfrak{C}_1\times \mathfrak{C}_2$ to $\mathfrak{C}_3$ with $\mcl{F}(x_0,y_0)={\bf{0}}$,
\item[(ii)] whenever a sequence $\{(x_k,y_k)\}\subset \mathfrak{C}_1\times \mathfrak{C}_2$ near $(x_0,y_0)$ converges to $(x_*,y_*)$ in $\mathfrak{B}_1\times\mathfrak{B}_2$, the sequence $\{\mcl{F}(x_k,y_k)\}\!\subset\!\mathfrak{C}_3$ converges to $\mcl{F}(x_*,y_*)$ in $\mathfrak{B}_3$,
\item[(iii)] $\mcl{F}$, as a mapping from $\mathfrak{B}_1\times\mathfrak{B}_2$ to $\mathfrak{B}_3$ also as a mapping from $\mathfrak{C}_1\times\mathfrak{C}_2$ to $\mathfrak{C}_3$, is Fr\'{e}chet differentiable at $(x_0,y_0)$,
\item[(iv)]the partial Fr\'{e}chet derivative $D_x\mcl{F}(x_0,y_0)$, as a  mapping from  $\mathfrak{B}_1\times\mathfrak{B}_2$ to $\mathfrak{B}_3$ also as a mapping from $\mathfrak{C}_1\times\mathfrak{C}_2$ to $\mathfrak{C}_3$,   is invertible.
\end{itemize}
Then there is a small neighborhood $\mcl{U}_2(y_0)$ of $y_0$ in $\mathfrak{C}_2$ so that, for any given $y\in \mcl{U}_2(y_0)$, there exists $x_*=x_*(y)$ satisfying $\mcl{F}(x^*(y),y)={\bf{0}}$.
\end{proposition}
In \cite{Sz}, the author used the Brouwer fixed point theorem to prove \emph{the right inverse function theorem} for finite dimensional Banach spaces without assuming continuous differentiability of a mapping. As an analogy, we apply the Schauder fixed point theorem to prove Proposition \ref{lemma5-1}. The detailed proof is given in Appendix \ref{appendix A}.

To apply Proposition \ref{lemma5-1} to $\mcl{P}$, we need to verify that $\mcl{P}$ is continuous, Fr\'{e}chet differentiable at $(Id,\vphi_0^-,p_0^-,v_c)$ and a partial Fr\'{e}chet derivative of $\mcl{P}$ at the point is invertible.
\subsection{Fr\'{e}chet differentiability of $\mcl{P}$ at $\zeta_0=(Id,\vphi_0^-,p_0^-,v_c)$}
To prove Fr\'{e}chet differentiability of $\mcl{P}$ at $\zeta_0$, it suffices to consider how to compute the partial Fr\'{e}chet derivative of $\mcl{P}$ with respect to $v_{ex}$ at $\zeta_0$ because the other partial derivatives of $\mcl{P}$ can be obtained similarly.

For $\mcl{P}$ defined in \eqref{def-p}, let us define $\mcl{R}$ and $\mcl{Q}$ by
\begin{equation}
\label{5-2}
\begin{split}
&\mcl{R}:(\Psi,\vphi_-,p_-,v_{ex})\mapsto (B_0-\frac 12|D\Psi^{-1}\grad\vphi|^2)^{\frac{\gam}{\gam-1}}|_{\Gam_{ex}}\\
&\mcl{Q}:(\Psi,\vphi_-,p_-,v_{ex})\mapsto \frac{\mcl{P}(\Psi,\vphi_-,p_-,v_{ex})}{\mcl{R}(\Psi,\vphi_-,p_-,v_{ex})}.
\end{split}
\end{equation}
If $\mcl{R}$ and $\mcl{Q}$ are Fr\'{e}chet differentiable at $\zeta_0$ then the Fr\'{e}chet derivative $D\mcl{P}$ of $\mcl{P}$ at $\zeta_0$ is given by
\begin{equation}
\label{5-3}
D\mcl{P}_{\zeta_0}=\mcl{R}(\zeta_0)D\mcl{Q}_{\zeta_0}+\mcl{Q}(\zeta_0)D\mcl{R}_{\zeta_0}.
\end{equation}

To simplify notations, we denote as $D_v\mcl{P},D_v\mcl{Q}$ and $D_v\mcl{R}$ respectively for the partial Fr\'{e}chet derivatives of $\mcl{P},\mcl{Q}$ and $\mcl{R}$ at $\zeta_0$ with respect to $v_{ex}$ hereafter.

\begin{lemma}
\label{frecheD}
Fix any $\alp\in(0,1)$, and set $\mathfrak{A}:= C^{1,\alp}_{(-\alp,\der\Lambda)}(\Lambda)$.  As mappings from $\mcl{B}_{\sigma}^{(1)}(\zeta_0)$ to $\mcl{B}_{C\sigma}^{(2)}(p_c)$,
$\mcl{Q},\mcl{R}$ and $\mcl{P}$ are Fr\'{e}chet differentiable at $\zeta_0$. In particular, the partial Fr\'{e}chet derivatives of $\mcl{Q},\mcl{R}$ and $\mcl{P}$ at $\zeta_0$ with respect to $v_{ex}$ are given by
\begin{align}
\label{new1}
&D_v\mcl{R}: w\in \mathfrak{A}\mapsto a_1w\in\mathfrak{A}\\
\label{5-15}
&D_v\mcl{Q} w\in \mathfrak{A}\mapsto  a_2\psi_0(r_s,\cdot)\in \mathfrak{A}\\
&\label{5-16}
D_v\mcl{P}w\in \mathfrak{A}\mapsto  \mcl{Q}(\zeta_0)a_1 w+\mcl{R}(\zeta_0)a_2\psi_0(r_s,\cdot)\in \mathfrak{A}
\end{align}
for any $w\in \mathfrak{A}$ with
\begin{equation}\label{5-A2}
\begin{split}
&a_1\!=-\frac{2\gam(B_0-\frac 12|\der_r\vphi_0^+|^2)^{\frac{\gam}{\gam-1}}\der_r\vphi_0^+}{(\gam+1)r_1^{n-1}(K_0-|\der_r\vphi_0^+|^2)}|_{r=r_1}\\
&a_2\!=\frac{1}{(B_0\!-\!\frac 12|\nabla\vphi_0^+(r_s)|^2)^{\frac{\gam}{\gam-1}}}\Bigl(\frac{\frac{d}{dr}(p_{s,0}\!-\!p_0^+)}{\der_r(\vphi_0^-\!-\!\vphi_0^+)}\!+\!\frac{\gam p_0^+\mu_0\der_r\vphi_0^+}{(\gam\!-\!1)(B_0\!-\!\frac 12|\der_r\vphi_0^+|^2)}|_{r=r_s}\Bigr)
\end{split}
\end{equation}
with $\mu_0$ defined in Lemma \ref{lemma-1}, and
\begin{equation}
\label{ps0}
p_{s,0}=\rho_0^-|\nabla\vphi_0^-\cdot\hat r|^2+p_0^--\rho_0^-K_0
\end{equation}
 where $K_0$ is defined in \eqref{k0}. Here, $\psi_0$ is a unique solution to an elliptic boundary problem specified in the proof.

\begin{proof} In this proof, the constants $\eps_0$ and $C$ appeared in various estimates depend on the data in sense of Remark \ref{remark-cons}.

Since \eqref{5-16} follows from \eqref{new1} and \eqref{5-15}, it suffices to prove \eqref{new1} and \eqref{5-15}.

If a bounded linear mapping $L:\mathfrak{A}\to \mathfrak{A}$ satisfies
\begin{equation}
\label{fre}
\|\mcl{Q}(Id,\vphi_0^-,p_0^-,v_c+\eps w)-\mcl{Q}(Id,\vphi_0^-,p_0^-,v_c)-\eps Lw\|_{1,\alp,\Lambda}^{(-\alp,\der\Lambda)}\le o(\eps)
\end{equation}
for any $w\in \mathfrak{A}$ with $\|w\|_{1,\alp,\Lambda}^{(-\alp,\der\Lambda)}=1$ and if $o(\eps)$ is independent of $w$, then we have $L=D_v\mcl{Q}$. To find such $L$, we fix $w\in \mathfrak{A}$ with $\|w\|_{1,\alp,\Lambda}^{(-\alp,\der\Lambda)}=1$, and compute the G\^{a}teaux derivative $
\frac{d}{d\eps}|_{\eps=0}\mcl{Q}(Id,\vphi_0^-,p_0^-,v_c+\eps w).
$

Fix a sufficiently small constant $\eps_0$, and, for $\eps\in[-\eps_0,\eps_0]\setminus\{0\}$, let $\vphi_{\eps}$ be a solution indicated in Proposition \ref{theorem3-1} with $\Psi=Id, (\vphi_-,p_-)=(\vphi_0^-,p_0^-)$ and $v_{ex}=v_c+\eps w$ on $\Gam_{ex}$ in \eqref{3-0-7}.

Let us denote the subsonic domain of $\vphi_{\eps}$ as $\mcl{N}^+_{\eps}$, the transonic shock as $S_{\eps}=\{r=f_{\eps}(\vartheta)\}$, $\Phi_{f_{\eps}}$ as $\Phi_{\eps}$, and $G_{\vphi_{\eps}}$ as $G_{\eps}$ where $\Phi_{f_{\eps}}$ and $G_{\vphi_{\eps}}$ are defined in \eqref{ver3-1} and \eqref{3-2-20}. For $W^*_{\eps}$, defined as in \eqref{c}, let $E_{\eps}\in C^{1,\alp}_{(-\alp,\Gam_w)}(\mcl{N}^+_{r_s})$ be a solution to
\begin{align}
\label{5-4}
&W^*_{\eps}\cdot DE_{\eps}=0\;\text{in}\;\;\mcl{N}^+_{r_s}\\
\label{5-5}
&E_{\eps}(r_s,\vartheta)=\frac{p_{s,\eps}}{(B_0-\frac 12|\nabla\vphi_{\eps}|^2)^{\frac{\gam}{\gam-1}}}(f_{\eps}(\vartheta),\vartheta)
\end{align}
with
$p_{s,\eps}=\rho_0^-(\grad\vphi_0^-\cdot\nu_{s,\eps})^2+p_0^--\rho_0^-K_{s,\eps}$ where $D=(\der_r,\der_{\vartheta_1},\cdots,\der_{\vartheta_{n-1}})$ and $\nu_{s,\eps}$ is the unit normal vector field on $S_{\eps}$ toward $\mcl{N}_{\eps}^+$, and $K_{s,\eps}$ is defined by (\ref{3-0-21}) with $({\rho}_-,{\vphi}_-,{p}_-)=(\rho_0^-,\vphi_0^-,p_0^-)$ and $\Psi=Id$. By Proposition \ref{proposition4-0-1}, the solution $E_{\eps}$ to (\ref{5-4}) with (\ref{5-5}) is unique, and \eqref{5-2} implies
$
\mcl{Q}(Id,\vphi_0^-,p_0^-,v_c+\eps w)=E_{\eps}|_{\Gam_{ex}}.
$
Similarly, $\mcl{Q}(\zeta_0)$ is the exit value of the solution $E_0^+$ to
\begin{align*}
%\label{5-6}
&W^*_0\cdot DE_0^+=0\;\;\text{in}\;\;\mcl{N}^+_{r_s},
\quad
E_0^+(r_s,\vartheta)=\frac{p_{s,0}(r_s)}{(B_0-\frac 12|\der_r\vphi_0^+(r_s)|^2)^{\frac{\gam}{\gam-1}}}
\end{align*}
where $W^*_0$ is defined by \eqref{c} with $(\Psi,\vphi,{\rho}_-,{\vphi}_-,{p}_-)=(Id,\vphi_0^+,\rho_0^-,\vphi_0^-,p_0^-)$, and $p_{s,0}$ is defined in \eqref{ps0}. Since $\vphi_0^+$ is a radial function, we have $W_0^*=(1,0,\cdots,0)=\hat r$, so $E_0^+$ is a constant in $\mcl{N}^+_{r_s}$, and this implies
\begin{equation}
\label{d}
\frac{d}{d\eps}|_{\eps=0}\mcl{Q}(Id,\vphi_0^-,p_0^-,v_c+\eps w)\!=\frac{d}{d\eps}|_{\eps=0}E_{\eps}(r_1,\cdot)\!=\frac{d}{d\eps}|_{\eps=0}(E_{\eps}\!-E_0^+)(r_1,\cdot).
\end{equation}

Note that, by Proposition \ref{theorem3-1} and Lemma \ref{lemma-reg}, if $\eps_0>0$ is sufficiently small, then whenever $0<|\eps|<\eps_0$, there exists a unique $\psi_{\eps}\in C^{2,\alp}(\mcl{N}^+_{\eps})$ to satisfy
$$
\vphi_{\eps}=\vphi_0^++\eps\psi_{\eps}\;\text{in}\;\;\mcl{N}^+_{\eps},
$$
and moreover, there is a constant $C>0$  so that there holds
\begin{equation}
\label{d10-1}
\begin{split}
&\|\psi_{\eps}\|_{2,\alp,\mcl{N}^+_{\eps}}^{(-1-\alp,\ol{S}_{\eps}\cup\ol{\Gam}_{ex})}\le C,\\
&\|W^*_{\eps}-W^*_0\|_{1,\alp,\mcl{N}^+_{r_s}}^{(-1-\alp,\Gam_w)}+\|W^*_{\eps}-W^*_0\|_{1,\alp,\mcl{N}^+_{r_s}}^{(-1-\alp,\ol{S}_0\cup\ol{\Gam}_{ex})}\le C|\eps|.
\end{split}
\end{equation}

By \eqref{d}, we formally differentiate $W_{\eps}\cdot D(E_{\eps}-E_0^+)=0$ with respect to $\eps$, and apply \eqref{4-0-11},\eqref{4-0-7} and \eqref{d10-1} to obtain
\begin{equation}
\label{limit-eqn}
W_0\cdot\frac{d}{d\eps}|_{\eps=0}(E_{\eps}-E_0^+)=0\quad\text{in}\;\;\mcl{N}^+_{r_s}.
\end{equation}
Then \eqref{limit-eqn} implies
\begin{equation}
\label{limit-eqn2}
\begin{split}
&\phantom{=}\frac{d}{d\eps}|_{\eps=0}\mcl{Q}(Id,\vphi_0^-,p_0^-,v_c+\eps w)=\frac{d}{d\eps}|_{\eps=0}(E_{\eps}-E_0^+)(r_s)\\
&=\frac{d}{d\eps}|_{\eps=0}\Bigl(\frac{p_{s,\eps}}{(B_0-\frac 12|D\vphi_{\eps}|^2)^{\frac{\gam}{\gam-1}}}-\frac{p_0^+}{(B_0-\frac 12|D\vphi_0^+|^2)^{\frac{\gam}{\gam-1}}}\Bigr)|_{(f_{\eps}(\theta),\theta)}
\end{split}
\end{equation}
where $p_{s,\eps}$ is given in \eqref{5-5}.

Set $p_s(\eps,f_{\eps}(\theta),\theta):=p_{s,\eps}(f_{\eps}(\theta),\theta)$. By \eqref{d10-1}, \eqref{2-p-bc} and \eqref{normal}, $p_s(\eps,\cdot,\cdot)$ is partially differentiable with respect to $\eps$ at $\eps=0$. Then this provides
\begin{equation}
\label{derivative}
\frac{d}{d\eps}|_{\eps=0}(p_{s,\eps}\!-\!p_0^+)(f_{\eps}(\theta),\theta)=\!\frac{\der}{\der \eps}|_{\eps=0}p_s(\eps,r_s,\theta)\!+\!(p_{s,0}-\!p_0^+)'(r_s)\frac{d}{d\eps}|_{\eps=0}f_{\eps}(\theta).
\end{equation}

First of all, An explicit calculation using \eqref{normal} and \eqref{d10-2} implies
\begin{equation}
\label{bb-1}
\frac{\der}{\der\eps}|_{\eps=0}p_s(\eps,r_s,\theta)=0.
\end{equation}

Set $\til{\psi}_{\eps}:=\psi_{\eps}\circ\Phi_{\eps}^{-1}$ then \eqref{ver3-1} implies $\psi_{\eps}(f_{\eps}(\vartheta),\vartheta)=\til{\psi_{\eps}}(r_s,\vartheta)$. We note that $\eps\til{\psi}_{\eps}$ solves \eqref{d4-2'} with $\phi=\vphi_0^++\eps\psi_{\eps}(=\vphi_{\eps})$ and $\til g_3=\eps w$. Since all the coefficients, including $\wtil{F}, \til g_1$ and $\til g_2$ in \eqref{d4-2'} smoothly depend on $\phi$ and $D\phi$, $\til{\psi}_{\eps}$ converges to $\psi_0$ in the norm of $C^{2,\alp}_{(-1-\alp,\Gam_w)}(\mcl{N}^+_{r_s})$ where $\psi_0$ is a unique solution to
\begin{align}
\label{5-13a}
&\der_{j}(a_{jk}(x,0)\der_k\psi_0)=0\;\;\text{in}\;\;\mcl{N}^+_{r_s},\quad\quad\;\;
D\psi_0\cdot \hat r-\mu_0\psi_0=0\;\;\text{on}\;\;S_{0}\\
\label{5-13d}
&(a_{jk}(x,0)\der_k\psi_0)\cdot\nu_w=0\;\;\text{on}\;\;\Gam_{w,r_s},\;\;
(a_{jk}(x,0)\der_k\psi_0)\cdot\hat{r}=w\;\;\text{on}\;\;\Gam_{ex},
\end{align}
where $a_{jk}, \mu_0$ are defined in  \eqref{3-1-2} and Lemma \ref{lemma-1}, Also, by a method similar to the proof of Lemma \ref{lemma-frechet}, we have
\begin{equation}
\label{d10-2}
\|\psi_{\eps}\circ\Phi_{\eps}^{-1}-\psi_0\|_{2,\alp,\mcl{N}^+_{r_s}}^{(-1-\alp,\Gam_w)}\le C|\eps|.
\end{equation}
By \eqref{3-0-4}, for any $\eps\in[-\eps_0,\eps_0]\setminus\{0\}$, we have $(\vphi_0^--\vphi_0^+)(f_{\eps}(\vartheta))=\eps\psi_{\eps}(f_{\eps}(\vartheta),\vartheta)$. Combining this with \eqref{derivative}, \eqref{d10-2} and \eqref{bb-1}, we obtain
\begin{equation}
\label{derivative3}
\frac{d}{d\eps}|_{\eps=0}(p_{s,\eps}-p_0^+)(f_{\eps}(\theta),\theta)=\frac{(p_{s,0}-p_0^+)(r_s)}{\der_r(\vphi_0^--\vphi_0^+)(r_s)}\psi_0(r_s,\theta).
\end{equation}
Similarly, by using \eqref{5-13a}, \eqref{d10-2}, we also obtain

\begin{equation}
\label{derivative2}
\begin{split}
&\phantom{=}\frac{d}{d\eps}|_{\eps=0}\Bigl(\frac{1}{(B_0-\frac 12|\nabla\vphi_{\eps}|^2)^{\frac{\gam}{\gam-1}}}-\frac{1}{(B_0-\frac 12|\der_r\vphi_0^+|^2)^{\frac{\gam}{\gam-1}}}\Bigr)(f_{\eps}(\vartheta),\vartheta)\\
&=\frac{\gam\der_r\vphi_0^+(r_s)}
{(\gam-1)(B_0-\frac 12(\der_r\vphi_0^+(r_s))^2)^{\frac{2\gam-1}{\gam-1}}}\der_r\psi_0(r_s,\vartheta)\\
&=\frac{\gam\der_r\vphi_0^+(r_s)\mu_0}{(\gam-1)(B_0-\frac 12(\der_r\vphi_0^+(r_s))^2)^{\frac{2\gam-1}{\gam-1}}}\psi_0(r_s,\vartheta)=:a_3\psi_0(r_s,\vartheta).
\end{split}
\end{equation}
Set $L w:=a_2\psi_0(r_s,\cdot)$ for $a_2$ defined in \eqref{5-A2}.

Since the constant $C$ in \eqref{4-0-11},\eqref{4-0-7}, \eqref{d10-1} and \eqref{d10-2} only depends on the data in sense of Remark \ref{remark-cons} , combining \eqref{derivative3} with \eqref{derivative2}, we conclude that the mapping $L$ satisfies \eqref{fre} with $o(\eps)$ being independent of $w$. Thus we conclude $D_v\mcl{Q}=L$. Since it is easier to verify \eqref{new1} so we omit the proof.
\end{proof}
\end{lemma}

\begin{corollary}
\label{corollary5-2}
$D_v\mcl{Q}:C^{1,\alp}_{(-\alp,\der\Lambda)}(\Lambda)\to C^{1,\alp}_{(-\alp,\der\Lambda)}(\Lambda)$ is a compact mapping.
\begin{proof}
For any $w\in C^{1,\alp}_{(-\alp,\der\Lambda)}(\Lambda)$, let $\psi_0$ be a unique solution to (\ref{5-13a}), (\ref{5-13d}). Then there are constants $\wtil{C}, C>0$ depending only on the data to satisfy
\begin{equation*}
%\label{5-17}
\|D_v\mcl{Q}w\|_{2,\alp,\Lambda}^{(-1-\alp,\der\Lambda)}\le \wtil{C}\|\psi_0\|_{2,\alp,\mcl{N}^+_{r_s}}^{(-1-\alp,\Gam_w)}\le C\|w\|_{1,\alp,\Lambda}^{(-\alp,\der\Lambda)}.
\end{equation*}
\end{proof}
\end{corollary}

\subsection{Local invertibility of $\mcl{P}$ near $\zeta_0=(Id,\vphi_0^-,p_0^-,v_c)$}\label{subsec-localinv}
For a local invertibility of $\mcl{P}$, by Lemma \ref{lemma5-1}, it remains to show that $D_v\mcl{P}$ is invertible and $\mcl{P}$ is continuous in a Banach space.

Note that for any $\alp\in(0,1)$,  $\mcl{B}^{(1)}_{\sigma}(\zeta_0)$ and $\mcl{B}_{C\sigma}^{(2)}(p_c)$ are compact and convex subsets of $\mathfrak{B}_{(1)}:=C^{2,{\alp}/{2}}(\mcl{N},\R^n)\times C^{3,{\alp}/{2}}(\mcl{N}^-_{r_s+\delta})\times C^{2,{\alp}/{2}}(\mcl{N}^-_{r_s+\delta})\times C^{1,{\alp}/{2}}_{-{\alp}/{2},\der\Lambda}(\Lambda)$ and $\mathfrak{B}_{(2)}:=C^{1,{\alp}/{2}}_{(-{\alp}/{2},\der\Lambda)}(\Lambda)$ respectively. We define, for $z=(z_1,z_2,z_3,z_4)$,
\begin{equation}
\label{def-12}
\begin{split}
&\|z\|_{(1)}\!\!:
=\!\|z_1\|_{2,\alp/2,\mcl{N}}\!\!+\!\!\|z_2\|_{3,\alp/2,\mcl{N}^-_{r_s+\delta}}
\!\!\!\!\!\!+\!\|z_3\|_{2,\alp/2,\mcl{N}^-_{r_s+\delta}}\!\!\!+\!\|z_4\|_{1,\alp/2,\Lambda}^{(-\alp/2,\der\Lambda)}\\
&\|z_4\|_{(2)}:=\|z_4\|_{1,\alp/2,\Lambda}^{(-\alp/2,\der\Lambda)},
\end{split}
\end{equation}
and regard $\mcl{P}$ as a mapping from a subset of $\mathfrak{B}_{(1)}$ to a subset of $\mathfrak{B}_{(2)}$.
\begin{lemma}
\label{lemma-cont}
For any $\alp\in(0,1)$, $\mcl{P}:\mcl{B}_{\sigma}^{(1)}(\zeta_0)\to \mcl{B}_{C\sigma}(p_c)$ is continuous in sense that if $\zeta_j$ converges to $\zeta_{\infty}$ in $\mathfrak{B}_{(1)}$ as $j$ tends to $\infty$, then $\mcl{P}(\zeta_j)$ converges to $\mcl{P}(\zeta_{\infty})$ in $\mathfrak{B}_{(2)}$.
\begin{proof}
Take a sequence $\{\zeta_j=(\Psi_j,\vphi_-^{(j)},p_-^{(j)},v_j)\}$ in $\mcl{B}^{(1)}_{\sigma}(\zeta_0)$ so that $\zeta_j$ converges to $\zeta_{\infty}\in \mcl{B}^{(1)}_{\sigma}(\zeta_0)$ in $\mathfrak{B}_{(1)}$ as $j$ tends to $\infty$. For each $j\in\mathbb{N}$, let $E_j $ and $E_{\infty}$ be the solutions to the transport equation of \eqref{2-p-nn} with $W^*\!=\!W^*_j$ and $W^*\!=\!W_{\infty}^*$ respectively in \eqref{c} where the initial conditions $E_{j}(r_s,\vartheta), E_{\infty}(r_s,\vartheta)$ in \eqref{4-0-5} are determined by the data $\zeta_j$ and $\zeta_{\infty}$ respectively. Then, we have $\mcl{Q}(\zeta_j)=E_j|_{\Gam_{ex}}$ and $\mcl{Q}(\zeta_{\infty})=E_{\infty}|_{\Gam_{ex}}$. By the convergence of $\zeta_j$ in $\mathfrak{B}_{(1)}$, $W^*_j$ converges to $W^*_{\infty}$ in $C^{1,\alp/2}_{(-\alp/2,\ol{S_0}\cup\ol{\Gam_{ex}})}(\mcl{N}^+_{r_s})\cap C^{1,\alp/2}_{(-\alp/2,\Gam_w)}(\mcl{N}^+_{r_s})$ as $j$ tends to $\infty$. Because of $\{\zeta_j\}\subset \mcl{B}^{(1)}_{\sigma}(\zeta_0)$, by Proposition \ref{theorem3-1} and \ref{proposition4-0-1}, $\{E_j\}$ is  bounded in $C^{1,\alp}_{(-\alp,\Gam_w)}(\mcl{N}_{r_s}^+)$ so there is a subsequence of $\{E_j\}$, which we simply denote as $\{E_j\}$ itself, that converges to, say $\wtil E_{\infty}$, in $C^{1,\alp/2}_{(-\alp/2,\Gam_w)}(\mcl{N}^+_{r_s})$. Then $\wtil E_{\infty}$ satisfies
\begin{equation*}
W^*_{\infty}\cdot D\wtil E_{\infty}=0\;\;\text{in}\;\;\mcl{N}^+_{r_s},\quad
\wtil E_{\infty}(r_s,\vartheta)=E_{\infty}(r_s,\vartheta)\;\;\text{for all}\;\;\vartheta\in\Lambda,
\end{equation*}
and thus $\wtil E_{\infty}=E_{\infty}$ holds in $\mcl{N}^+_{r_s}$. Moreover, this holds true for any subsequence of $\{E_j\}$ that converges in $C^{1,\alp/2}_{(-\alp/2,\Gam_w)}(\mcl{N}^+_{r_s})$. Therefore, we conclude that $\mcl{Q}(\zeta_j)$ converges to $\mcl{Q}(\zeta_{\infty})$ in $\mathfrak{B}_{(2)}$ i.e., $\mcl{Q}$ is continuous. One can similarly prove that $\mcl{R}$ is continuous. By the definition of $\mcl{P}$ in \eqref{def-p}, we finally conclude that $\mcl{P}$ is continuous.
\end{proof}
\end{lemma}

To prove the invertibility of $D_v\mcl{P}:C^{1,\alp}_{(-\alp,\der\Lambda)}(\Lambda)\to C^{1,\alp}_{(-\alp,\der\Lambda)}(\Lambda)$, let us make the following observation.

For a given $w\in C^{1,\alp}_{(-\alp,\der\Lambda)}(\Lambda)$, let $\psi^{(w)}$ be the solution to (\ref{5-13a}),(\ref{5-13d}). Define a mapping $T$ by
$
T:w\mapsto -\frac{a_2}{a_1}\psi^{(w)}(r_s,\cdot)
$, then $T$ maps $C^{1,\alp}_{(-\alp,\der\Lambda)}(\Lambda)$ into itself, and we can write $D_v\mcl{P}$ in Lemma \ref{frecheD} as
\begin{equation}
\label{kl}
D_v\mcl{P} w=a_1(I-T)w.
\end{equation}
By Corollary \ref{corollary5-2}, $T$ is compact. Then, by the Fredholm alternative theorem, either $D_v\mcl{P}w=0$ has a nontrivial solution $w$ or $D_v\mcl{P}$ is invertible. In the following lemma,we show that $D_v\mcl{P}$ satisfies the latter hence $D_v\mcl{P}$ is invertible.
\begin{lemma}
\label{corollary5-1} For $\alp\in(0,1)$,
$D_v\mcl{P}:C^{1,\alp}_{(-\alp,\der\Lambda)}(\Lambda)\to C^{1,\alp}_{(-\alp,\der\Lambda)}(\Lambda)$ is invertible.
\begin{proof}
The proof is divided into two steps.

\textbf{Step 1.} Let $a_1$ and $a_2$ be defined as in {\rm{(\ref{5-A2})}}. We claim $a_2<0$ and $-\frac{a_2}{a_1}<0$.

By Lemma \ref{lemma-1}, we can write $a_2$ as
\begin{equation*}
a_2=\frac{\der_r(p_{s,0}-p_0^+)-\rho_0^-\der_r\vphi_0^-\der_r(\frac{K_0}{\der_r\vphi_0^-}-\der_r\vphi_0^+)}{\der_r(\vphi_0^--\vphi_0^+)(r_s)}
\mid_{r=r_s}.
\end{equation*}
Then, by \eqref{2-3-6}, \eqref{2-3-8} and \eqref{3-1-12}, we obtain
\begin{equation*}
a_2=
\frac{(n-1)\rho_0^-(r_s)(K_0-(\der_r\vphi_0^-(r_s))^2)}{r_s(B_0-\frac 12|\grad\vphi_0^+(r_s)|)^{\frac{\gam}{\gam-1}}\der_r(\vphi_0^--\vphi_0^+)(r_s)}<0.
\end{equation*}
Since $a_1<0$ is obvious from \eqref{5-16}, there holds $-\frac{a_2}{a_1}<0$.

\textbf{Step 2.} \emph{Claim. $D_v\mcl{P} w=0$ if and only if $w=0$.}

If $w=0$ then it is obviously $D_v\mcl{P}w$=0.

According to Remark \ref{laplacian-sphere},  in spherical coordinates, (\ref{5-13a}) and (\ref{5-13d}) are expressed as
\begin{equation}
\label{5-24}
\begin{split}
&\der_r\{k_1(r)\der_r\psi\}+k_2(r)\triangle_{S^{n-1}}\psi=0\;\;\text{in}\;\;\mcl{N}^+_{r_s},\;\;\der_r\psi-\mu_0\psi=0\;\;\text{on}\;\;S_0,\\
&\der_{\nu_w}\psi=0\;\;\text{on}\;\;\Gam_{w,r_s},\;\;\quad\quad\qquad\qquad\qquad\quad\;\;\der_r\psi=\frac{w}{k_1(r_1)}\;\;\text{on}\;\;\Gam_{ex}
\end{split}
\end{equation}
where $k_1$ and $k_2$ are defined as in (\ref{d8-1}).

Suppose $D_v\mcl{P}w=0$ for some $w\in C^{1,\alp}_{(-\alp,\der\Lambda)}(\Lambda)$. Then, by \eqref{kl}, there holds
\begin{equation}
\label{g}
w(\vartheta)=-\frac{a_2}{a_1}\psi^{(w)}(r_s,\vartheta)\;\;\text{for all}\;\;\vartheta\in\Lambda.
\end{equation}
We denote $\psi^{(w)}$ as $\psi$ for the simplicity hereafter. Let $\nu_*$ be a unit normal vector field on $\der\Lambda\subset \mathbb{S}^{n-1}$. Then $\nu_*$ is perpendicular to $\Gam_w$ so $w$ satisfies
\begin{equation*}
\der_{\nu_*} w=-\frac{a_2}{a_1}\der_{\nu_w}\psi(r_s,\cdot)=0\;\;\text{on}\;\der\Lambda.
\end{equation*}
Then $w$ has an eigenfunction expansion
\begin{equation*}
%\label{5-19}
w=\sum_{j=1}^{\infty}c_j\eta_j\;\;\text{in}\;\;\Lambda\;\;\text{with}\;\;c_j=\int_{\Lambda}w\eta_j
\end{equation*}
where every $\eta_j$ is an eigenfunction of the eigenvalue problem
\begin{equation}
\label{5-egn}
-\triangle_{\mathbb{S}^{n-1}} \eta_j=\lambda_j\eta_j\;\;\text{in}\;\;\Lambda,\quad\quad
\der_{\nu_*}\eta_j=0\;\;\text{on}\;\;\der\Lambda
\end{equation}
with $\|\eta_j\|_{L^2(\Lambda)}=1$. We note that all the eigenvalues of {\rm{(\ref{5-egn})}} are non-negative real numbers.

We claim that \eqref{g} implies  $c_j=0$ for all  $j\in\mathbb{N}$. Once this is verified, then we obtain $w=0$.

For each $N\in\mathbb{N}$, let us set $w_N:=\sum_{j=1}^Nc_j\eta_j$ and $\psi_N(r,\vartheta):=\sum_{j=1}^N q_j(r)\eta_j(\vartheta)$ for one variable functions $q_j$ satisfying
\begin{align}
\label{5-20a}
&(k_1q'_j)'-\lambda_jk_2q_j=0\;\;\text{in}\;\;(r_s,r_1),\\
\label{5-20b}
&q_j'(r_s)-\mu_0 q_j(r_s)=0,\quad q_j'(r_1)=\frac{c_j}{k_1(r_1)}.
\end{align}
for $j=1,\cdots,N$. Then $\psi_N(r_s,\cdot)$ converges to $\psi(r_s,\cdot)$ in $L^2(\Lambda)$ so that \eqref{g} implies

$$q_j(r_s)=-\frac{a_1}{a_2}c_j\;\;\text{for all}\;\;j\in\mathbb{N}.$$

Suppose that $c_j\neq 0$ for some $j\in\mathbb{N}$. Without loss of generality, we assume $c_j>0$. By the sign of $\frac{a_1}{a_2}$ considered in Step 1, we have $q_j'(r_s)<0, q_j(r_s)<0$ and $q_j'(r_1)>0$. We will derive a contradiction.
\begin{itemize}
\item[(i)] If $\lambda_j=0$ then there is a constant $m_0$ so that $k_1q_j'=m_0$ for all $r\in[r_s,r_1]$. Since $k_1$ is strictly positive in $[r_s,r_1]$, we have $sgn\; m_0=sgn\; q_k'(r_1)$ and $sgn\; m_0=sgn\; q_j'(r_s)$. Then $sgn\;q'_j(r_s)=sgn\;q_j'(r_1)$ must hold. But this contradicts to the assumption of $q_j'(r_s)<0<q_j'(r_1)$.
\item[(ii)] Assume $\lambda_j>0$. By (\ref{5-20a}), we have
\begin{equation*}
%\label{5-21}
\int_{r_s}^{r_1}\lambda_jk_2q_j dr=\int_{r_s}^{r_1}(k_1q_j')'dr=k_1(r_1)q_j'(r_1)-k_1(r_s)q_j'(r_s)>0.
\end{equation*}
Then $q_j(r_s)<0$ implies $\max_{[r_s,r_1]}q_j>0$ so there exists $t_1\in[r_s,r_1]$ satisfying $q_j(t_1)=0$. Also, by the intermediate value theorem, there exists $t_2\in[r_s,r_1]$ satisfying $q_j'(t_2)=0$.

If $t_1=t_2$ then the uniqueness theorem of second-order linear ODEs implies $q_j\equiv 0$ on $[r_s,r_1]$. But this is a contradiction. If $t_1\neq t_2$, assume $t_1<t_2$ without loss of generality. Then, by the maximum principle, $q_j\equiv 0$ on $[t_1,t_2]$ so $q_j\equiv 0$ on $[r_s,r_1]$. This is a contradiction as well.
\end{itemize}
Thus every $c_j$ must be $0$ hence $w=0$. By the Fredholm alternative theorem, we finally conclude that $D_v\mcl{P}$ is invertible.
\end{proof}
\end{lemma}

\begin{proposition}
\label{proposition5-1}
For any $\alp\in(0,1)$, there exists a positive constant $\sigma_1(\le\sigma_3)$ depending on the data in sense of Remark \ref{remark-cons} so that whenever $0<\sigma\le \sigma_1$, for a given $(\Psi,\vphi_-,p_-,p_{ex})$ with
$\varsigma_1,\varsigma_2,\varsigma_3$ in Theorem \ref{main-thm'} satisfying $0\le \varsigma_l\le \sigma$ for $l=1,2,3$, there exists a function $v^*_{ex}=v_{ex}^*(\Psi,\vphi_-,p_-,p_{ex})\in C^{1,\alp}_{(-\alp,\der\Lambda)}(\Lambda)$ so that there holds
\begin{equation}
\label{star}
 \mcl{P}(\Psi,\vphi_-,p_-,v_{ex}^*)=p_{ex}.
\end{equation}
Moreover such a $v_{ex}^*$ satisfies
 \begin{equation}
 \label{5-23}
 \|v_{ex}^*-v_c\|_{1,\alp,\Lambda}^{(-\alp,\der\Lambda)}
 \le C(\varsigma_1+\varsigma_2+\varsigma_3).
 \end{equation}
 \begin{proof}
Following the definition \eqref{sets}, define, for the constant $C$ in \eqref{4-0-11}, $\mcl{P}_*:\mcl{B}^*:=\mcl{B}^{(1)}(\zeta_0)\times\mcl{B}^{(13)}_{C\sigma}(p_c)\to \mcl{B}^{(13)}_{2C\sigma}(0)$ by
 \begin{equation*}
 \mcl{P}_*(\Psi,\vphi_-,p_-,v_{ex},p_{ex}):=\mcl{P}(\Psi,\vphi_-,p_-,v_{ex})-p_{ex}.
 \end{equation*}
 Then, by Lemma \ref{lemma-cont} and the definition of $\mcl{P}_*$ above, for $\{(\zeta-_j,p_j)\}\subset\mcl{B}^* $, if $\lim_{j\to\infty}\|\zeta_j-\zeta_{\infty}\|_{(1)}+\|p_j-p_{\infty}\|_{(2)}=0$ for some $(\zeta_{\infty},p_{infty})\in \mcl{B}^*$ where $\|\cdot\|_{k=1,2}$ are defined in \eqref{def-12}, then $\lim_{j\to\infty}\|\mcl{P}_*(\zeta_j,p_j)-\mcl{P}_*(\zeta_{\infty},p_{\infty})\|=0$ holds, i.e., $\mcl{P}_*$ is continuous in the sense similar to Lemma \ref{lemma-cont}. Besides, by Lemma \ref{frecheD}, $\mcl{P}_*$ is Fr\'{e}chet differentiable at $(\zeta_0,p_c)$. In particular, we have $D_v\mcl{P}_*(\zeta_0,p_c)=D_v\mcl{P}$ so it is invertible. So if $\sigma>0$ is sufficiently small depending on the data and $\alp$ then Lemma \ref{lemma5-1} implies \eqref{star}. (\ref{5-23}) is a direct result from Remark \ref{remark-norm}.
\end{proof}
\end{proposition}

\begin{proof}[\textbf{Proof of Theorem \ref{main-thm'}}]
Fix $\alp\in(0,1)$. By Proposition \ref{proposition5-1}, for a given $\Psi,\vphi_-,$ $p_-,p_{ex}$ with $\varsigma_1,\varsigma_2,\varsigma_3$ in Theorem \ref{main-thm'} satisfying $0\le \varsigma_l\le \sigma_1$ for $l=1,2,3$, there exists a $v^*_{ex}\in C^{1,\alp}_{(-\alp,\der\Lambda)}(\Lambda)$ so that there holds \eqref{star}. For such a $v^*_{ex}$, Proposition \ref{theorem3-1} and \ref{proposition4-0-1} imply that there exists a unique $(\rho^*,\vphi^*,p^*)$ in $\mcl{N}^+$ satisfying \eqref{2-b}, \eqref{3-0-7}, Proposition \ref{theorem3-1}, Lemma \ref{lemma-reg} and Proposition \ref{proposition4-0-1}. So $(\rho,\vphi,p)=(\rho^*,\vphi^*,p^*)\chi_{\mcl{N}^+}+(\rho_-,\vphi_-,p_-)\chi_{\mcl{N}\setminus \mcl{N}^+}$ is a transonic shock solution satisfying Theorem \eqref{main-thm'} (a)-(d). Particularly, Theorem \ref{main-thm'} (c) is a result from \eqref{3-0-19}, \eqref{4-0-11} and \eqref{5-23}. Here, $\chi_{\Om}=1$ in $\Om$, and $0$ in $\Om^C$.
\end{proof}

\section{Uniqueness}\label{sec:nozzle5}
\subsection{Proof of Theorem \ref{theorem-uniq}}\label{subsec-uniq-proof}
Fix $\zeta_*:=(\Psi,\vphi_-,p_-)$ and $p_{ex}$ with satisfying Theorem \ref{main-thm'} (i)-(iii), and suppose that $v_1,v_2\in C^{1,\alp}_{(-\alp,\der\Lambda)}(\Lambda)$ satisfy
\begin{equation}
\label{d11-2}
\mcl{P}(\zeta_*,v_1)=p_{ex}=\mcl{P}(\zeta_*,v_2).
\end{equation}
If $\|v_j-v_c\|_{1,\alp,\Lambda}^{(-\alp,\der\Lambda)}\le \sigma_1$ in Theorem \ref{main-thm'}, then by \eqref{d11-2} and Lemma \ref{corollary5-1}, we have

\begin{equation}
\label{d11-1}
v_1-v_2=-D_v\mcl{P}^{-1}(\mcl{P}(\zeta_*,v_1)-\mcl{P}(\zeta_*,v_2)-D_v\mcl{P}(v_1-v_2)).
\end{equation}
If $\mcl{P}$ were continuously Fr\'{e}chet differentiable, (\ref{d11-1}) would imply
$$\|v_1-v_2\|_{1,\alp,\Lambda}^{(-\alp,\der\Lambda)}\le C\sigma_1\|v_1-v_2\|_{1,\alp,\Lambda}^{(-\alp,\der\Lambda)}$$
so that if $\sigma_1$ is sufficiently small,  then $v_1=v_2$ holds, which provides Theorem \ref{theorem-uniq} by Proposition \ref{theorem3-1} and \ref{proposition4-0-1}. But it is unclear whether $\mcl{P}$ is differentiable at other points near $\zeta_0=(Id,\vphi_0^-,p_0^-,v_c)$ for the following reason.

For $v_j(j=1,2)$ chosen above, let $W^*_j$ be defined as in \eqref{c} associated with $\zeta_*=(\Psi,\vphi_-,p_-)$ and $v_j$, and let $E_j\in C^{1,\alp}_{(-\alp,\Gam_w)}(\mcl{N}^+_{r_s})$ be the unique solution to \eqref{2-p-nn} with \eqref{4-0-5}. Then, $E_1-E_2$ satisfies $W_1^*\cdot D(E_1-E_2)=(W^*_2-W^*_1)\cdot D(E_2-E_0^+)$ in $\mcl{N}^+_{r_s}$ with $D=(\der_r,\der_{\vartheta_1},\cdots, \der_{\vartheta_{n-1}})$ by Lemma \ref{lemma-ent-0}. Then, by the method of characteristics used in section \ref{subsec-char},for any $\vartheta\in \Lambda$, $E_1-E_2$ has an expression of
\begin{equation}
\label{l1}
\begin{split}
&(E_1-E_2)(r_1,\vartheta)\\
&=(E_1-E_2)|_{X(2r_1-r_s;r_1,\vartheta)}-\!\!\underset{=:\mcl{I}_*}{\underbrace{\int_{r_1}^{2r_1-r_s}\!\!\!\!(W_2^*\!-\!W_1^*)\cdot \! D(E_2\!-\!E_0^+)|_{X(t;r_1,\vartheta)}dt}}
\end{split}
\end{equation}
where $X(t;r,\vartheta)$ solves \eqref{4-0-9'} with $W=W_1^*$.

For $\mcl{Q}$ in \eqref{5-2} to be Fr\'{e}chet differentiable in sense of Lemma \ref{frecheD}, we need to obtain a uniform bound of $\frac{\|\mcl{I}_*\|_{1,\alp,\Lambda}^{(-\alp,\der\Lambda)}}{\|v_1-v_2\|_{1,\alp,\Lambda}^{(-\alp,\der\Lambda)}}$ for all $v_1$ sufficiently close to $v_2$ with $v_1\neq v_2$, but we only know, by Lemma \ref{lemma4-0-1}, $D(E-E_0^+)\in C^{\alp}_{(1-\alp,\Gam_w)}(\mcl{N}^+_{r_s})$.

But still, by \eqref{new-14} and \eqref{4-0-9'}, it is likely that we have
\begin{equation*}
\|\mcl{I}_*\|_{\alp,\Lambda}^{(1-\alp,\der\Lambda)}\le C\|v_1-v_2\|_{\alp,\Lambda}^{(1-\alp,\der\Lambda)}
\end{equation*}
for a constant $C$ depending on the data in sense of Remark \ref{remark-cons}. Then \eqref{d11-1} may imply $\|v_1-v_2\|_{\alp,\Lambda}^{(1-\alp,\der\Lambda)}\le C\sigma_1\|v_1-v_2\|_{\alp,\Lambda}^{(1-\alp,\der\Lambda)}$, and so $v_1=v_2$ hold by reducing $\sigma_1>0$ in Theorem \ref{main-thm'} if necessary. For that reason, we will prove the following lemma in this section.

\begin{lemma}
\label{corollary1-lemma6.1}For $\alp\in(\frac 12,1)$,
$D_v\mcl{P}$ defined in Lemma \ref{frecheD} is a bounded linear mapping from $C^{\alp}_{(1-\alp,\der\Lambda)}(\Lambda)$ to itself and that $D_v\mcl{P}$ is invertible in $C^{\alp}_{(1-\alp,\der\Lambda)}(\Lambda)$.
\end{lemma}

\begin{lemma}
\label{lemma7-1}
Fix $\alp\in(\frac 12,1)$. Fix $\zeta_*=(\Psi,\vphi_-,p_-)$ and $v_1,v_2\in C^{1,\alp}_{(-\alp,\der\Lambda)}(\Lambda)$ with $\varsigma_1,\varsigma_2,$ defined in Theorem \ref{main-thm'} and $\varsigma_4^{(j)}:=\|v_j-v_c\|_{1,\alp,\Lambda}^{(-\alp,\der\Lambda)}$. Then there is a constant $C,\sigma_*>0$ depending on the data in sense of Remark \ref{remark-cons} so that whenever $0<\varsigma_1,\varsigma_2,\varsigma_4^{(1)},\varsigma_4^{(2)}\le \sigma_*$ are satisfied there holds
\begin{equation}
\label{new-34}
\begin{split}
&\phantom{\le}\|\mcl{P}(\zeta_*,v_1)-\mcl{P}(\zeta_*,v_2)-D_v\mcl{P}(v_1-v_2)\|_{\alp,\Lambda}^{(1-\alp,\der\Lambda)}\\
&\le C(\varsigma_1+\varsigma_2+\sum_{j=1}^2\|v_j-v_c\|_{1,\alp,\Lambda}^{(-\alp,\der\Lambda)})\|v_1-v_2\|_{\alp,\Lambda}^{(1-\alp,\der\Lambda)}.
\end{split}
\end{equation}
\end{lemma}

Before we prove Lemma \ref{corollary1-lemma6.1} and Lemma \ref{lemma7-1}, we assume that those lemmas hold true, and prove Theorem \ref{theorem-uniq} first.
\begin{proof}
[\textbf{Proof of Theorem \ref{theorem-uniq}}] Fix $\alp\in(\frac 12,1)$ and also fix $\Psi, \vphi_-, p_-, p_{ex}$ with $\varsigma_1, \varsigma_2, \varsigma_3$ in Theorem \ref{main-thm'} satisfying $0<\varsigma_l\le \sigma$ for a constant $\sigma\in(0,\sigma_1]$ where $\sigma_1$ is as in Theorem \ref{main-thm'}.

Suppose that there are two functions $v_1,v_2\in C^{1,\alp}_{(-\alp,\der\Lambda)}(\Lambda)$ satisfying \eqref{star}, \eqref{5-23} with $v^*_{ex}=v_j$ for $j=1,2$. Then, by \eqref{5-23}, \eqref{d11-1}, Lemma \ref{corollary1-lemma6.1}, Lemma \ref{lemma7-1}, there holds
\begin{equation}
\label{star2}
\|v_1-v_2\|_{\alp,\Lambda}^{(1-\alp,\der\Lambda)}\le C(\varsigma_1+\varsigma_2+\varsigma_3)\|v_1-v_2\|_{\alp,\Lambda}^{(1-\alp,\der\Lambda)}.
\end{equation}
Since $C$ depends on the data in sense of Remark \ref{remark-cons}, there exists a constant $\sigma_2\in (0,\sigma_1]$ with the same depndence as $C$ so that we have
\begin{equation}
\label{star3}
C(\varsigma_1+\varsigma_2+\varsigma_3)\le 3C\sigma_2<1.
\end{equation}
Then, whenever $0<\sigma\le \sigma_2$, from \eqref{star2} and \eqref{star3}, we conclude $v_1=v_2$.
\end{proof}

The rest of section \ref{sec:nozzle5} is devoted to the proof of Lemma \ref{corollary1-lemma6.1} and Lemma \ref{lemma7-1}.
\subsection{Proof of Lemma \ref{corollary1-lemma6.1}}
Lemma \ref{corollary1-lemma6.1} follows easily from:

\begin{lemma}
\label{prop7-1}
Fix $\alp\in (\frac 12,1)$.
For $F=(F_j)_{j=1}^n\in C^{\alp}_{(1-\alp,\der S_0\cup \der\Gam_{ex})}(\mcl{N}^+_{r_s},\R^n)$, and $g_1, g_3\in C^{\alp}_{(1-\alp,\der\Lambda)}(\Lambda)$, $g_2\in C^{\alp}_{(1-\alp,\der S_0\cup\der\Gam_{ex})}(\Gam_{w,r_s})$, the linear boundary problem
\begin{align}
\label{bd1}
&\der_j(a_{jk}(x,0)\der_k\psi)=\der_j F_j\;\;\text{in}\;\;\mcl{N}^+_{r_s}\\
\label{bd2}
&D\psi\cdot\hat r-\mu_0 \psi=g_1\;\;\text{on}\;\;S_0\\
\label{bd3}
&(a_{jk}(x,0)\der_k\psi)\cdot\nu_w=g_2\;\;\text{on}\;\;\Gam_{w,r_s}\\
\label{bd4}
&(a_{jk}(x,0)\der_k\psi)\cdot\hat r=g_3\;\;\text{on}\;\;\Gam_{ex}
\end{align}
has a unique weak solution $\psi\in C^{1,\alp}(\ol{\mcl{N}^+_{r_s}}\setminus(\der S_0\cup \der \Gam_{ex}))\cap C^{\alp}(\ol{\mcl{N}^+_{r_s}})$. Moreover, there is a constant $C$ depending only the data in sense of Remark \ref{remark-cons} so that there holds
\begin{equation}
\label{weak-est-gen}
\begin{split}
&\|\psi\|_{1,\alp,\mcl{N}^+_{r_s}}^{(-\alp,\der S_0\cup\der\Gam_{ex})}\le C(\|F\|_{\alp,\mcl{N}^+_{r_s}}^{(1-\alp,\der S_0\cup\der\Gam_{ex})}\\
&\phantom{aaaaaaaaaaaaaaaaaaaaaaaaa}
+\sum_{j=1,3}\|g_j\|_{\alp,\Lambda}^{(1-\alp,\der\Lambda)}+\|g_2\|_{\alp,\Gam_{w,r_s}}^{(1-\alp,\der S_0\cup\der \Gam_{ex})}).
\end{split}
\end{equation}
\end{lemma}
 The proof of Lemma \ref{prop7-1} is given in Appendix \ref{appendix B}.
\begin{proof}[\textbf{Proof of Lemma \ref{corollary1-lemma6.1}}]
By \eqref{5-16}, for $D_v\mcl{P}$ to be well-defined in $C^{\alp}_{(1-\alp,\der\Lambda)}(\Lambda)$, the linear boundary problem \eqref{5-13a}, \eqref{5-13d} must have a unique solution $\psi^{(w)}\in C^2(\mcl{N}^+_{r_s})\cap C^0(\ol{\mcl{N}^+_{r_s}})$ for $w\in C^{\alp}_{(1-\alp,\der\Lambda)}(\Lambda)$. This condition is satisfies by Lemma \ref{prop7-1}. Furthermore, \eqref{weak-est-gen} implies that $D_v\mcl{P}$ maps $C^{\alp}_{(1-\alp,\der\Lambda)}(\Lambda)$ into itself, and it is bounded in  $C^{\alp}_{(1-\alp,\der\Lambda)}(\Lambda)$ as well. Also, the mapping $w\mapsto a_2\mcl{R}(\zeta_0)\psi^{(w)}(r_s,)$ is compact where $\mcl{R}, a_2$ are defined in \eqref{5-2} and \eqref{5-A2}.

Suppose $D_v\mcl{P}w_*=0$ for some $w_*\in C^{\alp}_{(1-\alp,\der\Lambda)}(\Lambda)$, then by \eqref{5-16}, we have
\begin{equation}
\label{new-35}
w=-\frac{\mcl{R}(\zeta_0)a_2}{\mcl{Q}(\zeta_0)a_1}\psi^{(w)}(r_s,\cdot).
\end{equation}
By \eqref{weak-est-gen}, $\psi^{(w)}$ is in $C^{\alp}(\ol{\mcl{N}^+_{r_s}})$, then \eqref{new-35} implies $w_*\in C^{\alp}(\ol{\Lambda})$, and again this provides $w_*=-\frac{\mcl{R}(\zeta_0)a_2}{\mcl{Q}(\zeta_0)a_1}\psi^{(w_*)}(r_s,\cdot)\in C^{1,\alp}(\ol{\Lambda})$. This can be checked by arguments similar to the proof of Lemma \ref{lemma-linear-existence}. So one can follow the proof of Lemma \ref{corollary5-1} to show $w_*=0$. By the Fredholm alternative theorem, we conclude that $D_v\mcl{P}$ is invertible in $C^{\alp}_{(1-\alp,\der\Lambda)}(\Lambda)$.
\end{proof}
\subsection{Proof of Lemma \ref{lemma7-1}}
\label{subsec:proof-theorem-uniq}
To prove Lemma \ref{lemma7-1}, we need the following corollaries of Lemma \ref{prop7-1}. In this section, we fix $\Psi, \vphi_-, p_-$ with $\varsigma_1, \varsigma_2$ defined in Theorem \ref{main-thm'} satisfying $0<\varsigma_l\le \sigma (l=1,2)$ for a constant $\sigma>0$.

\begin{corollary}
\label{corollary2-lemma6.1} Fix $\alp\in(\frac 12,1)$. For $j=1,2$, let $\vphi_j$ be the solution stated in Proposition \ref{theorem3-1} with $v_{ex}=v_j$ satisfying $\varsigma_4^{(j)}=\|v_j-v_c\|_{1,\alp,\Lambda}^{(-\alp,\der\Lambda)}\le \sigma$ for $\sigma\in(0,\sigma_3]$, and with the transonic shock $r=f_j(x')$. Set $\phi_j:=(\vphi_j-\vphi_0^+)\circ G^{-1}_{\vphi_j}$ for $G_{\vphi_j}$ defined in \eqref{3-2-20}. Then, there exist constant $C,\sigma^{\flat}>0$ depending on the data in sense of Remark \ref{remark-cons} so that whenever $0<\sigma\le \sigma^{\flat}$ if $0<\varsigma_1, \varsigma_2, \varsigma_4^{(j)}(j=1,2)\le \sigma$ are satisfied, then there hold
\begin{align}
\label{new-at-1}
&\|\phi_1-\phi_2\|_{1,\alp,\mcl{N}^+_{r_s}}^{(-\alp,\corners)}\le C\|v_1-v_2\|_{\alp,\Lambda}^{(1-\alp,\der\Lambda)},\\
\label{new-at-2}
&\|f_1-f_2\|_{1,\alp,\Lambda}^{(-\alp,\der\Lambda)}\le C\|v_1-v_2\|_{\alp,\Lambda}^{(1-\alp,\der\Lambda)}.
\end{align}
\begin{proof}
\textbf{Step 1.} For $G_{\vphi}$ defined in \eqref{3-2-20}, let us write $G^{-1}_{\vphi}(t,x')=(\omega^{(\vphi)}(t,x'),x')$ then we have
\begin{equation*}
\upsilon^{(\vphi)}(\omega^{(\vphi)}(t,x'),x')=t\;\;\text{for all}\;\;(t,x')\in[r_s,r_1]\times \ol{\Lambda}.
\end{equation*}
For convenience, we denote as $\upsilon_j, \omega_j$ for $\upsilon^{(\vphi_j)}, \omega^{(\vphi_j)}(j=1, 2)$, then there holds
\begin{equation}
\label{new-at-6}
\upsilon_1(\omega_1(t,x'),x')-\upsilon_2(\omega_2(t,x'),x')=0\;\;\text{for all}\;\;(t,x')\in[r_s,r_1]\times\ol{\Lambda}
\end{equation}
By the definition of $\upsilon_j$ in \eqref{3-2-20}, a direct computation shows that \eqref{new-at-6} is equivalent to
\begin{equation*}
[q(r,x')]_{r=\omega_2(t,x')}^{\omega_1(t,x')}=-k(\phi_1-\phi_2)(t,x')+k[\phi_j(t,x')\chi(\omega_j(t,x'))]_{j=2}^1
\end{equation*}
where we set
\begin{align*}
&q(r,x'):=k(\vphi_0^-(r)-\vphi_0^+(r)+\psi_-(r,x'))(1-\chi(r))+r\chi(r)\\
&\psi_-(r,x'):=\vphi_-(r,x')-\vphi_0^-(r).
\end{align*}
By the choice of $k$ in \eqref{3-2-10}, if $\varsigma_2$ in Theorem \ref{proposition2-2-1} is sufficiently small,  then there is a positive constant $q_0$ satisfying $\der_rq\ge q_0>0$ for all $r\in[r_s,r_1]$ so that we can express $\omega_1-\omega_2$ as
\begin{equation}
\label{new-at-7}
(\omega_1-\omega_2)(t,x')=
\frac{-k(\phi_1-\phi_2)(t,x')+k[\phi_j(t,x')\chi(\omega_j(t,x'))]_{j=2}^1}{\int_0^1\der_rq(a\omega_1(t,x')+(1-a)\omega_2(t,x'))da}
\end{equation}
for $(t,x')\in[r_s,r_1]\times \ol{\Lambda}$. $q_0$ depends only on the data in sense of Remark \ref{remark-cons}. Furthermore, denoting as $G_j, G^{-1}_j$ for $G_{\vphi_j}, G^{-1}_{\vphi_j}(j=1, 2)$, since we have
$
D(G_j\circ G_{j}^{-1})=DG^{-1}_j\circ DG_j(G_j^{-1})=I_n$ for $j=1, 2,$
we obtain, in $\mcl{N}^+_{r_s}$,
\begin{equation}
\label{new-at-8}
\begin{split}
&[DG_j(G_j^{-1})]=-DG_1(G_1^{-1})\!\circ\!(DG_1^{-1}\!-\!DG_2^{-1})\!\circ\! DG_2(G_2)^{-1}.
\end{split}
\end{equation}
with setting $[DG_j(G_j^{-1})]:=DG_1{(G_1^{-1})}\!-\!DG_2(G_2^{-1})$.

Let us set $\Gam_{cns}:=\der S_0\cup \der\Gam_{ex}$. Combining \eqref{new-at-8} with \eqref{new-at-7}, we obtain
\begin{align}
\label{new-at-9}
&\|\omega_1\!-\!\omega_2\|_{1,\alp,\mcl{N}^+_{r_s}}^{(-\alp,\Gam_{cns})}\le C(\|\phi_1\!-\!\phi_2\|_{1,\alp,\mcl{N}^+_{r_s}}^{(-\alp,\Gam_{cns})}\!+\!
\|(\omega_1\!-\!\omega_2)\phi_2\|_{1,\alp,\mcl{N}^+_{r_s}}^{(-\alp,\Gam_{cns})}),\\
\label{new-at-10}
&\|[DG_j(G_j^{-1})]\|_{\alp,\mcl{N}^+_{r_s}}^{(1-\alp,\Gam_{cns})}
\le C(\|\phi_1\!-\!\phi_2\|_{1,\alp,\mcl{N}^+_{r_s}}^{(-\alp,\Gam_{cns})}\!+\!
\|\phi_2(\omega_1\!-\!\omega_2)\|_{1,\alp,\mcl{N}^+_{r_s}}^{(-\alp,\Gam_{cns})}).
\end{align}
\quad\\
\textbf{Step 2.} For convenience, let us denote as $a_{kl}^{(j)}, F_k^{(j)}, b_1^{(j)}, \mu^{(j)},$ $g_m^{(j)}$ for $a_{kl}, F_k, b_1,$ $\mu_f,$ $g_m$ defined in \eqref{3-1-2}, \eqref{3-1-3}, \eqref{3-1-16} and \eqref{new-11} associated with $\Psi, \vphi_-, p_-, \vphi_j-\vphi_0^+, f_j$. Then, by \eqref{3-1-1}-\eqref{3-1-3}, \eqref{3-1-17}, \eqref{new-11} and \eqref{def-gs}, for each $j=1, 2$, $\phi_j$ satisfies
\begin{equation}
\label{new-at-17}
\begin{split}
&\der_k(\til a_{kl}(y,0)\der_l\phi_j)=\der_k\til F^{(j)}_k\!+\!\der_k([\til a_{kl}(y,0)\!-\!\til a^{(j)}_{kl}]\der_l\phi_j)\!=:\der_k\mathfrak F_k^{(j)}\;\;\text{in}\;\;\mcl{N}^+_{r_s},\\
&\der_r\phi_j-\til\mu_0\phi_j=\til g_1^{(j)}\!+\!(\hat r\!-\!\til b^{(j)}_1)\!\cdot \!D\phi_j\!-\!(\til\mu_0\!-\!\til\mu^{f(j)})\phi_j\!=:\mathfrak g_1^{(j)}\;\;\text{on}\;\;S_0\\
&(\til a_{kl}(y,0)\der_l\phi_j)\cdot\nu_w=\til g_2^{(j)}\!+\!([\til a_{kl}(y,0)\!-\!\til a^{(j)}_{kl}]\der_l\phi_j)\!\cdot\! \nu_w=:\mathfrak g_2^{(j)}\;\;\text{on}\;\;\Gam_{w,r_s}\\
&(\til a_{kl}(y,0)\der_l\phi_j)\cdot\hat r=(v_j-v_c)\!+\!([\til a_{kl}(y,0)\!-\!\til a^{(j)}_{kl}]\der_l\phi_j)\!\cdot\!\hat r=:\mathfrak g_3^{(j)}\;\;\text{on}\;\;\Gam_{ex}
\end{split}
\end{equation}
where $\til a_{kl}(y,0), \til{\mu}_0$ are obtained from \eqref{bd1}-\eqref{bd4} through the change of variables $(r,x')\mapsto G_{\vphi_0^+}(r,x')$, and where $\til a^{(j)}_{kl}, \til b^{(j)}_1, \til{\mu}^{(j)}, \til F^{(j)}_k, \til g_m^{(j)}$ for $k,l\in\{1,\cdots,n\},$ $m\in\{1,2\}$ are obtained from $a_{kl}^{(j)}, b^{(j)}_1,\mu^{(j)}, F^{(j)}_k,  g^{(j)}_m$ through the change of variables $(r,x')\mapsto G_j(r,x')$, and thus they smoothly depend on $D\Psi, (\rho_-,p_-,\vphi_-),$ $D\phi_j$ and $DG_j(G_j^{-1})$.

By  \eqref{3-1-17}, \eqref{new-11}, \eqref{def-gs}, and by the definitions of $\mathfrak{F}_k^{(j)}, \mathfrak{g}_m^{(j)}$ in \eqref{new-at-17}, one can directly check
\begin{equation}
\label{new-at-12}
\begin{split}
&\|\mathfrak{F}^{(1)}-\mathfrak{F}^{(2)}\|_{\alp,\mcl{N}^+_{r_s}}^{(1-\alp,\corners)}\le C\kappa_*\delta(\phi_1,\phi_2)\\
&\|\mathfrak{g}_1^{(1)}-\mathfrak{g}_1^{(2)}\|_{\alp,\Lambda}^{(1-\alp,\der\Lambda)}\le C\kappa_*\delta(\phi_1,\phi_2)\\
&\|\mathfrak{g}_2^{(1)}-\mathfrak{g}_2^{(2)}\|_{\alp,\Gam_{w,r_s}}^{(1-\alp,\corners)}\le C\kappa_*\delta(\phi_1,\phi_2)\\
&\|\mathfrak{g}_3^{(1)}-\mathfrak{g}_3^{(2)}\|_{\alp,\Lambda}^{(1-\alp,\der\Lambda)}\le C(\|v_1-v_2\|_{\alp,\Lambda}^{(1-\alp,\der\Lambda)}+\kappa_*)\delta(\phi_1,\phi_2)
\end{split}
\end{equation}
where we set
\begin{equation*}
\begin{split}
&\kappa_*:=\varsigma_1+\varsigma_2
+\sum_{j=1}^2\|\phi_j\|_{1,\alp,\mcl{N}^+_{r_s}}^{(-\alp,\corners)},\\
&\delta(\phi_1,\phi_2):=\|\phi_1-\phi_2\|_{\alp,\mcl{N}^+_{r_s}}^{(1-\alp,\corners)}\\
&\phantom{aaaaa}
+\|DG_1(G_1^{-1})-DG_2(G_2^{-1})\|_{\alp,\mcl{N}^+_{r_s}}^{(1-\alp,\corners)}+\|G_1^{-1}-G_2^{-1}\|_{1,\alp,\mcl{N}^+_{r_s}}^{(-\alp,\corners)}.
\end{split}
\end{equation*}
Moreover, \eqref{new-at-9} and \eqref{new-at-10} imply
\begin{equation}
\label{new-at-13}
\delta(\phi_1,\phi_2)\le C(\|\phi_1-\phi_2\|_{1,\alp,\mcl{N}^+_{r_s}}^{(-\alp,\corners)}+
\|(\omega_1-\omega_2)\phi_2\|_{1,\alp,\mcl{N}^+_{r_s}}^{(-\alp,\corners)}).
\end{equation}
By Proposition \ref{theorem3-1} and \eqref{3-2-20}, each $\phi_j$ satisfies
\begin{equation}
\label{new-at-14}
\|\phi_j\|_{2,\alp,\mcl{N}^+_{r_s}}^{(-1-\alp,\corners)}\le C(\varsigma_1+\varsigma_2+\|v_j-v_c\|_{1,\alp,\Lambda}^{(-\alp,\der\Lambda)})
\end{equation}
for a constant $C$ depending on the data in sense of Remark \ref{remark-cons}. By \eqref{new-at-9}, \eqref{new-at-13} and \eqref{new-at-14}, there exists a constant $\sigma^{\flat}\in (0,\sigma_3]$ so that whenever $0<\varsigma_1,\varsigma_2, \varsigma_4^{(j)}\le \sigma^{\flat}$ are satisfied, there holds
\begin{equation}
\label{new-at-15}
\delta(\phi_1,\phi_2)\le C\|\phi_1-\phi_2\|_{1,\alp,\mcl{N}^+_{r_s}}^{(-\alp,\corners)}.
\end{equation}
Here, the constants $C$ and $\sigma^{\flat}$ depend only on the data in sense of Remark \ref{remark-cons}.
Then, applying Lemma \ref{prop7-1} to $\phi_1-\phi_2$, by \eqref{new-at-17}, \eqref{new-at-12} and \eqref{new-at-15}, we obtain \eqref{new-at-1} reducing $\sigma^{\flat}$ if necessary.
\quad\\
\textbf{Step 3.} By arguing similarly as \eqref{3-0-4}, and by the definition of $\phi_j$, setting $\psi_-=\vphi_--\vphi_0^-$, we have
\begin{equation}
\label{new-at-16}
\begin{split}
(f_1-f_2)(\vartheta)=\frac{[\phi_j(r_s,\vartheta)-\psi_-(f_j(\vartheta),\vartheta)]_{j=2}^1}
{\int_0^1\der_r(\vphi_0^--\vphi_0^+)(af_1(\vartheta)+(1-a)f_2(\vartheta))da}
\end{split}
\end{equation}
for $x\in\Lambda$.
By \eqref{new-at-16} and \eqref{new-at-14}, we obtain
\begin{equation}
\label{star4}
\begin{split}
\|f_1\!-\!f_2\|_{1,\alp,\Lambda}^{(\!-\alp,\der\Lambda\!)}
\!\!\le \!C\!(\|\!\phi_1\!-\!\phi_2\|_{1,\alp,\mcl{N}^+_{r_s}}^{(\!-\alp,\corners\!)}
\!+\!(\!\varsigma_2\!+\!\varsigma_4^{(1)}\!+\!\varsigma_4^{(2)}\!)\|f_1\!-\!f_2\|_{1,\alp,\Lambda}^{(\!-\alp,\der\Lambda\!)}\!).
\end{split}
\end{equation}
We reduce $\sigma^{\flat}$ further to obtain \eqref{new-at-2} by  \eqref{new-29} and \eqref{star4}.
\end{proof}
\end{corollary}
For any given $\alp\in(\frac 12,1)$, let $\psi\in C^2(\mcl{N}^+_{r_s})\cap C^{\alp}(\ol{\mcl{N}^+_{r_s}})$ be be a unique solution to \eqref{new-at-17} with $\mathfrak{F}_k=\mathfrak{g}_1=\mathfrak{g}_2=0$ and $\mathfrak{g}_3=v_1-v_2$.
\begin{corollary}
\label{corollary3-lemma6.1}Fix $\alp\in (\frac 12,1)$.
Let $\phi_1$ and $\phi_2$ be as in Corollary \ref{corollary2-lemma6.1}, Then, there are constants $C,\sigma^{\natural}>0$ depending only on the data in sense of Remark \ref{remark-cons} so that whenever $0<\varsigma_1,\varsigma_2,\varsigma_4^{(1)},\varsigma_4^{(2)}\le \sigma^{\natural}$ are satisfied, then there holds,
\begin{equation}
\label{new-at-3}
\|\phi_1-\phi_2-\psi\|_{1,\alp,S_0}^{(-\alp,\der S_0\cup\der\Gam_{ex})}\le C(\varsigma_1+\varsigma_2+\varsigma_4^{(1)}+\varsigma_4^{(2)})\|v_1-v_2\|_{\alp,\Lambda}^{(1-\alp,\der\Lambda)}.
\end{equation}
\begin{proof}
By \eqref{new-at-17}, and by the definition of $\psi$, $\psi_*:=\phi_1-\phi_2-\psi$ satisfies
So $\phi_1-\phi_2-u$ satisfies
\begin{align*}
&\der_k(\til a_{kl}(y,0)\der_l\psi_*)=\der_k(\mathfrak{F}^{(1)}-\mathfrak{F}^{(2)})\;\;\text{in}\;\;\mcl{N}^+_{r_s}\\
&\der_r\psi_*-\til\mu_0\psi_*=\mathfrak{g}^{(1)}_1-\mathfrak{g}^{(2)}_1\;\;\text{on}\;\;S_0\\
&(\til a_{kl}(y,0)\der_l\psi_*)\cdot\nu_w=\mathfrak{g}^{(1)}_2-\mathfrak{g}^{(2)}_2\;\;\text{on}\;\;\Gam_{w,r_s}\\
&(\til a_{kl}(y,0)\der_l\psi_*)\cdot\hat r=\mathfrak{g}^{(1)}_3-\mathfrak{g}^{(2)}_3-(v_1-v_2)\;\;\text{on}\;\;\Gam_{ex}.
\end{align*}
By the definition of $\mathfrak{g}_3^{(j)}$ in \eqref{new-at-17}, we emphasize that $\mathfrak{g}_3^*:=\mathfrak{g}^{(1)}_3-\mathfrak{g}^{(2)}_3-(v_1-v_2)$ satisfies
\begin{equation*}
\|\mathfrak{g}_3^*\|_{\alp,\Lambda}^{(1-\alp,\der\Lambda)}
\le C(\varsigma_1+\varsigma_2+\varsigma_4^{(1)}+\varsigma_4^{(2)})\delta(\phi_1,\phi_2)
\end{equation*}
where $\delta(\phi_1,\phi_2)$ is defined in \eqref{new-at-13}. Then, by Lemma \ref{prop7-1}, Corollary \ref{corollary2-lemma6.1} and \eqref{new-at-12}, we conclude that if $\sigma^{\natural}$ is chosen sufficiently small depending on the data in sense of Remark \ref{remark-cons}, then $\psi_*=\phi_1-\phi_2-\psi$ satisfies
\begin{equation}
\label{new-at-18}
\|\psi_*\|_{1,\alp,\mcl{N}^+_{r_s}}^{(-\alp,\corners)}\le C(\varsigma_1+\varsigma_2+\varsigma_4^{(1)}+\varsigma_4^{(2)})\|v_1-v_2\|_{\alp,\Lambda}^{(1-\alp,\der\Lambda)}
\end{equation}
for a constant $C$ depending only on the data in sense of Remark \ref{remark-cons}.
\end{proof}
\end{corollary}

Fix $(\Psi,\vphi_-,p_-)$ in $B_{\sigma}^{(1)}(Id,\vphi_0^-,p_0^-)$. For simplicity, let us denote as $\mcl{P}(v)$ for $\mcl{P}(\Psi,\vphi_-,p_-,v)$ .

\begin{proof}[\textbf{Proof of Lemma \ref{lemma7-1}}] \textbf{Step 1.}For convenience, for each $j=1,2$, let us denote as $\mcl{P}_j, \mcl{Q}_j$ and $\mcl{R}_j$ for $\mcl{P}(\zeta_*,v_j), \mcl{Q}(\zeta_*,v_j)$ and $\mcl{R}(\zeta_*,v_j)$ respectively where $\mcl{Q}, \mcl{R}$ are defined in \eqref{5-2}. Then we can write

\begin{equation}
\label{new-12}
\mcl{P}_1-\mcl{P}_2-D_v\mcl{P}(v_1-v_2)=l_1\mcl{R}_1+l_2\mcl{Q}_2+l_3
\end{equation}
with
\begin{align*}
&l_1=\mcl{Q}_1-\mcl{Q}_2-D_v\mcl{Q}(v_1-v_2),\quad l_2=\mcl{R}_1-\mcl{R}_2-D_v\mcl{R}(v_1-v_2),\\
&l_3=(\mcl{Q}_2-\mcl{Q}_0)D_v\mcl{R}(v_1-v_2)+(\mcl{R}_1-\mcl{R}_0)D_v\mcl{Q}(v_1-v_2),
\end{align*}
where we denote $\mcl{Q}_0=\mcl{Q}(Id,\vphi_0^-,p_0^-,v_c)$, $\mcl{R}_0=\mcl{R}(Id,\vphi_0^-,p_0^-,v_c)$. We estimate $l_1\mcl{R}_1, l_2\mcl{Q}_2, l_3$ separately.

By Lemma \ref{corollary1-lemma6.1}, Proposition \ref{theorem3-1}, \eqref{new-13} and \eqref{new-14}, we obtain
\begin{equation}
\label{new14}
\|l_3\|_{\alp,\Lambda}^{(1-\alp,\der\Lambda)}\le C(\varsigma_1+\varsigma_2+\varsigma_4^{(1)}+\varsigma_4^{(2)})\|v_1-v_2\|_{\alp,\Lambda}^{(1-\alp,\der\Lambda)}.
\end{equation}

\textbf{Step 2.}
Let $\psi\in C^2(\mcl{N}^+_{r_s})\cap C^{\alp}(\ol{\mcl{N}^+_{r_s}})$ be a unique solution to \eqref{5-13a},\eqref{5-13d} with $w=v_1-v_2$, then by \eqref{l1} and the definition of $D_v\mcl{Q}$ in \eqref{5-15}, we have
\begin{equation}
\label{newex-l2}
l_1(\vartheta)
\!=\!(\!E_1\!-\!E_2-\!a_2\psi\!)|_{X(2r_1\!-\!r_s;r_1,\vartheta)}\!+\!a_2(\psi|_{X(2r_1\!-\!r_s;r_1,\vartheta)}\!-\!\psi(r_s,\vartheta)\!)\!
+\!I_*(\vartheta)
\end{equation}
where $a_2, I_*$ are defined in \eqref{5-A2}, \eqref{l1} respectively.
First of all, by \eqref{new-13}, \eqref{new-14}, Corollary \ref{corollary2-lemma6.1} and \eqref{new-at-9}, we obtain

\begin{equation}
\label{new-15}
\|I_*\|_{\alp,\Lambda}^{(1-\alp,\der\Lambda)}
\le C(\varsigma_1+\varsigma_2+\varsigma_4^{(2)})\|v_1-v_2\|_{\alp,\Lambda}^{(1-\alp,\der\Lambda)}.
\end{equation}

Secondly, for $V_0$ defined in Lemma \ref{lemma-V}, let us set $W^*_0:=\frac{V_0}{V_0\cdot\hat r}$ so that we have $W_0^*=(1,0,\cdots,0)$ in $(r,\vartheta)$ coordinates which is a spherical coordinate system we specified in section \ref{sec:nozzle3}. Then, by Lemma \ref{lemma-V}, there holds
\begin{equation}
\label{star5}
\|W^*_1-W^*_0\|_{1,\alp,\mcl{N}^+_{r_s}}^{(-\alp,\corners)}\le C(\varsigma+\varsigma_2+\varsigma_4^{(2)})
\end{equation}
for a constant $C$ depending on the data in sense of Remark \ref{remark-cons} where $W^*_1$ is defined in the paragraph of \eqref{l1}. Thus, by \eqref{star5} and Lemma \ref{prop7-1}, we obtain
\begin{equation}
\label{new-18}
\begin{split}
&\phantom{\le}\|\psi|_{X_1(2r_1-r-s;r_1,\cdot)}-\psi(r_s,\cdot)\|_{\alp,\Lambda}^{(1-\alp,\der\Lambda)}\\
&\le C\|D\psi\|_{\alp,\mcl{N}^+_{e_s}}^{(1-\alp,\corners)}\|W^*_1-W^*_0\|_{1,\alp,\mcl{N}^+_{r_s}}^{(-\alp,\corners)}\\
&\le C(\varsigma_1+\varsigma_2+\varsigma_4^{(1)})\|v_1-v_2\|_{\alp,\Lambda}^{(1-\alp,\der\Lambda)}
\end{split}
\end{equation}
for a constant $C$ depending only on the data.

To complete step 2, it remains to estimate $\|(E_1-E_2-a_2\psi)(r_s,\cdot)\|_{\alp,\Lambda}^{(1-\alp,\der\Lambda)}$ because $X_1(2r_1-r_s;r_1,\vartheta)$ lies on $S_0=\{r=r_s\}\cap \mcl{N}$ for all $\vartheta\in\Lambda$.

By the definition of $a_2$ in \eqref{5-A2}, $a_2\psi$ is expressed as
\begin{equation}
\label{new-19}
a_2\psi(r_s,\cdot)=a_2^{(1)}\psi(r_s,\cdot)+a_2^{(2)}\psi(r_s,\cdot)
\end{equation}
with
\begin{equation}
\label{a2s}
a_2^{(1)}=\frac{1}{(B_0-\frac 12(\der_r\vphi_0^+(r_s))^2)^{\frac{\gam}{\gam-1}}}\frac{\der_r(p_{s,0}-p_0^+)(r_s)}{\der_r(\vphi_0^--\vphi_0^+)(r_s)},\quad
a_2^{(2)}=p_0^+(r_s)a_3
\end{equation}
where $a_3$ is defined in \eqref{derivative2}. So, we need to decompose $E_1-E_2$ into two parts so that one is comparable to $a_2^{(1)}\psi$, and the other is comparable to $a_2^{(2)}\psi$.

For each $j=1,2$, by \eqref{4-0-5}, we have
\begin{equation}
\label{new-29}
E_j(r_s,\vartheta)=\frac{\rho_-(M_j\nabla\vphi_-\cdot\nu_s^{(j)})^2+p_--\rho_-K_s^{(j)}}{(B_0-\frac 12|D\Psi^{-1}(\Psi)\nabla\vphi_j|^2)^{\frac{\gam}{\gam-1}}}=:\frac{p_{s}^{(j)}}{U_j}
\end{equation}
with
\begin{align*}
&M_j=\frac{|\nabla(\vphi_--\vphi_j)|(D\Psi^{-1})^T(D\Psi^{-1})}{|D\Psi^{-1}\nabla(\vphi_--\vphi_j)|}\\
&K_s^{(j)}=\frac{2(\gam-1)}{\gam+1}(\frac 12(M_j\nabla\vphi_-\cdot\nu_s^{(j)})^2+\frac{\gam p_-}{(\gam-1)\rho_-})
\end{align*}
where every quantity is evaluated at $(f_j(\vartheta),\vartheta)$. Here, $\vphi_j$ and $f_j$ are as in Corollary \ref{corollary2-lemma6.1}, and $\nu_s^{(j)}$ is the unit normal on $S_j=\{r=f_j(\vartheta)\}$ toward the subsonic domain $\mcl{N}^+_{f_j}$ of $\vphi_j$ for $j=1,2$.

Set
$
U_0(r):=(B_0-\frac 12|\der_r\vphi_0^+(r)|^2)^{\frac{\gam}{\gam-1}}.
$ Then, for each $j=1,2$, we can write $E_j(r_s,\vartheta)-E_0^+(f_j(\vartheta))=\frac{p_s^{(j)}}{U_j}-\frac{p_{s,0}}{U_0}+\frac{p_0^+-p_{s,0}}{U_0}$
where $p_{s,0}, E_0^+$ is defined in \eqref{2-2-6} and Lemma \ref{lemma-ent-0} respectively, and where all the quantities are evaluated at $(f_j(\vartheta),\vartheta)$.
%\begin{equation}
%\label{U0}
%p_{s,0}(r)=(\rho_0^-(\der\vphi_0^-)^2+p_0^--\rho_0^-K_0)(r),\;\;
%U_0(r):=(B_0-\frac 12|D\vphi_0^+(r)|^2)^{\frac{\gam}{\gam-1}},
%\end{equation}
Then, by Lemma \ref{lemma-ent-0}, we can express $(E_1-E_2)(r_s,\vartheta)$ as
\begin{align}
\label{new-20}
&(E_1-E_2)(r_s,\vartheta)=z_1+z_2,\\
\label{new-21}
&z_1=[\frac{p_s^{(j)}}{U_j}-\frac{p_{s,0}(f_{j}(\vartheta))}{U_0(f_j(\vartheta))}]_{j=2}^1,\;\;
z_2=[\frac{(p_{s,0}-p_0^+)(f_j(\vartheta))}{U_0(f_j(\vartheta))}]_{j=2}^1.
\end{align}
Here, $[F_j]_{j=2}^1$ is defined by $[F_j]_{j=2}^1:=F_1-F_2$.

We claim that, for a constant $C$ depending on the data in sense of Remark \ref{remark-cons}, there hold
\begin{align}
\label{new-22}
&\|z_1-a_2^{(2)}\psi(r_s,\cdot)\|_{\alp,\Lambda}^{(1-\alp,\der\Lambda)}
\le C\varsigma_*\|v_1-v_2\|_{\alp,\Lambda}^{(1-\alp,\der\Lambda)},\\
\label{new-24}
&\|z_2-a_2^{(1)}\psi(r_s,\cdot)\|_{\alp,\Lambda}^{(1-\alp,\der\Lambda)}
\le C\varsigma_*\|v_1-v_2\|_{\alp,\Lambda}^{(1-\alp,\der\Lambda)}.
\end{align}
where we set $\varsigma_*:=\varsigma_1+\varsigma_2+\varsigma_4^{(1)}+\varsigma_4^{(2)}$.

\textbf{Step 3.}\textbf{(Verification of \eqref{new-24})} We rewrite $z_2$ in \eqref{new-21} as
\begin{equation}
\label{new-23}
\begin{split}
&\phantom{=}z_2(\vartheta)\\
&=\frac{[(p_{s,0}\!-\!p_0^+)|_{r=f_j(\vartheta)}]_{j=2}^1}{U_0(r_s)}\!+\!
[(p_{s,0}\!-\!p_0^+)|_{r=f_j(\vartheta)}(\frac {1}{U_0(f_j(\vartheta))}\!-\!\frac{1}{U_0(r_s)})]_{j=2}^1\\
&=:z_2^{(1)}+z_2^{(2)}.
\end{split}
\end{equation}
By \eqref{2-2-6}, Proposition \ref{theorem3-1} and Corollary \eqref{corollary2-lemma6.1}, we obtain

\begin{equation}
\label{new-27}
\|z_2^{(2)}\|_{\alp,\Lambda}^{(1-\alp,\der\Lambda)}\le C\varsigma_*\|v_1-v_2\|_{\alp,\Lambda}^{(1-\alp,\der\Lambda)}.
\end{equation}
By \eqref{new-at-16}, we can write
\begin{equation*}
z_2^{(1)}(\vartheta)-a_2^{(1)}\psi(r_s,\vartheta)
=\frac{\til a_2^{(1)}(\vartheta)[\phi_j(r_s,\vartheta)-\psi_-(f_j(\vartheta),\vartheta)]_{j=2}^1-a_2^{(1)}\psi_(r_s,\vartheta)}{U_0(r_s)}
\end{equation*}
with
\begin{equation*}
\til a_2^{(1)}(\vartheta)
=\frac{\int_0^1\der_r(p_{s,0}-p_0^+)(af_1(\vartheta)+(1-a)f_2(\vartheta))da}{U_0(r_s)\int_0^1\der_r(\vphi_0^--\vphi_0^+)(af_1(\vartheta)+(1-a)f_2(\vartheta))da}.
\end{equation*}
where each $\phi_j$ is as in Corollary \ref{corollary2-lemma6.1}, and where we denote as $\psi_-=\vphi_--\vphi_0^-$. By Corollary \ref{corollary2-lemma6.1} and Corollary \ref{corollary3-lemma6.1}, there holds
\begin{equation}
\label{new-26}
\|z_2^{(1)}-a_2^{(1)}\psi(r_s,\cdot)\|_{\alp,\Lambda}^{(1-\alp,\der\Lambda)}
\le C\varsigma_*\|v_1-v_2\|_{\alp,\Lambda}^{(1-\alp,\der\Lambda)}.
\end{equation}
The constants $C$ in \eqref{new-27}, \eqref{new-26} depend only on the data in sense of Remark \ref{remark-cons}. Then, \eqref{new-24} is obtained by \eqref{new-27} and \eqref{new-26}.

\textbf{Step 4. (Verification of \eqref{new-22})} By \eqref{5-13a} and \eqref{derivative3}, it is easy to see that for $a_2^{(2)}\psi(r_s,\cdot)$ in \eqref{new-19}, there holds
\begin{equation}
\label{jewel1}
a_2^{(2)}\psi(r_s,\cdot)=\frac{\gam p_0^+(r_s)\der_r\vphi_0^+(r_s)}{(\gam-1)(B_0-\frac 12(\der_r\vphi_0^+(r_s))^2)^{\frac{2\gam-1}{\gam-1}}}\der_r\psi(r_s,\cdot)=a_4\der_r\psi(r_s,\cdot).
\end{equation}

Similarly to step 3 above, we rewrite $z_1$ in \eqref{new-21} as
\begin{equation}
\label{jewel2}
\begin{split}
z_1&=p_0^+(r_s)[\frac{1}{U_j^*}-\frac{1}{U_0(f_j(\vartheta))}]_{j=2}^1
+[(\frac{p_s^
{(j)}}{U_j}-\frac{p_0^+}{U_j^*})+\frac{p_0^+(r_s)-p_{s,0}(f_j(\vartheta))}{U_0(f_j(\vartheta))}]_{j=2}^1\\
&=:z_1^{(1)}+z_1^{(2)}
\end{split}
\end{equation}
where we set
\begin{equation*}
U_j^*:=(B_0-\frac 12|\nabla\vphi(f_j(\vartheta),\vartheta)|^2)^{\frac{\gam}{\gam-1}}.
\end{equation*}

By \eqref{2-2-6}, Proposition \ref{theorem3-1}, \eqref{new-29} and Corollary \ref{corollary2-lemma6.1}, $z_1^{(2)}$ in \eqref{jewel2} satisfies
\begin{equation}
\label{new-30}
\|z_1^{(2)}\|_{\alp,\Lambda}^{(1-\alp,\der\Lambda)}\le C\varsigma_*\|v_1-v_2\|_{\alp,\Lambda}^{(1-\alp,\der\Lambda)}.
\end{equation}

For each $j=1,2$, we can write
\begin{equation}
\label{jewel3}
\begin{split}
&\phantom{=}p_0^+[\frac{1}{U_j^*}-\frac{1}{U_0(f_j(\vartheta))}]\\
&=\!\int_0^1\!\!\frac{\gam p_0^+(r_s)\beta_4^{(j)}(\vartheta)\!\cdot\! \nabla(\vphi_j\!-\!\vphi_0^+)|_{(f_j(\vartheta),\vartheta)}}
{(\gam\!-\!1)(B_0\!-\!\frac 12(|\der_r\vphi_0^+(f_j(\vartheta))|^2\!+\!2t\beta_4^{(j)}\!\cdot\!\nabla\!(\vphi_j\!-\!\vphi_0^+)|_{(f_j(\vartheta),\vartheta)}))^{\frac{2\gam-1}{\gam-1}}}dt\\
&=a_4^{(j)}\cdot \nabla\!(\vphi_j\!-\!\vphi_0^+)|_{(f_j(\vartheta),\vartheta)}
\end{split}
\end{equation}
with $\beta_4^{(j)}(\vartheta)=\der_r\vphi_0^+(f_j(\vartheta))\hat r+\frac 12\nabla(\vphi_j-\vphi_0^+)|_{(f_j(\vartheta),\vartheta)}$.

By Proposition \ref{theorem3-1}, we have, for a constant $C$ depending on the data in sense of Remark \ref{remark-cons}
\begin{equation}
\label{jewel4}
\|a_4^{(j)}-a_4\hat r\|_{1,\alp,\Lambda}^{(-\alp,\der\Lambda)}\le C\varsigma_*.
\end{equation}
So, if we write $z_1^{(1)}-a_4\der_r\psi(r_s,\cdot)$ as
\begin{equation*}
\begin{split}
&\phantom{=}z_1^{(1)}-a_4\der_r\psi(r_s,\vartheta)\\
&=\!a_4([\der_r(\vphi_j\!-\!\vphi_0^+)|_{(f_j(\vartheta,\vartheta))}]_{j=2}^1\!-\!\psi(r_s,\!\vartheta)\!)
+[(a_4^{(j)}\!-\!a_4)\!\cdot\!\nabla(\vphi_j\!-\!\vphi_0^+)|_{(f_j(\vartheta),\vartheta)}]_{j=2}^1,
\end{split}
\end{equation*}
then, by Proposition \ref{theorem3-1}, \eqref{3-2-20}, Corollary \ref{corollary2-lemma6.1}, Corollary \ref{corollary3-lemma6.1}, \eqref{jewel4}, there is a constant $C$ depending on the data in sense of Remark \ref{remark-cons} so that there holds
\begin{equation}
\label{jewel5}
\|z_1^{(1)}-a_4\der_r\psi(r_s,\vartheta)\|_{\alp,\Lambda}^{(1-\alp,\der\Lambda)}\le C\varsigma_*\|v_1-v_2\|_{\alp,\Lambda}^{(1-\alp,\der\Lambda)}.
\end{equation}
Thus \eqref{new-22} holds true.  Finally, \eqref{new-22} and \eqref{new-24} imply, for a constant $C$ depending on the data in sense of Remark \ref{remark-cons},
\begin{equation}
\label{new-32}
\|\mcl{R}_1l_1\|_{\alp,\Lambda}^{(1-\alp,\der\Lambda)}\le
C\varsigma_*\|v_1-v_2\|_{\alp,\Lambda}^{(1-\alp,\der\Lambda)}.
\end{equation}
Similarly, one can also show
\begin{equation}
\label{new-33}
\|\mcl{Q}_2l_2\|_{\alp,\Lambda}^{(1-\alp,\der\Lambda)}\le
C\varsigma_*\|v_1-v_2\|_{\alp,\Lambda}^{(1-\alp,\der\Lambda)}.
\end{equation}
As \eqref{new-33} can be checked more simply than \eqref{new-32}, we omit the details for the verification of \eqref{new-33}.

Finally, combining \eqref{new14}, \eqref{new-32} and \eqref{new-33} all together, we obtain \eqref{new-34}, and this competes the proof of Lemma \ref{lemma7-1}.

\end{proof}

\appendix
\section{Proof of Proposition \ref{lemma5-1}}\label{appendix A}
 Before we prove Proposition \ref{lemma5-1}, we first prove the following lemma.

\begin{lemma}
[Right Inverse Mapping Theorem]
\label{theoremB-1}
Let $\mathfrak{C}_1$ be a compactly imbedded subspace of a Banach space $\mathfrak{B}_1$. Also, suppose that  $\mathfrak{C}_2$ is a subspace of another Banach space $\mathfrak{B}_2$.  For a given point $(x_0,y_0)\in\mathfrak{C}_1\times \mathfrak{C}_2$,  suppose that a mapping $\mcl{H}$ satisfies the following conditions:
\begin{itemize}
\item[(i)]  $\mcl{H}$ maps a neighborhood of $x_0$ in $\mathfrak{B}_1$ to $\mathfrak{B}_2$, and maps a neighborhood of $x_0$ in $\mathfrak{C}_1$ to  $\mathfrak{C}_2$, with $\mcl{H}(x_0)=y_0$,
\item[(ii)]  whenever a sequence $\{x_k\}\subset\mathfrak{C}_1$ near $x_0$ converges to $x_*$ in $\mathfrak{B}_1$, the sequence $\{\mcl{H}(x_k)\}\subset\mathfrak{C}_2$ converges to $\mcl{H}(x_*)$ in $\mathfrak{B}_2$,
\item[(iii)] $\mcl{H}$,as a mapping from $\mathfrak{B}_1$ to $\mathfrak{B}_2$, also as a mapping from $\mathfrak{C}_1$ to $\mathfrak{C}_2$,  is Fr\'{e}chet differentiable at $x_0$,
\item[(iv)] the Fr\'{e}chet derivative $D\mcl{H}(x_0)$, as a mapping from $\mathfrak{B}_1$ to $\mathfrak{B}_2$, also as a mapping from $\mathfrak{C}_1$ to $\mathfrak{C}_2$, is invertible.
\end{itemize}
Then there is a small neighborhood $\mcl{V}(x_0)$ of $x_0$ in $\mathfrak{C}_1$, and a small neighborhood $\mcl{W}(y_0)$ of $y_0$ in $\mathfrak{C}_2$ so that $\mcl{H}:\mcl{V}(x_0)\to\mcl{W}(y_0)$ has a right inverse i.e., there is a mapping $\mcl{H}^{-1}_{Right}:\mcl{W}(y_0)\to\mcl{V}(x_0)$ satisfying
\begin{equation}
\label{right}
\mcl{H}\circ \mcl{H}_{Right}^{-1}=id
\end{equation}
and moreover $\mcl{H}_{Right}^{-1}$ satisfies $\mcl{H}_{Right}^{-1}\circ \mcl{H}(x_0)=x_0$.
\begin{proof}
Denote $D\mcl{H}(x_0)$ as $D\mcl{H}_0$. For $t\in\R_+$, and $v,w\in \mathfrak{C}_1$, define
\begin{equation*}
\mcl{H}^*(t,v,w):=\frac 1t\bigl(\mcl{H}(x_0+t(v+w))-\mcl{H}(x_0)\bigr)-D\mcl{H}_0 w.
\end{equation*}
For $r>0$ and $a_0\in\mathfrak{C}_1$, set
$
B_r(a_0):=\{y\in\mathfrak{C}_1:\|y-a_0\|_{\mathfrak{C}_1}\le r\}.
$ For a fixed $(t,w)\in (0,\eps_0]\times \der B_1({{0}})$, we define
\begin{equation*}
\mcl{H}^*_{t,w}(v):=D\mcl{H}_0^{-1}(\mcl{H}^*)(t,v,w)-v.
\end{equation*}
If the constant $\eps_0$ is sufficiently small then we have $\mcl{H}^*_{t,w}$ maps $B_1(0)$ into itself. Since $B_1(0)$ is bounded in $\mathfrak{C}_1$, ${B}_1(0)$ is a compact convex subset of $\mathfrak{B}_1$.
By Lemma \ref{theoremB-1} (ii) and (iii), the Schauder fixed point theorem applies to $\mcl{H}^*_{t,w}$. Thus $\mcl{H}^*_{t,w}$  has a fixed point $v_*=v_*(t,w)$ satisfying
$$\mcl{H}^*_{t,w}(v_*)=v_*.$$

By the definition of $\mcl{H}^*_{t,w}$, $v_*$ satisfies
\begin{equation}
\label{B-1}
\begin{split}
\mcl{H}^*_{t,w}(v_*)=v_* &\Leftrightarrow \mcl{H}^*(t,v_*,w)=0\\
&\Leftrightarrow \mcl{H}(x_0+t(v_*+w))=\mcl{H}(x_0)+tD\mcl{H}_0w.
\end{split}
\end{equation}
We claim
\begin{equation}
\label{B-2}
\mcl{H}_{Right}^{-1}(z):=x_0+t(v_*+w)
\end{equation}
for $z=\mcl{H}(x_0)+tD\mcl{H}_0w$. Obviously, $\mcl{H}_{Right}^{-1}(z)$ in \eqref{B-2} satisfies \eqref{right}. So it remains to check if $\mcl{H}_{Right}^{-1}$  in\eqref{B-2} is well-defined.

Set $\mcl{W}(y_0):=\{y_0+tD\mcl{H}_0w:t\in[0,\eps_0],w\in \der B_1({{0}})\}$. Since $D\mcl{H}_0$ is invertible, $\mcl{W}(y_0)$ covers a small neighborhood of $y_0$ in $\mathfrak{C}_2$. Suppose $
y_0+t_1D\mcl{H}_0 w_1=y_0+t_2D\mcl{H}_0 w_2$
for $t_1,t_2\in[0,\eps_0]$ and $w_1,w_2\in \der B_1({0})$. If $t_1=0$ then $t_2$ must be $0$ as well. If $t_1\neq 0$ then we have $D\mcl{H}_0w_1=\frac{t_2}{t_1}D\mcl{H}_0w_2$ and this implies $w_1=\frac{t_2}{t_1}w_2$. So we obtain $t_1=t_2$ and $w_1=w_2$.

By (\ref{B-2}), we also have
$
\mcl{H}^{-1}_{Right}\circ\mcl{H}(x_0)=\mcl{H}_{Right}^{-1}(\mcl{H}(x_0)+0w)=x_0.
$
\end{proof}
\end{lemma}
\begin{remark}
\label{remark-norm}
By \eqref{B-2}, there holds
\begin{equation*}
\|\mcl{H}^{-1}(z)-x_0\|_{\mathfrak{C}_1}\le 2t=2\frac{\|z-\mcl{H}(x_0)\|_{\mathfrak{C}_2}}{\|D\mcl{H}_0w\|_{\mathfrak{C}_2}}
\le 2\|D\mcl{H}_0^{-1}\|\|z-\mcl{H}(x_0)\|_{\mathfrak{C}_2}
\end{equation*}
where we set $\|D\mcl{H}_0^{-1}\|:=\sup_{\|v\|_{\mathfrak{C}_2}=0}\|D\mcl{H}_0^{-1}\|_{\mathfrak{C}_1}$.

\end{remark}
Proposition \ref{lemma5-1} can be easily proven from Lemma \ref{theoremB-1}.
\begin{proof}
[\textbf{Proof of Proposition \ref{lemma5-1}}]
Let us define $\mcl{H}(x,y)$ by
\begin{equation}
\label{V}
\mcl{H}(x,y):=(\mcl{F}(x,y),y).
\end{equation}
In Lemma \ref{theoremB-1}, replacing $\mathfrak{B}_1$ by $\mathfrak{B}_1\times \mathfrak{B}_2$ , $\mathfrak{C}_1$ by $\mathfrak{C}_1\times\mathfrak{C}_2$, $\mathfrak{B}_2$ by $\mathfrak{B}_3$,  $\mathfrak{C}_2$ by $\mathfrak{C}_3$, $x_0$ by $(x_0,y_0)$ and $y_0$ by $({\mathbf{0}},y_0)$, $\mcl{H}$ in \eqref{V} satisfies Lemma \ref{theoremB-1} (i)-(iii).
In particular, the Fr\'{e}chet derivative of $\mcl{H}$ at $(x_0,y_0)$, which we write as $D\mcl{H}_0$, is given by
\begin{equation*}
D\mcl{H}_0(x,y)=(D_x\mcl{F}_0x+D_y\mcl{F}_0y,y)
\end{equation*}
where we denote $D_{(x,y)}\mcl{F}(x_0,y_0)$ as $(D_x\mcl{F}_0,D_y\mcl{F}_0)$. By Proposition \ref{lemma5-1} (iv), $\mcl{H}$ satisfies Lemma \ref{theoremB-1} (iv). Therefore, by Lemma \ref{theoremB-1}, $\mcl{H}$ has a right inverse $\mcl{H}_{Right}^{-1}$ in a small neighborhood of $(x_0,y_0)$ in $\mathfrak{C}_1\times\mathfrak{C}_2$. In particular, for $y\in\mathfrak{C}_2$ with $\|y-y_0\|_{\mathfrak{C}_2}$ being sufficiently small, $\mcl{H}_{Right}^{-1}$ is well defined. Write $\mcl{H}_{Right}^{-1}({{0}},y)=(x^*(y),y)$ then we obtain
$
\mcl{F}(x^*(y),y)={{0}}.
$
\end{proof}

\section{Proof of Lemma \ref{prop7-1}}\label{appendix B}
\begin{proof}[\textbf{Proof of Lemma \ref{prop7-1}}]\quad\\
\textbf{Step 1.} Fix $F^*=(F_j)_{j=1}^n\in C^{\alp}_{(1-\alp,\corners)}(\mcl{N}^+_{r_s},\R^n)$, and $g_1^*, g_3^*\in C^{\alp}_{(1-\alp,\der\Lambda)}(\Lambda)$, $g^*_2\in C^{\alp}_{(1-\alp,\corners)}(\Gam_{w,r_s})$. Then, there are sequences $\{F^{(k)}=(F^{(k)}_j)_{j=1}^n\}\subset C^{\alp}(\ol{\mcl{N}^+_{r_s}},\R^n),$ and $\{g^{(k)}_1\}, \{g^{(k)}_3\}\subset C^{\alp}(\ol{\Lambda})$, $\{g^{(k)}_2\}\subset C^{\alp}(\ol{\Gam_{w,r_s}})$ satisfying
\begin{equation*}
\begin{split}
&\lim_{k\to \infty}\|F^{(k)}-F^*\|_{\alp,\mcl{N}^+_{r_s}}^{(1-\alp,\corners)}=\lim_{k\to \infty}\|g_1^{(k)}-g_1^*\|_{\alp,\Lambda}^{(1-\alp,\der\Lambda)}\\
&=\lim_{k\to \infty}\|g_2^{(k)}-g_2^*\|_{\alp,\Gam_{w,r_s}}^{(1-\alp,\corners)}=
\lim_{k\to\infty}\|g^{(k)}_3-g_3^*\|_{\alp,\Lambda}^{(1-\alp,\der\Lambda)}=0.
\end{split}
\end{equation*}
For each $k$, let $\psi^{(k)}\in C^2(\mcl{N}^+_{r_s})\cap C^0(\ol{\mcl{N}^+_{r_s}})$ be the unique solution to \eqref{bd1}-\eqref{bd4} with $F=F^{(k)}, g_l=g^{(k)}_l$ for $l=1,2,3$. The unique existence of $\psi^{(k)}$ can be proven by the method of continuity as in the proof of Lemma \ref{lemma-linear-existence}.

 If \eqref{weak-est-gen} holds true for $F\in C^{\alp}(\ol{\mcl{N}^+_{r_s}},\R^n), g_1, g_2\in C^{\alp}(\ol{\Lambda}), g_2\in C^{\alp}(\ol{\Gam_{w,r_s}})$, then we have
\begin{equation}
\label{new-at-19}
\lim_{k,m\to\infty}\|\psi^{(k)}-\psi^{(m)}\|_{1,\alp,\mcl{N}^+_{r_s}}^{(-\alp,\corners)}=0.
\end{equation}
By \eqref{new-at-19}, $\psi^{(k)}$ converges to a function $\psi^*$ in $C^2$ for any compact subset of $\mcl{N}^+_{r_s}$, so $\psi^*$ satisfies \eqref{bd1} with $F=F^*$ in $\mcl{N}^+_{r_s}$. Also, $\psi^*$ satisfies \eqref{bd2}-\eqref{bd4} on the relative interior of each corresponding boundary with $g_l=g^*_l$ for $l=1,2,3$. Since $\{\psi^{(k)}\}$ is a Cauchy sequence in $C^{1,\alp}_{(-\alp,\corners)}(\ol{\mcl{N}^+_{r_s}})$, we have $\psi^*\in C^{1,\alp}_{(-\alp,\corners)}(\mcl{N}^+_{r_s})\subset C^{\alp}(\ol{\mcl{N}^+_{r_s}})$. Then, by \cite[Corollary 2.5]{Lie-1}, $\mu_0>0$ in \eqref{bd2} implies that $\psi^*\in C^2(\mcl{N}^+_{r_s})\cap C^0(\ol{\mcl{N}}^+_{r_s})$ is a unique solution to \eqref{bd1}-\eqref{bd4}. Moreover, $\psi^*$ satisfies \eqref{weak-est-gen} by \eqref{new-at-19}. Thus, it suffices to prove \eqref{weak-est-gen} for
\begin{equation}
\label{new-at-20}
F\in C^{\alp}(\ol{\mcl{N}^+_{r_s}},\R^n),\;  g_1, g_2\in C^{\alp}(\ol{\Lambda}),\;  g_2\in C^{\alp}(\ol{\Gam_{w,r_s}}).
\end{equation}
So, for the rest of proof, we assume \eqref{new-at-20}.

\textbf{Step 2.} By the local scalings as in the proof of Lemma \ref{lemma-linear-existence}, one can easily check that if $\psi\in C^2(\mcl{N}^+_{r_s})\cap C^0(\ol{\mcl{N}^+_{r_s}})$ solves \eqref{bd1}-\eqref{bd4}, then there is a constant $C$ depending on $\mcl{N}^+_{r_s}$ and $\alp$ so that there holds
\begin{equation}
\label{new-dec-16}
\|\psi\|_{1,\alp,\mcl{N}^+_{r_s}}^{(-\alp,\corners)}\le C(|\psi|_{\alp,\mcl{N}^+_{r_s}}+K(\psi, F, g_1,g_2,g_3))
\end{equation}
where we set
\begin{equation*}K(\psi, \!F,\! g_1,\!g_2,\!g_3)\!:=\|F\|_{\alp,\mcl{N}^+_{r_s}}^{(1-\alp,\corners\!)}+\|g_2\|_{\alp,\Gam_{w,r_s}}^{(1-\alp,\corners\!)}
\!+\!\sum_{l=1,3}\!\|g_l\|_{\alp,\Lambda}^{(1-\alp,\der\Lambda\!)}.
\end{equation*}
So it suffices to estimate $|\psi|_{0,\alp,\mcl{N}^+_{r_s}}$.

Since $\mcl{N}^+_{r_s}$ is a cylindrical domain with the cross-section $\Lambda$ in $(r,x')-$coordinates in \eqref{spherical-coord}, there is a constant $\kappa_0>0$ depending only on $n, \Lambda$ to satisfy, for any $x_0\in \mcl{N}^+_{r_s}$ and $R>0$,
\begin{equation*}
\frac{1}{\kappa_0}\le\frac{vol(B_R(x_0)\cap\mcl{N}^+_{r_s})}{R^n}\le \kappa_0.
\end{equation*}
So, if there is a constant $R_*>0$ and $M>0$ satisfying
\begin{equation}
\label{new-at-22}
\sup_{x_0\in \mcl{N}^+_{r_s}, R>0}\frac{1}{R^{n-2+2\alp}}\int_{B_R(x_0)\cap\mcl{N}^+_{r_s}}|D\psi|^2\le M^2\;\;\text{for all}\;\;R\in(0,R_*],
\end{equation}
then there holds
\begin{equation}
\label{new-at-21}
|\psi|_{\alp,\mcl{N}^+_{r_s}}\le C(\frac{1}{R_*^{\alp}}|\psi|_{0,\mcl{N}^+_{r_s}}+M)
\end{equation}
assuming $R_*<1$ without loss of generality where the constant $C$ depends only on $n,\alp,\Lambda, r_s, r_1$. So, we devote the rest of proof to find a constant $M>0$ satisfying \eqref{new-at-22}.

To obtain \eqref{new-at-22}, we need to consider the three cases: (i) $B_R(x_0)\subset\mcl{N}^+_{r_s}$, (ii) $x_0\in \der\mcl{N}^+_{r_s}\setminus(\corners)$ and $B_R(x_0)\cap (\corners)=\emptyset$, (iii) $x_0\in \corners$. More general cases can be treated to these three cases. Also, case (i), (ii) are easier to handle than case (iii). So we only consider case (iii).

\textbf{Step 3.} Assume $x_0\in \der S_0$ because the case of $x_0\in\der\Gam_{ex}$ can be treated more simply. For a fixed point $x_0\in \der S_0$, let $x_0=(r_s,x_0')$ in the $(r,x')$-coordinates given in \eqref{spherical-coord}. Then, there is a constant $R_0>0$ depending on $\Lambda$, and a smooth diffeomorphism $h$ defined in a neighborhood $\mcl O_{x_0'}$ of $x_0'$ in $\ol{\Lambda}$ so that $h$ flattens $\der\Lambda$ near $x_0'$, and moreover the followings hold:
\begin{equation*}
\begin{split}
&(a)\;\; h(x_0')=0\in \R^{n-1},\quad\quad (b)\;\;\text{for any}\;\; x'\in \mcl{O}_{x_0'}\cap \der\Lambda,\;\; h(x')\in \R^{n-2}\times \{0\},\\
&(c)\;\;\text{for any}\;\; (r,x')\in B_{R_0}(x_0)\cap \mcl{N}^+_{r_s},\;\; (r,h(x'))\in (r_s,r_1)\times \R^{n-2}\times \R_+.
\end{split}
\end{equation*}

For a constant $R>0$ and $y_0:=(r_s,h(x_0'))=(r_s,0)\in \R^n$, let us set
\begin{equation}
\label{sets-1}
\begin{split}
&\mcl{D}_{R}:=B_{R}(y_0)\cap \{(r,h(x')):(r,x')\in B_{R_0}(x_0)\cap\mcl{N}^+_{r_s}\},\\
&\Sigma_{\Gam_w,R}:={B_{R}(y_0)}\cap \{(r,h(x')): (r,x')\in B_{R_0}(x_0)\cap \Gam_w\},\\
&\Sigma_{S_0,R}:={B_{R}(y_0)}\cap\{(r,h(x')): (r,x')\in B_{R_0}(x_0)\cap S_0\}.
\end{split}
\end{equation}

For $(r,x')\in \mcl{N}^+_{r_s}$, let us write $y=(r,y')=(r,h(x'))$, and $\phi(y)=\psi(r,x')$. By \eqref{bd1}-\eqref{bd3}, $\phi$ satisfies
\begin{equation}
\label{jewel6}
\begin{split}
&\der_k(\til a_{kl}(y,0)\der_l\phi)=\der_k\til F_k\;\;\text{in}\;\;\mcl{D}_R,\quad\der_r\phi-\mu_0\phi=\til g_1\;\;\text{on}\;\;\Sigma_{S_0,R}\\
&(\til a_{kl}(y,0)\der_l\phi)\cdot\til\nu_w=\til g_2\;\;\text{on}\;\;\Sigma_{\Gam_w,R}
\end{split}
\end{equation}
where $\til a_{kl}, \til F_k, \til g_1, \til g_2$ are obtained from $a_{kl}, F_k, g_1,g_2$ in\eqref{bd1}-\eqref{bd4} through the change of variables $(r,x')\mapsto (r,h(x'))$. So $[\til a_{kl}]_{k,l=1}^n$ is strictly positive, and the regularity of $\til a_{kl}, \til F_k, \til g_1, \til g_2$ are same as the regularity of $a_{kl}, F_k, g_1, g_2$ in $\mcl{D}_R$. Here, $\til\nu_w$ is the inward unit normal to $\Sigma_{\Gam_w,R}$. Since $h$ is smooth, and $\Lambda$ is compact, there is a constant $R_1>0$ depending on $n,\Lambda$ so that we have
\begin{equation*}
\der\mcl{D}_{2R_1}\setminus(\Sigma_{S_0,2R_1}\cup\Sigma_{\Gam_w,2R_1})\subset \der B_{2R_1}(y_0).
\end{equation*}

\textbf{Step 4.} Fix $R\in(0,R_1]$, and set $\Sigma_{\der B,R}:=\der\mcl{D}_R\cap \der B_R(y_0)$. Write $\phi$ in \eqref{jewel6} as $\phi=u+w$ for a weak solution $u$ to
\begin{equation}
\label{jewel7}
\begin{split}
&\der_k(\til a_{kl}(y_0,0)\der_lu)=0\;\;\text{in}\;\;\mcl{D}_R,\quad\;\; u=\phi\;\;\text{on}\;\;\Sigma_{\der B,R}\\
&\der_ru-\mu_0 u=0\;\;\text{on}\;\;\Sigma_{S_0,R},\quad \;\;\;\;\;(\til a_{kl}(y_0,0)\der_lu)\cdot\til{\nu}_w=0\;\;\text{on}\;\;\Sigma_{\Gam_w,R}.
\end{split}
\end{equation}
Such $u\in W^{1,2}(\mcl{D}_R)$ uniquely exists by Lemma \ref{lemma-1}, and by a basic estimate for harmonic functions as in \cite[Lemma 3.10]{Ha-L}, there is a constant $C$ depending on the data in sense of Remark \ref{remark-cons} to satisfy
\begin{equation}
\label{new-at-24}
\int_{\mcl{D}_{\varrho_1}}|Du|^2dy\le C(\frac{\varrho_1}{\varrho_2})^n\int_{\mcl{D}_{\varrho_2}}|Du|^2dy
\;\;\text{for}\;\;0<\varrho_1\le \varrho_2\le R.
\end{equation}

\textbf{Step 5.}
Then, by \eqref{jewel6} and \eqref{jewel7}, $w=\phi-u$ satisfies
\begin{equation}
\label{jewel8}
\begin{split}
&\der_k(\til a_{kl}(y_0,0)\der_lw)=\der_k([\til a_{kl}(y_0,0)-\til a_{kl}(y,0)]\der_l\phi+\til F_k)\;\;\text{in}\;\;\mcl{D}_R,\\
&\der_rw-\mu_0 w=\til g_1\;\;\text{on}\;\;\Sigma_{S_0,R},\\
&(\til a_{kl}(y_0,0)\der_l w)\cdot \til\nu_w=([\til a_{kl}(y_0,0)-\til a_{kl}(y,0)]\der_l \phi)\cdot \til\nu_w+\til g_2\;\;\text{on}\;\;\Sigma_{\Gam_w,R},\\
&w=0\;\;\text{on}\;\;\Sigma_{\der B,R}.
\end{split}
\end{equation}
By \eqref{new-40}, \eqref{d8-1} and \eqref{jewel8}, there is a constant $\lambda, C>0$ depending on the data in sense of Remark \ref{remark-cons} so that there holds
\begin{equation}
\label{new-at-27}
\lambda\! \int_{\mcl{D}_R}\!|Dw|^2 \;dy\le C(\!\int_{\mcl{D}_R}\!\!|h_1||Dw|dy\!+\!\int_{\Sigma_{S_0,R}}\!\!\!\!|h_2||w|\; dA_y\!+\!\int_{\Sigma_{\Gam_w,R}}\!\!\!\!|h_3||w|dA_y)
\end{equation}
where we set
\begin{equation}
\label{cor}
\begin{split}
&|h_1(y)|=\delta^*(y_0,y)|D\phi(y)|+|\til F(y)|\\
&|h_2(y)|=\delta^*(y_0,y)|D\phi(y)|+|\til F(y)|+|w(y)|+|\til g_1(y)|\\
&|h_3(y)|=\delta^*(y_0,y)|D\phi(y)|+|\til F(y)|+|\til g_2(y)|
\end{split}
\end{equation}
with $\delta^*(y_0,y)=(\sum_{k,l=1}^n|\til a_{kl}(y_0,0)-\til a_{kl}(y,0)|^2)^{1/2}$.

Before we proceed further, we first note that, by the definition of $\phi, \til F=(\til F_k)_{k=1}^n, \til g_1, \til g_2$ in \eqref{jewel6}, for any $\beta\in[0,\alp]$, there is a constant $C(\beta, \Lambda)>0$ depending on $\beta, \Lambda$ so that we have
\begin{equation}
\label{new-dec-1}
\begin{split}
&\|\phi\|_{1,\beta.\mcl{D}_R}^{(-\alp,\localcorners)}\le C(\beta,\Lambda)\|\psi\|_{1,\beta,\mcl{N}^+_{r_s}}^{(-\alp,\corners)}\\
&\|\til F\|_{\beta.\mcl{D}_R}^{(1-\alp,\localcorners)}\le C(\beta,\Lambda)\|F\|_{\beta,\mcl{N}^+_{r_s}}^{(1-\alp,\corners)}\;\;\text{with}\;\;
F=(F_k)_{k=1}^n\\
&\|\til g_1\|{\beta,\Sigma_{S_0,R}}^{(1-\alp,\localcorners)}\le C(\beta,\Lambda)\|g_1\|_{\beta,\Lambda}^{(1-\alp,\der\Lambda)}\\
&\|\til g_2\|_{\beta,\Sigma_{\Gam_w,R}}^{(1-\alp,\localcorners)}\le C(\beta,\Lambda)\|g_2\|_{\beta,\Gam_{w,r_s}}^{(1-\alp,\corners)}.
\end{split}
\end{equation}
Let us set $d(y):=dist(y,\localcorners)$ and
\begin{equation}
\label{Ms}
\begin{split}
&M_1(R):=R\|\psi\|_{1,0,\mcl{N}^+_{r_s}}^{(-\alp,\corners)}+\|F\|_{0,\mcl{N}^+_{r_s}}^{(1-\alp,\corners)},\\
&M_2(R):=M_1(R)+\| g_1\|_{0,\Lambda}^{(1-\alp,\der\Lambda)}
+\|g_2\|_{0,\Sigma_{\Gam_w,R}}^{(1-\alp,\corners)}.
\end{split}
\end{equation}
Using \eqref{cor}, \eqref{new-dec-1} and the continuity of $\til a_{kl}(y,0)$ in $\delta^*(y_0,y)$, we get, for any $y\in \mcl{D}_{R_1}$,
\begin{equation}
\label{hs}
\begin{split}
&|h_1(y)|\le M_1(R)[d(y)]^{-1+\alp},\\
&||h_2(y)|-|w(y)||\le M_2(R)[d(y)]^{-1+\alp},\quad |h_3(y)|\le M_2(R)[d(y)]^{-1+\alp}.
\end{split}
\end{equation}
By the Poincar\'{e} inequality with scaling, and \eqref{hs}, for any $\eps>0$, we have
\begin{equation}
\label{new-dec-7}
\begin{split}
&\int_{\Sigma_{S_0,R}}|h_2||w|dA_y\\
&\le C(\eps,n,\alp)R[M_2(R)]^2\int_{\Sigma_{S_0,R}}[d(y)]^{2(-1+\alp)}dA_y+(\eps+C(n)R)\int_{\mcl{D}_R}|Dw|^2dy.
\end{split}
\end{equation}
By \eqref{sets-1}, for any $R\le R_1$, $\Sigma_{S_0,R}$ satisfies
\begin{equation}
\label{subset}
\Sigma_{S_0,R}\subset\{(r,y_1',\cdots,y_{n-1}'):r=r_s, y'_{n-1}\in(0,R), (y'_1,\cdots,y'_{n-2})\in B^{(n-2)}_R(0)\}
\end{equation}
where $B_R^{k}(0)$ is a ball in $\R^k$ of radius $R$ centered at $0\in \R^k$ . This provides
\begin{equation}
\label{ds}
\int_{\Sigma_{S_0,R}}[d(y)]^{2(-1+\alp)}dA_y\le C(n)vol(B_R^{(n-2)}(0))\int_0^R t^{2(-1+\alp)}dt \le C(n,\alp)R^{n-3+2\alp}.
\end{equation}
Note that $\alp>\frac 12$ is an important condition in \eqref{ds}.
By \eqref{new-dec-7} and \eqref{ds}, we obtain
$$
\int_{\Sigma_{S_0,R}}|h_2||w|dA_y\le C(\eps,n,\alp)R^{n-2+2\alp}[M_2(R)]^2+(\eps+C(n)R)\int_{\mcl{D}_R}|Dw|^2dy.
$$

One can treat $\int_{\mcl{D}_R}|h_1||Dw|dy, \int_{\Sigma_{\Gam_w,R}}|h_3||w|dA_y$ similarly, and combine all the results together so that \eqref{new-at-27} implies
\begin{equation}
\label{new-dec-10}
\lambda\int_{\mcl{D}_R}\!\!|Dw|^2dy\le C(\eps,n,\alp)R^{n-2+2\alp}[M_2(R)]^2\!+\!(\eps\!+\!C(n)R)\!\int_{\mcl{D}_R}\!\!|Dw|^2 dy.
\end{equation}
 Choose $\eps=\frac{\lambda}{10}$, and reduce $R_1$ if necessary so that \eqref{new-dec-10} implies
\begin{equation}
\label{new-dec-11}
\int_{\mcl{D}_R}|Dw|^2dy\le C(n,\alp,\lambda)R^{n-2+2\alp}[M_2(R)]^2\;\;\text{for all}\;\;R\in(0,R_1].
\end{equation}
\textbf{Step 6.}
By \eqref{new-at-24} and \eqref{new-dec-11}, whenever $0<\varrho\le R\le R_1$, $\phi$ satisfies
\begin{equation*}
\int_{\mcl{D}_{\varrho}}|D\phi|^2dy
\le C\bigl((\frac{\varrho}{R})^n\int_{\mcl{D}_R}|D\phi|^2 dy+R^{n-2+2\alp}[M_2(R)]^2\bigr)
\end{equation*}
So, by \cite[Lemma 3.4]{Ha-L} and \eqref{new-dec-1}, we obtain
\begin{equation}
\label{new-dec-12}
\begin{split}
\int_{\mcl{D}_R}|D\phi|^2dy
&\le CR^{n-2+2\alp}(\frac{R^{2-2\alp}}{R_1^n}\|D\phi\|^2_{L^2(\mcl{D}_{R})}+[M_2(R)]^2)\\
&\le CR^{n-2+2\alp}\bigl((\|\psi\|_{1,0,\mcl{N}^+_{r_s}}^{(-\alp,\corners)})^2+[M_2(R_1)]^2\bigr).
\end{split}
\end{equation}
We note that the estimate \eqref{new-dec-12} holds true for any choice of $x_0=h^{-1}(y_0)\in \mcl{N}^+_{r_s}$ and any $R\in(0,R_1]$, moreover $R_1$ is uniform for all $x_0\in \mcl{N}^+_{r_s}$. Since $\Lambda$ is compact, and $\der\Lambda$ is smooth, we can take a constant $R^*>0$ depending on $\Lambda$ so that, by \eqref{new-dec-12}, the solution $\psi$ to \eqref{bd1}-\eqref{bd4} satisfies
\begin{equation}
\label{new-dec-13}
\int_{B_R(x_0)\cap \mcl{N}^+_{r_s}}|D\psi|^2
\le CR^{n-2+2\alp}(\|\psi\|_{1,0,\mcl{N}^+_{r_s}}^{(-\alp,\corners\!)}\!+K_0)^2
\end{equation}
with $K_0=\|F\|_{0,\mcl{N}^+_{r_s}}^{(1-\alp,\corners\!)}
\!+\!\!\sum_{l=1,3}\!\|g_l\|_{0,\Lambda}^{(1-\alp,\der\Lambda\!)}\!+\!\|g_2\|_{0,\Gam_{w,r_s}}^{(1-\alp,\corners\!)}$.

By the interpolation inequality in \cite[Lemma 6.34]{GilbargTrudinger}, for any $\eps>0$, we have
$$
\|\psi\|_{1,0,\mcl{N}^+_{r_s}}^{(-\alp,\corners)}\le C(\eps,\alp)|\psi|_{0,\mcl{N}^+_{r_s}}+\eps\|\psi\|_{1,\alp,\mcl{N}^+_{r_s}}^{(-\alp,\corners)}.
$$
So, assuming $R^*<1$ without loss of generality, by \eqref{new-dec-16}, \eqref{new-at-21}, we obtain
\begin{equation*}
\|\psi\|_{1,\alp,\mcl{N}^+_{r_s}}^{(-\alp,\corners)}\le C(|\psi|_{0,\mcl{N}^+_{r_s}}+K(F,g_1,g_2,g_3))
\end{equation*}
where $K(F,g_1,g_2,g_3)$ is defined after \eqref{new-dec-16}.

Since $\mu_0>0$ in \eqref{bd2} by Lemma \ref{lemma-1}, by the uniqueness of solution(\cite[Corollary 2.5]{Lie-1}) for \eqref{bd1}-\eqref{bd4}, we have
$
|\psi|_{0,\mcl{N}^+_{r_s}}\le CK(F,g_1,g_2,g_3).
$
The proof is complete.
\end{proof}

\bigskip
{\bf Acknowledgments.} The authors thank Marshall Slemrod for helpful comments on this work. Mikhail Feldman's research was supported in part by the National Science Foundation under Grants
DMS-0800245 and DMS-0354729.

\end{document}